\documentclass[10pt,a4paper]{CMbook}
\usepackage{ifpdf}
\ifpdf%
\PassOptionsToPackage{pdftex}{graphicx}
\PassOptionsToPackage{pdftex}{color}
\else%
\PassOptionsToPackage{dvips}{graphicx}%
\fi

\usepackage{graphicx,color}
\usepackage{xr}
\usepackage{eurosym}
\usepackage{%amsmath,
amsfonts,amssymb}
\usepackage{mathrsfs}
\usepackage{oldstyle}
\usepackage[T1]{fontenc}
\usepackage[utf8]{inputenc}
\usepackage[francais]{babel}
\usepackage[all]{xy}
\usepackage{pb-diagram}
\usepackage[french,nohints]{minitoc}
\usepackage{multicol}
\usepackage{CMmacros2}
\usepackage{makeidx}\makeindex\makeglossary
\usepackage{lmodern}
\usepackage{microtype}
\usepackage{array}
\usepackage{stmaryrd}
\usepackage{float}
\usepackage{wasysym}

\usepackage{bm}

\DeclareFontFamily{U}{wasy}{}
\DeclareFontShape{U}{wasy}{m}{n}{<5> <6> <7> <8> <9> gen * wasy
      <10> <10.95> <12> <14.4> <17.28> <20.74> <24.88>wasy10  }{}
\DeclareFontShape{U}{wasy}{b}{n}{ <-10> ssub * wasy/m/n
 <10> <10.95> <12> <14.4> <17.28> <20.74> <24.88>wasyb10 }{}
\DeclareFontShape{U}{wasy}{bx}{n}{<5> <6> <7> <8> <9> gen * wasy
 <10> <10.95> <12> <14.4> <17.28> <20.74> <24.88>wasyb10}{}
\DeclareSymbolFont{wasy}{U}{wasy}{m}{n}
\SetSymbolFont{wasy}{bold}{U}{wasy}{bx}{n}
\DeclareFontFamily{U}{lasy}{}
\DeclareFontShape{U}{lasy}{m}{n}{ <5> <6> <7> <8> <9> gen * lasy
      <10> <10.95> <12> <14.4> <17.28> <20.74> <24.88>lasy10  }{}
\DeclareFontShape{U}{lasy}{b}{n}{ <5> <6> <7> <8> <9> gen * lasy
      <10> <10.95> <12> <14.4> <17.28> <20.74> <24.88>lasyb10  }{}
\DeclareFontFamily{U}{stmry}{}
\DeclareFontShape{U}{stmry}{m}{n}
   {  <5> <6> <7> <8> <9> <10> gen * stmary
      <10.95><12><14.4><17.28><20.74><24.88>stmary10%
   }{}
\DeclareFontShape{U}{stmry}{b}{n}
   {  <5> <6> <7> <8> <9> <10> gen * stmary
      <10.95><12><14.4><17.28><20.74><24.88>stmary10%
   }{}

\ifpdf
\makeatletter
\def\toclevel@chapternb{0}
\makeatother
\usepackage[bookmarksopen=false,pdftex=true,breaklinks=true,backref=page,%
    pagebackref=true,plainpages=false,%
    hyperindex=true,pdfstartview=FitH,colorlinks=true,
    pdfpagelabels=true,linkcolor=blue,citecolor=red]%
    {hyperref}

\else

\fi

\DeclareSymbolFont{lasy}{U}{lasy}{m}{n}
\SetSymbolFont{lasy}{bold}{U}{lasy}{b}{n}
\let\Box\undefined
\DeclareMathSymbol\Box{\mathord}{lasy}{"32}

\title{Calculs matriciels sur les domaines de Prüfer}
\author{Gema-Maria D\'\i az-Toca, Henri Lombardi}
\date{\today }
\DeclareMathAlphabet{\mathpzc}{OT1}{pzc}{m}{it}

\renewcommand\chaptername{}

\begin{document}

\let\oldref\ref
\renewcommand{\ref}[1]{\hbox{\oldref{#1}}}

\makeatletter
\def\U@@@ref#1\relax#2\relax{\let\showsection\relax#2}
\newcommand{\U@@ref}[5]{\U@@@ref#1}
\newcommand{\U@ref}[1]{\NR@setref{#1}\U@@ref{#1}}
\newcommand{\@iref}{\@ifstar\@refstar\U@ref}
\newcommand{\iref}[1]{\hbox{\@iref{#1}}}
\def\l@chapter#1#2{\ifnum\c@tocdepth>\m@ne\addpenalty{-\@highpenalty}\vskip1.0em\@plus\p@
\setlength\@tempdima{2.8em}\begingroup\parindent\z@\rightskip\@pnumwidth\parfillskip-\@pnumwidth\penalty-2000\leavevmode\bfseries
\advance\leftskip\@tempdima\hskip-\leftskip{\boldmath#1}\nobreak\hfil\nobreak\hb@xt@\@pnumwidth{\hss}\par\penalty\@highpenalty\endgroup\fi}
\def\l@chapterbis#1#2{\ifnum\c@tocdepth>\m@ne\addpenalty{-\@highpenalty}\vskip1.0em\@plus\p@
\setlength\@tempdima{2em}\begingroup\parindent\z@\rightskip\@pnumwidth\parfillskip-\@pnumwidth\penalty-2000\leavevmode\bfseries
\advance\leftskip\@tempdima\hskip-\leftskip{\boldmath#1}\nobreak\hfil\nobreak\hb@xt@\@pnumwidth{\hss#2}\par\penalty\@highpenalty\endgroup\fi}
\def\toclevel@chapterbis{0}
\def\l@section{\@dottedtocline{1}{2.8em}{1.2em}}
\def\l@subsection{\@dottedtocline{2}{4.7em}{1em}}
\def\l@subsubsection{\@dottedtocline{3}{6.5em}{1em}}
\def\contentsline#1#2#3#4{\ifx\\#4\\\csname
l@#1\endcsname{#2}{#3}\else\ifHy@linktocpage\csname
l@#1\endcsname{{#2}}{\hyper@linkstart{link}{#4}{#3}\hyper@linkend }
\else\csname l@#1\endcsname{\hyper@linkstart{link}{#4}{#2}\hyper@linkend
}{#3}\fi\fi}
\makeatother

%%%%%%%%%%%%%%%%%%%%%%%%%%%%%%
%!TEX root =  CohenPrufer.tex

\date{2012}

\ifx\mesMacrosDejaChargees\undefined\let\chargeMesMacros\relax
\else\let\chargeMesMacros \fi
\chargeMesMacros

\gdef\mesMacrosDejaChargees{}

\makeatletter
\def\smalll{\small\smallbigs}
\def\smallbigs{\def\big##1{{\raise-0,25pt\hbox{\Large$\mathsurround-0,5pt\scriptstyle##1$}}}%
\def\Big##1{{\hbox{$\left##1\vbox to7.5\p@{}\right.\n@space$}}}%
\def\bigg##1{{\hbox{$\left##1\vbox to9.5\p@{}\right.\n@space$}}}%
\def\Bigg##1{{\hbox{$\left##1\vbox to12.5\p@{}\right.\n@space$}}}}
\makeatother
%%%%%%%%%%%%%%%%smalltimes
\def\smalltimes{\raise0.75pt\hbox{\boldmath$\mathsurround0pt\scriptstyle\times$}}
\def\ssmalltimes{\raise0.5pt\hbox{$\mathsurround0pt\scriptscriptstyle\times$}}
%%%%%%%%%%%%%%%%%

\FrenchFootnotes

%%%%%%%%%%%%%%%%%%%%%%%%%%%%%%%%%%%%%%%%%%%%%%%%%%%%%%%%%%%%%

\renewcommand{\labelenumii}{\textit{\theenumii.}}
\renewcommand{\labelenumi}{\textit{\theenumi.}}

%\def\emptyline{\vspace{12pt}}
%\DefineParaStyle{Maple Output}
%\DefineParaStyle{Maple Output}
%\DefineParaStyle{Warning}
%\DefineCharStyle{2D Math}
%\DefineCharStyle{2D Output}
%\DefineCharStyle{Maple 2D Input}
%\DefineCharStyle{Maple Maple Input}
%\DefineCharStyle{Maple 2D Output}
%\DefineCharStyle{Maple 2D Math}

%:  index
\newcommand \ix[1] {\index{#1}\emdf{#1}}
\newcommand \ixc[2] {\index{#1!#2}\emdf{#1}}
\newcommand \ixx[2] {\index{#1!#2}\emdf{#1 #2}}
\newcommand \ixy[2] {\index{#2!#1}\emdf{#1 #2}}
\newcommand \ixd[2] {\index{#2!#1}\emdf{#1}}
\newcommand \ixe[2] {\index{#2@#1}\emdf{#1}}
\newcommand \ixf[3] {\index{#3@#2!#1}\emdf{#1}}
\newcommand \ixg[3] {\index{#3@#1!#2}\emdf{#1}}

\newcommand \iplg {\index{principe local-global!de base} }
\newcommand \iplgz {\index{principe local-global!de base}}

%: lecteur, lectrice, masculin feminin

\newcounter{MF}

\newcommand\stMF{\stepcounter{MF}}

\newcommand{\lec}{\stMF\ifodd\value{MF}lecteur \else
lectrice \fi}
\newcommand{\lecz}{\stMF\ifodd\value{MF}lecteur\else lectrice\fi}

\newcommand{\lecs}{\stMF\ifodd\value{MF}lecteurs \else
lectrices \fi}
\newcommand{\lecsz}{\stMF\ifodd\value{MF}lecteurs\else
lectrices\fi}

\newcommand{\alec}{\stMF\ifodd\value{MF}au lecteur \else%
\`a la lectrice \fi}
\newcommand{\alecz}{\stMF\ifodd\value{MF}au lecteur\else%
\`a la lectrice\fi}

\newcommand{\dlec}{\stMF\ifodd\value{MF}du lecteur \else%
de la lectrice \fi}
\newcommand{\dlecz}{\stMF\ifodd\value{MF}du lecteur\else%
de la lectrice\fi}

\newcommand{\llec}{\stMF\ifodd\value{MF}le lecteur \else la lectrice \fi}
\newcommand{\llecz}{\stMF\ifodd\value{MF}le lecteur\else la lectrice\fi}

\newcommand{\Llec}{\stMF\ifodd\value{MF}Le lecteur \else La lectrice \fi}
\newcommand{\Llecz}{\stMF\ifodd\value{MF}Le lecteur\else La lectrice\fi}

\newcommand{\lui}{\ifodd\value{MF}lui \else
elle \fi}
\newcommand{\luiz}{\ifodd\value{MF}lui\else
elle\fi}

\newcommand{\f}{\ifodd\value{MF}f \else
ve \fi}

\newcommand{\x}{\ifodd\value{MF}x \else
se \fi}

\newcommand{\il}{\ifodd\value{MF}il \else
elle \fi}

\newcommand{\Il}{\ifodd\value{MF}Il \else
Elle \fi}

\newcommand{\eux}{\ifodd\value{MF}eux \else elles \fi}
\newcommand{\euxz}{\ifodd\value{MF}eux\else elles\fi}

\makeatletter
\newcommand{\la}{\@ifstar{\ifodd\value{MF}la\else
le\fi}{\stMF\ifodd\value{MF}la \else le \fi}}
\makeatother
%  \la fait le ou la sans changer le genre,
%  \la* fait le ou la en changeant le genre

\newcommand{\e}{\ifodd\value{MF} \else e \fi}
\newcommand{\ez}{\ifodd\value{MF}\else e\fi}
\newcommand{\es}{\ifodd\value{MF}s \else es \fi}
\newcommand{\esz}{\ifodd\value{MF}s\else es\fi}
%%%%%%%%%%%%%%%%%%%%%%%%%%%%%%%%%%%%%%%%%%%%%%%%%%%%%%%%%%%%%%%%%%%%%%%%%%%

%:Macros a parametres
\newcounter{solex}[chapter]

\newcommand\tbf [1]{\texttt{\textbf{#1}}}
\newcommand\tsbf [1]{\textsf{\textbf{#1}}}

\newcommand\emdf[1]{\textbf{\textit{#1}}}
\newcommand\emdi[1]{\emdf{#1}\index{#1}}

\newcommand \vab[2]{[\,#1\;#2\,]}

\newcommand \fraC[2] {{{#1}\over {#2}}}
\newcommand \tra[1] {{\,^{\rm t}\!#1}}
\newcommand \tralst[1] {\tra\,{\lst{#1}}}
\newcommand \lra[1] {\langle{#1}\rangle}
\newcommand \lrA[1] {\big\langle{#1}\big\rangle}
\newcommand \gen[1] {\left\langle{#1}\right\rangle}
\newcommand \geN[1] {\big\langle{#1}\big\rangle}
\newcommand \Zz[1] {{\ZZ/{#1}\ZZ}}
\newcommand \Z[1]  {{\ZZ/{#1}\ZZ}}

\newcommand \uci[1] {{\buildrel{\circ}\over{#1}}}

\newcommand\sur[1]{\!\left/#1\right.}
\newcommand\suR[1]{\big/#1 }
\renewcommand\matrix[1]{{\begin{array}{ccccccccccccccccccccccccc} #1 \end{array}}}  

\newcommand\abs[1]{\left|{#1}\right|}
\newcommand\abS[1]{\big|{#1}\big|}
\newcommand\vala[1]{\abs{#1}}
\newcommand\aqo[2]{#1\sur{\gen{#2}}\!}
\newcommand\aqO[2]{#1\suR{\geN{#2}}\!}
\newcommand\bloc[4]{\left[\matrix{#1 & #2 \cr #3 & #4}\right]}

\newcommand{\Dpp}[2]{{{\partial #1}\over{\partial #2}}}

\newcommand\carray[2]{{\left[\begin{array}{#1} #2 \end{array}\right]}}
\newcommand\cmatrix[1]{\left[\matrix{#1}\right]}
\newcommand\Cmatrix[2]{\setlength{\arraycolsep}{#1}\left[\matrix{#2}\right]}
\newcommand\cmatriX[1]{\Cmatrix{2pt}{#1}}
\newcommand\clmatrix[1]{{\left[\begin{array}{lllllll} #1 \end{array}\right]}}
\newcommand\dmatrix[1]{\abs{\matrix{#1}}}

\newcommand\mat[2]{\MM_{#2}(#1)}

\newcommand\so[1]{\left\{{#1}\right\}}
\newcommand\sO[1]{\big\{{#1}\big\}}
\newcommand\sotq[2]{\so{\,#1\mid#2\,}}
\newcommand\sotQ[2]{\sO{\,#1\mid#2\,}}
\newcommand\frt[1]{\!\left|_{#1}\right.\!}
\newcommand\sims[1]{\buildrel{#1}\over \sim}
\newcommand\norme[1]{\abs{\abs{#1}}}
\newcommand\formule[1]{{\left\{ \begin{array}{lll} #1 \end{array}\right.}}
\newcommand\formul[2]{{\left\{ \begin{array}{#1} #2 \end{array}\right.}}
\newcommand{\Vec}[1]{{\overrightarrow{#1}}}

\newcommand \bu[1] {{{#1}\bul}}
\newcommand \ci[1] {{{#1}^\circ}}
\newcommand \wi[1] {\widetilde{#1}}
\newcommand \wh[1]{{\widehat{#1}}}
\newcommand \ov[1] {\overline{#1}}
\newcommand \und[1] {\underline{#1}}

% indice d'un sous groupe
\newcommand \idg[1] {|\,#1\,|}
\newcommand \idG[1] {\big|\,#1\,\big|}

% dimension d'une extension
\newcommand \dex[1] {[\,#1\,]}
\newcommand \deX[1] {\big[\,#1\,\big]}

% liste
\newcommand \lst[1] {[\,#1\,]}
\newcommand \lsT[1] {\big[\,#1\,\big]}

\newcommand \thref[1] {\thoz~\ref{#1}}
\newcommand \thrf[1] {\thoz~\ref{#1}}
\newcommand \paref[1] {page~\pageref{#1}}
\newcommand{\pref}[1]{\textup{\hbox{\normalfont(\iref{#1})}}}
\newcommand \egrf[1] {\egtz~\pref{#1}}
\newcommand \egref[1] {\egtz~\pref{#1}  \paref{#1}}
\newcommand \eqrf[1] {\'equation~\pref{#1}}
\newcommand \eqref[1] {\'equation~\pref{#1} \paref{#1}}

\newcommand\lrb[1] {\llbracket #1 \rrbracket}
\newcommand\lrbn {\lrb{1..n}}
\newcommand\lrbl {\lrb{1..\ell}}
\newcommand\lrbm {\lrb{1..m}}
\newcommand\lrbk {\lrb{1..k}}
\newcommand\lrbh {\lrb{1..h}}
\newcommand\lrbp {\lrb{1..p}}
\newcommand\lrbq {\lrb{1..q}}
\newcommand\lrbr {\lrb{1..r}}
\newcommand\lrbs {\lrb{1..s}}

\newcommand \DCOL[1] {\begin{multicols}{2}\smalll #1 \end{multicols}}
\newcommand \TCOL[1] {\begin{multicols}{3}\footnotesize #1 \end{multicols}}

\newcommand\oge{\leavevmode\raise.3ex\hbox{$\scriptscriptstyle\langle\!\langle\,$}}
\newcommand\feg{\leavevmode\raise.3ex\hbox{$\scriptscriptstyle\,\rangle\!\rangle$}}

\DeclareRobustCommand{\guig}{\mbox{{\usefont{U}{lasy}%
{\if b\expandafter\@car\f@series\@nil b\else m\fi}{n}%
\char40\kern-0.20em\char40}~}}
\DeclareRobustCommand{\guid}{\mbox{~\usefont{U}{lasy}%
{\if b\expandafter\@car\f@series\@nil b\else m\fi}{n}%
\char41\kern-0.20em\char41}}
\newcommand\gui[1]{\oge{#1}\feg}

%:  CAdre
\newcommand \CAdre[2]{%
\begin{center}
\begin{tabular}%
{|p{#1\textwidth}|}
\hline
\vspace{-1.5mm}
#2 
\vspace{1mm}\\ %
\hline
\end{tabular}
\end{center} 
}

%:  Grandcadre
\newcommand \Grandcadre[1]{%
\begin{center}
\begin{tabular}{|c|}
\hline
~\\[-3mm]
#1\\[-3mm]
~\\
\hline
\end{tabular}
\end{center}

}

\newcommand\boite[2]{\begin{minipage}[c]{#1\cm}
     \centering {#2} \end{minipage}}    
\newcommand\Boite[3]{\parbox[t][#1\cm][c]{#2\cm}{\boite{#2}{#3}}}

%%%%%%%%%%%%%%%%%%% compteur bidon %%%%%%%
\newcounter{bidon}
\newcommand{\rdb}{\refstepcounter{bidon}}

%: newtheorem 
\CMnewtheorem{theorem}{Th\'eor\`eme}{\itshape}   
\CMnewtheorem{thdef}{Th\'eor\`eme et d\'efinition}{\itshape}
\CMnewtheorem{plcc}{Principe local-global concret}{\itshape}
\CMnewtheorem{plca}{Principe local-global abstrait$\mathbf{^*}$}{\itshape}
\CMnewtheorem{proposition}{Proposition}{\itshape}
%\CMnewtheorem{prop}{Proposition}{\itshape}
\CMnewtheorem{propdef}{Proposition et d\'efinition}{\itshape}
\CMnewtheorem{lemma}{Lemme}{\itshape}
\CMnewtheorem{algori}{Algorithme}{\sfshape}
\CMnewtheorem{lem}{Lemme}{\itshape}
\CMnewtheorem{cor}{Corollaire}{\itshape}
\CMnewtheorem{corollary}{Corollaire}{\itshape}
\CMnewtheorem{fact}{Fait}{\itshape}
\CMnewtheorem{theoremc}{Th\'{e}or\`{e}me\etoz}{\itshape}
\CMnewtheorem{lemmac}{Lemme\etoz}{\itshape}
\CMnewtheorem{corollaryc}{Corollaire\etoz}{\itshape}
\CMnewtheorem{proprietec}{Propri\'{e}t\'{e}\etoz}{\itshape}
\CMnewtheorem{propositionc}{Proposition\etoz}{\itshape}
\CMnewtheorem{factc}{Fait\etoz}{\itshape}
\CMnewtheorem{lemnoyaux}{Lemme des noyaux}{\itshape}

\CMnewtheorem{remark}{Remarque}{}
\CMnewtheorem{remarks}{Remarques}{}
\CMnewtheorem{comment}{Commentaire}{}
\CMnewtheorem{comments}{Commentaires}{}
\CMnewtheorem{example}{Exemple}{}
\CMnewtheorem{examples}{Exemples}{}
\CMnewtheorem{defic}{D\'efinition\etoz}{}
\CMnewtheorem{definition}{D\'efinition}{}
\CMnewtheorem{defin}{D\'efinition}{}
\CMnewtheorem{definitions}{D\'efinitions}{}
\CMnewtheorem{definota}{D\'efinition et notation}{}
\CMnewtheorem{definotas}{D\'efinitions et notations}{}
\CMnewtheorem{convention}{Convention}{}
\CMnewtheorem{notation}{Notation}{} 
\CMnewtheorem{notations}{Notations}{} 
\CMnewtheorem{question}{Question}{}
%\CMnewtheorem{exercise}{Exercice}{\itshape}
%\CMnewtheorem{problem}{Probl\`eme}{}
\CMnewtheorem{algorithm}{Algorithme}{}

\newcounter{exercise}[chapter]

\newenvironment{exo}{\ifhmode\par\fi
\vskip-\lastskip\vskip1.5ex\mou\penalty-300 \relax
\everypar{}\noindent
\refstepcounter{exercise}{\bfseries Exercice \theexercise.}\relax
%\itshape
\ignorespaces}{\par\vskip-\lastskip\vskip1em}

\newenvironment{Exo}{\ifhmode\par\fi
\vskip-\lastskip\vskip1.5ex\mou\penalty-300 \relax
\everypar{}\noindent
\refstepcounter{exercise}{\bfseries Exercice$^\natural$\hspace{.4mm}\theexercise.}\relax
%\itshape
\ignorespaces}{\par\vskip-\lastskip\vskip1em}
\newcounter{problem}[chapter]
\newenvironment{problem}{\ifhmode\par\fi
\vskip-\lastskip\vskip1em\mou\penalty-300 \relax
\everypar{}\noindent
\refstepcounter{problem}{\bfseries Probl\`eme \theproblem.}\relax
\itshape\ignorespaces}{\par\vskip-\lastskip\vskip1em}
 
%: sectioN pour les sols d'exos
\newcommand\sectioN[1]{\stepcounter{solex}
\newpage\rdb\pagestyle{CMExercicesheadings}
\addcontentsline{toc}{section}{Solutions du chapitre \ref{#1}}
{\LARGE\bfseries\raggedright\setlength{\leftskip}{\widthof{\LARGE\bfseries\thesolex. }}\hspace*{-\leftskip}\thesolex. \nameref{#1}\par} \vspace{3mm}
\markboth{Solutions des exercices}{Solutions du chapitre \ref{#1}.}
\smalll}

%  pour les sols d'exo de l'annexe
\newcommand\sectiON[1]{\stepcounter{solex}
\newpage\rdb\pagestyle{CMExercicesheadings}
{\LARGE\bfseries\raggedright\setlength{\leftskip}
{\widthof{\LARGE\bfseries\ref{#1}.%\thesolex. 
               }}
\hspace*{-\leftskip}\ref{#1}. %\thesolex. 
\nameref{#1}\par} \vspace{3mm}
\markboth{Solutions des exercices}{Solutions de l'annexe \ref{#1}.}
\addcontentsline{toc}{section}{Solutions de l'annexe \ref{#1}}
\smalll}

\newtheorem{sugg}{Exercice possible}[section]

\newcommand{\eop}{\hfill \mbox{$\Box$}}

%:  SiModules, Intro, Exercices
\newcommand{\SiModules}[1]{#1}

\newcommand\entrenous[1] {\sibrouillon{
{\begin{flushleft}

\textsf{\textbf{Entre nous\,:}}
{\smalll\sf  #1\par}
\vspace{-.8mm}
\textsf{\textbf{Fin d'entre nous}}\end{flushleft}}}}

\newcommand \intro{\addcontentsline{toc}{section}{Introduction}
\subsection*{Introduction}}

\newcommand \Intro{\addcontentsline{toc}{section}{Introduction}\markright{Introduction}%
\pagestyle{CMExercicesheadings}\subsection*{Introduction} }

\newcommand{\bonbreak}{\penalty-1500}

\newcommand\Exercices{\bonbreak\rdb
\subsection*{Exercices}
\addcontentsline{toc}{subsection}{Exercices}

}

%:              modif M1
\newcommand\exer[1]{{\medskip \bonbreak\noindent\rdb \textbf{Exercice \ref{#1}}, \paref{#1}.}
\\
}

\newcommand\exeR[1]{{\medskip \bonbreak\noindent\rdb \textbf{Exercice \ref{#1}}, \paref{#1}.}
}

\newcommand\prob[1]{{\medskip \bonbreak\noindent\rdb \textbf{Probl\`eme \ref{#1}} \paref{#1}. }}

\newenvironment{proof}{\ifhmode\par\fi\vskip-\lastskip\vskip0.5ex\global\insidedemotrue
\everypar{}\noindent{\setbox0=\hbox{\emph{D\'emonstration.} %$\!\!$\DebP 
}\global\wdTitreEnvir\wd0\box0}
%: D\'emonstration.  a la place de \demoname  espace apres plus petit 
\ignorespaces}{\enddemobox\par\vskip.5em}
\def\enddemobox{\ifinsidedemo
\ifmmode\hbox{$\square$}\else
\ifhmode\unskip\else\noindent\fi\nobreak\null\nobreak\hfill
\nobreak$\square$\fi\fi\global\insidedemofalse}

\newenvironment{preuve}{
\trivlist \item[\hskip \labelsep{\it D\'emonstration.}]}{\hfill\mbox{$\Box$}
\endtrivlist}

\newenvironment{Proof}[1]{
\trivlist \item[\hskip \labelsep{\it #1.}]}{\hfill\mbox{$\Box$}
\endtrivlist}

\newcommand\facile{
\emph{D\'emonstration,} laiss\'ee \alecz.\hfill\mbox{$\Box$}

\vspace{.5em}
}

%\newcommand\facile{\begin{proof}
%La \dem est laiss\'ee \alecz.
%\end{proof}
%}

%:  Algorithmes 

\floatstyle{boxed}
\floatname{agc}{Circuit}
\newfloat{agc}{ht}{lag}[section]
\floatname{agC}{Circuit}
\newfloat{agC}{H}{lag}[section]

\newenvironment{algor}[1][]
{\par\smallskip\begin{agc}
\vskip 1mm
\begin{algorithm}{\bfseries#1}
\upshape\sffamily
}
{\end{algorithm}
\end{agc}
}

\newenvironment{algorH}[1][]
{\par\smallskip\begin{agC}\vskip 1mm
\begin{algorithm}{\bfseries#1}
\upshape\sffamily
}
{\end{algorithm}
\end{agC}
}

\newcommand\Vrai{\mathsf{Vrai}}
\newcommand\Faux{\mathsf{Faux}}
\newcommand\ET{\;\mathsf{ et }\;}
\newcommand\OU{\;\mathsf{ ou }\;}
\newcommand\aff{\,\leftarrow\,}
\def\pour#1#2#3{\textbf{Pour } $#1$ \textbf{ de } $#2$
     \textbf{ \`a } $#3$ \textbf{ faire }}
\def\por#1#2#3{\textbf{Pour } $#1$ \textbf{ de } $#2$
     \textbf{ \`a } $#3$  }
\def\sialors#1{\textbf{Si } $#1$ \textbf{ alors }}
\def\tantque#1{\textbf{Tant que } $#1$ \textbf{ faire }}
\newcommand\finpour{\textbf{fin pour}}
\newcommand\sinon{\textbf{sinon }}
\newcommand\sinonsi[1]{\textbf{sinon si } $#1$ \textbf{ alors }}
\newcommand\finsi{\textbf{fin si }}
\newcommand\fintantque{\textbf{fin tant que }}
\newcommand\Debut{\\[1mm] \textbf{D\'ebut }}
\newcommand\debut{\textbf{D\'ebut }}
\newcommand\Fin{\\ \textbf{Fin.}}
\newcommand\Entree{\\ \textbf{Entr\'ee : }}
\newcommand\Sortie{\\ \textbf{Sortie : }}
\newcommand\Varloc{\\ \textbf{Variables locales : }}
\newcommand\Repeter{\textbf{R\'ep\'eter }}
\newcommand\Retourner{\textbf{Retourner }}
\newcommand\jusqua{\textbf{jusqu'\`a ce que }}
\newcommand\hsz{\\ }
\newcommand\hsu{\\ \hspace*{4mm}}
\newcommand\hsd{\\ \hspace*{8mm}}
\newcommand\hst{\\ \hspace*{1,2cm}}
\newcommand\hsq{\\ \hspace*{1,6cm}}
\newcommand\hsc{\\ \hspace*{2cm}}
\newcommand\hsix{\\ \hspace*{2,4cm}}
\newcommand\hsept{\\ \hspace*{2,8cm}}
\newcommand\hsud[1]{\hsu$#1$}
\newcommand\Pro[3]{\\ \textbf{profondeur} $#1$ : #2 \hsud{#3}}
\newcommand\Eta[3]{\\ \textbf{\'Etape} $#1$ : #2 \hsud{#3}}
\newcommand\Etas[3]{\\ \textbf{\'Etapes} $#1$ : #2 \hsud{#3}}
\newcommand\Etap[3]{\\ \textbf{\'Etape} $#1$ :  #2 \hsu {#3}  }
\newcommand \Boucle {\textbf{Boucle}}
\newcommand \sortirboucle {\textbf{sortir de la boucle}}
\newcommand \finboucle {\textbf{fin de boucle}}

%:  en exposant ou en indice
\newcommand \bul{^\bullet}
\newcommand \eci{^\circ}
\newcommand \esh{^\sharp}
\newcommand \efl{^\flat}
\newcommand \epr{^\perp}
\newcommand \eti{^\times}
\newcommand \etl{^* }
\newcommand \eto{$^*\!$ }
\newcommand \etoz{$^*\!$}
\newcommand \sta{^\star}
\newcommand \dual{\uci{\phantom .}}
\newcommand \ista{_\star}

\newcommand\Abul {\gA\!\bul}
\newcommand\Ati {\gA^{\!\times}}
\newcommand\Asta {\gA^{\!\star}}
\newcommand\Atl {\gA^{\!*}}
\newcommand\Zbul {\gZ\bul}
\newcommand\Zti {\gZ^{\times}}
\newcommand\Zsta {\gZ^{\star}}
\newcommand\Ztl {\gZ^{*}}

\newcommand \iBA {_{\gB/\!\gA}}
\newcommand \iBK {_{\gB/\gK}}

%:  treillis distributifs, divisibilit\'e

\newcommand \divi {\mid}
\def \nedivi {\not\kern 2.5pt\mid}

\newcommand \vu {\vee} % sup dans les treillis
\newcommand \vi {\wedge} % inf dans les treillis
\newcommand \Vu {\bigvee\nolimits}
\newcommand \Vi {\bigwedge\nolimits}

%%%%%%%%%%%%%%%%%%%%%%%%%%%%%%%%%%%%%%%%%%%%%%%%%%%%%%%%%%%%%%%%%%
%:  du genre :  subsection* avec table des mati\`eres

%%%%%%%%%%%%%%%%%%%%%%%%%%%%%%%%%%%%%%%%%
\newcommand\subsec[1]{\bonbreak\rdb

 \subsection*{#1}
\addcontentsline{toc}{subsection}{#1}}
%%%%%%%%%%%%%%%%%%%%%%%%%%%%%%%%%%%%%%%%%

%%%%%%%%%%%%%%%%%%%%%%%%%%%%%%%%%%%%%%%%%
\newcommand\subsect[2]{\bonbreak\rdb

 \subsection*{#1}
\addcontentsline{toc}{subsection}{#2}}
%%%%%%%%%%%%%%%%%%%%%%%%%%%%%%%%%%%%%%%%%

%%%%%%%%%%%%%%%%%%%%%%%%%%%%%%%%%%%%%%%%%
\newcommand\subsubsec[1]{\bonbreak\rdb

 \subsubsection*{#1}
\addcontentsline{toc}{subsubsection}{#1}}
%%%%%%%%%%%%%%%%%%%%%%%%%%%%%%%%%%%%%%%%%

%%%%%%%%%%%%%%%%%%%%%%%%%%%%%%%%%%%%%%%%%
\newcommand\subsubsect[2]{\bonbreak\rdb

 \subsubsection*{#1}
\addcontentsline{toc}{subsubsection}{#2}}

%%%%%%%%%%%%%%%%%%%%%%%%%%%%%%%%%%%%%%%%%%%%%%%%%%%%%%%%%%%%%%%%%%%%%%%%%%%

%\newcommand\snic[1] {\vspace{.1cm}\noindent\centerline{$#1$}\vspace{.1cm}}

\newcommand \snic[1]{

{\centering$#1$\par}

}
\newcommand \snac[1]{

{\smalll\centering$#1$\par}

}
\newcommand \snuc[1]{

{\footnotesize\centering$#1$\par}
}

\newcommand \cli[1] {\left[{#1}\right]}

%:  Ahbon
\newdimen\oldleftskip
\newdimen\oldrightskip

\newcommand{\Ahbon}[1]{%
\medskip%
\oldleftskip=\leftskip%
\oldrightskip=\rightskip%
\leftskip=-\tabcolsep%
\rightskip=-\tabcolsep%
\begin{center}\begin{tabular}%
{|p{.9\textwidth}}
\smalll #1%
\end{tabular}\end{center}\par\medskip%
\leftskip=\oldleftskip%
\rightskip=\oldrightskip}

%   \adots
%   trois petits points en diagonale (axe so-ne) pour les matrices
%   (\ddots de plain: idem sur l'axe no-se)

\newcommand\adots{\mathinner{\mkern0mu\raise 1pt\hbox{\string.}\mkern
3mu\raise 4pt\hbox{\string.}\mkern 3mu\raise 7pt\hbox{\string.}}}

%: espacements
\newcommand\sms{\smallskip}
\newcommand\ms{\medskip}
\newcommand\bs{\bigskip}
\newcommand\sni{}
\newcommand\snii{}
\newcommand\mni{\ms \noindent}
\newcommand\bni{\bs \noindent}
\newcommand\noi{\noindent}
\newcommand\alb{\allowbreak}

\newcommand\et{\hbox{ et }}

\newcommand \eoe {\hbox{}\nobreak\hfill
\vrule width 1.4mm height 1.4mm depth 0mm \par \smallskip}

\newcommand\dsp{\displaystyle}
\newcommand\ndsp{\textstyle}
%:   fleches longues
\newcommand{\llongrightarrow}{\relbar\joinrel\mkern-1mu\lora}
\newcommand{\lllongrightarrow}{\relbar\joinrel\mkern-1mu\llra}
\newcommand{\llllongrightarrow}{\relbar\joinrel\mkern-1mu\lllra}
\newcommand{\lllllongrightarrow}{\relbar\joinrel\mkern-1mu\llllra}
\newcommand \lora {\longrightarrow}
\newcommand \llra {\llongrightarrow}
\newcommand \lllra {\lllongrightarrow}
\newcommand \llllra {\llllongrightarrow}
\newcommand \lllllra {\lllllongrightarrow}
\newcommand\simarrow{\vers{_\sim}}
\newcommand\isosim {\simarrow}
\newcommand\vers[1]{\buildrel{#1}\over \lora }
\newcommand\vvers[1]{\buildrel{#1}\over \llra }
\newcommand\vvvers[1]{\buildrel{#1}\over \lllra }
\newcommand\vvvvers[1]{\buildrel{#1}\over \llllra }
\newcommand\vvvvvers[1]{\buildrel{#1}\over \lllllra }

\newcommand\dar{\downarrow}
\newcommand\uar{\uparrow}

\newcommand \lmt {\longmapsto}
\newcommand \mt {\mapsto}

%:  leqslant
\renewcommand \le{\leqslant}
\renewcommand \leq{\leqslant}
\renewcommand \preceq{\preccurlyeq}
\renewcommand \ge{\geqslant}
\renewcommand \geq{\geqslant}
\renewcommand \succeq{\succurlyeq}

%:  Deux colonnes
\newcommand \DeuxCol[2]{%
\sni\mbox{\parbox[t]{.475\textwidth}{#1}%
\hspace{.05\textwidth}%
\parbox[t]{.475\textwidth}{#2}}}

\newcommand \Deuxcol[4]{%
\sni\mbox{\parbox[t]{#1\textwidth}{#3}%
\hspace{.05\textwidth}%
\parbox[t]{#2\textwidth}{#4}}}

\newcommand\bv{\,|\,}
\newcommand\lv{\,[\,}
\newcommand\rv{\,]\,}

\newcommand\vep{\varepsilon}

\newcommand\eqdf[1]{\buildrel{#1}\over =}
\newcommand\eqdefi{\eqdf{\rm def}}
\newcommand\eqdef{\buildrel{{\rm def}}\over \Longleftrightarrow }
\newcommand\som{\sum\nolimits}
\newcommand\Som{\sum\limits}

\newcommand\num{${\rm n}^{\rm o}\,$}

%%%%%%%%%%%%%%%%%%%%%%%%%%%%%%%%%%%%%%%%%%%%%%%%%%%%%%%%%%%%%%%%%%
%:         alphabets

%:  mathbb

\renewcommand \AA{\mathbb{A}}
\newcommand \BB{\mathbb{B}}
\newcommand \CC{\mathbb{C}}
\newcommand \DD{\mathbb{D}}
\newcommand \FF{\mathbb{F}}
\newcommand \GG{\mathbb{G}}
\newcommand \KK{\mathbb{K}}
\newcommand \MM{\mathbb{M}}
\newcommand \NN{\mathbb{N}}
\newcommand \QQ{\mathbb{Q}}
\newcommand \Qpos{\QQ^{>0}}
\newcommand \UU{\mathbb{U}}
\newcommand \ZZ{\mathbb{Z}}
\newcommand \ZB{\mathbb{ZB}}
\newcommand \RR{\mathbb{R}}
\newcommand \Rpos{\RR^{>0}}

\newcommand \PP{\mathbb{P}}
\newcommand \Ker {\MA{\mathrm{Ker}}}

\newcommand \FFp{\FF_p}
\newcommand \FFq{\FF_q}
\newcommand \Mn{\MM_n}
\newcommand \Mm{\MM_m}

\newcommand \GL {\mathbb{GL}}
\newcommand \GLn {{\GL_n}}
\newcommand \Gl {\mathbf{GL}}
\newcommand \Gln {{\Gl_n}}
\newcommand \SO {\mathbb{SO}}
\newcommand \SL {\mathbb{SL}}
\newcommand \SLn {{\SL_n}}
\newcommand \EE {\mathbb{E}}
\newcommand \En {\EE_n}

%\renewcommand \PP{\mathbb{P}}

%\renewcommand \HH{\mathbb{H}}

%:  mathcal

\newcommand\cA{\mathcal{A}}
\newcommand\cB{\mathcal{B}}
\newcommand\cC{\mathcal{C}}
\newcommand\cD{\mathcal{D}}
\newcommand\cE{\mathcal{E}}
\newcommand\cF{\mathcal{F}}
\newcommand\cG{\mathcal{G}}
\newcommand\cH{\mathcal{H}}
\newcommand\cI{\mathcal{I}}
\newcommand\cJ{\mathcal{J}}
\newcommand\cK{\mathcal{K}}
\newcommand\cL{\mathcal{L}}
\newcommand\cM{\mathcal{M}}
\newcommand\cO{\mathcal{O}}
\newcommand\cP{\mathcal{P}}
\newcommand\cR{\mathcal{R}}
\newcommand\cS{\mathcal{S}}
\newcommand\cZ{\mathcal{Z}}

\newcommand\NP{\mathcal{NP}}
\newcommand\FP{\mathcal{FP}}

%:  mathscr
\newcommand\scA{\mathscr{A}}
\newcommand\scB{\mathscr{B}}
\newcommand\scC{\mathscr{C}}
\newcommand\scL{\mathscr{L}}
\newcommand\scM{\mathscr{M}}
\newcommand\scR{\mathscr{R}}
\newcommand\scS{\mathscr{S}}

%:  mathfrak

\newcommand\fa{\mathfrak{a}}
\newcommand\fb{\mathfrak{b}}
\newcommand\fc{\mathfrak{c}}
\newcommand\fd{\mathfrak{d}}
\newcommand\fe{\mathfrak{e}}
\newcommand\fD{\mathfrak{D}}
\newcommand\ff{\mathfrak{f}}
\newcommand\fF{\mathfrak{F}}
\newcommand\ffg{\mathfrak{g}}
\newcommand\fh{\mathfrak{h}}
\newcommand\fI{\mathfrak{I}}
\newcommand\fm{\mathfrak{m}}
\newcommand\fn{\mathfrak{n}}
\newcommand\fN{\mathfrak{N}}
\newcommand\fp{\mathfrak{p}}
\newcommand\fq{\mathfrak{q}}
\newcommand\fr{\mathfrak{r}}
\newcommand\fP{\mathfrak{P}}
\newcommand\fx{\mathfrak{x}}
\newcommand\fS{\mathfrak{S}}

\newcommand\NPc{$\NP$\,-\,complet }
\newcommand\NPcz{$\NP$\,-\,complet}
\newcommand\NPcs{$\NP$\,-\,complets }
\newcommand\NPcsz{$\NP$\,-\,complets}

%:  mathbf

\newcommand\HC{\textbf{HC }}

\newcommand \ga{\mathbf{a}}
\newcommand \gb{\mathbf{b}}
\newcommand \gc{\mathbf{c}}
\newcommand \gh{\mathbf{h}}
\newcommand \gk{\mathbf{k}}
\newcommand \gl{\mathbf{l}}
\newcommand \gs{\mathbf{s}}
\newcommand \gv{\mathbf{v}}
\newcommand \gw{\mathbf{w}}
\newcommand \gA{\mathbf{A}}
\newcommand \gB{\mathbf{B}}
\newcommand \gC{\mathbf{C}}
\newcommand \gD{\mathbf{D}}
\newcommand \gE{\mathbf{E}}
\newcommand \gF{\mathbf{F}}
\newcommand \gG{\mathbf{G}}
\newcommand \gK{\mathbf{K}}
\newcommand \gL{\mathbf{L}}
\newcommand \gM{\mathbf{M}}
\newcommand \gP{\mathbf{P}}
\newcommand \gQ{\mathbf{Q}}
\newcommand \gR{\mathbf{R}}
\newcommand \gS{\mathbf{S}}
\newcommand \gT{\mathbf{T}}
\newcommand \gU{\mathbf{U}}
\newcommand \gV{\mathbf{V}}
\newcommand \gX{\mathbf{X}}
\newcommand \gW{\mathbf{W}}
\newcommand \ZK{\mathbf{Z}_\mathbf{K}}

\newcommand \gZ{\mathbf{Z}}
\newcommand \gZn{\gZ_n}

%:  mathrm
\newcommand\rA{\mathrm{A}}
\newcommand\rB{\mathrm{B}}
\newcommand\rC{\mathrm{C}}
\newcommand\rc{\mathrm{c}}
\newcommand\rd{\mathrm{d}}
\newcommand\rS{\mathrm{S}}
\newcommand\rD{\mathrm{D}}
\newcommand\rE{\mathrm{E}}
\newcommand\rG{\mathrm{G}}
\newcommand\rH{\mathrm{H}}
\newcommand\rI{\mathrm{I}}
\newcommand\rM{\mathrm{M}}
\newcommand\rN{\mathrm{N}}
\newcommand\rR{\mathrm{R}}
\newcommand\Sn{\mathrm{S}_{n}}
\newcommand\rT{\mathrm{T}}

%:  souligne tout pret
\newcommand \ua {{\underline{a}}}
\newcommand \ub {{\underline{b}}}
\newcommand \ud {{\underline{d}}}
\newcommand \udel {{\underline{\delta}}}
\newcommand \uf {{\underline{f}}}
\newcommand \ug {{\underline{g}}}
\newcommand \uh {{\underline{h}}}
\newcommand \um {{\underline{m}}}
\newcommand \us {{\underline{s}}}
\newcommand \ut {{\underline{t}}}
\newcommand \uu {{\underline{u}}}
\newcommand \ux {{\underline{x}}}
\newcommand \uy {{\underline{y}}}
\newcommand \uP {{\underline{P}}}
\newcommand \uS {{\underline{S}}}
\newcommand \uX {{\underline{X}}}
\newcommand \uY {{\underline{Y}}}
\newcommand \ual {{\underline{\alpha}}}
\newcommand \uga {{\underline{\gamma}}}
\newcommand \usi {{\underline{\sigma}}}
\newcommand \uxi {{\underline{\xi}}}
\newcommand \uze {{\underline{0}}}

%:  s\'equences x_1,\ldots,x_m toutes pretes
\newcommand \an {a_1,\ldots,a_n}
\newcommand \bn {b_1,\ldots,b_n}
\newcommand \bp {b_1,\ldots,b_p}
\newcommand \azn {a_0,\ldots,a_n}
\newcommand \bzn {b_0,\ldots,b_n}
\newcommand \czn {c_0,\ldots,c_n}
\newcommand \cp {c_1,\ldots,c_p}
\newcommand \sn {s_1,\ldots,s_n}
\newcommand \un {u_1,\ldots,u_n}
\newcommand \xk {x_1,\ldots,x_k}
\newcommand \Xk {X_1,\ldots,X_k}
\newcommand \xl {x_1,\ldots,x_\ell}
\newcommand \xm {x_1,\ldots,x_m}
\newcommand \xn {x_1,\ldots,x_n}
\newcommand \xp {x_1,\ldots,x_p}
\newcommand \yp {y_1,\ldots,y_p}
\newcommand \cq {c_1,\ldots,c_q}
\newcommand \gq {g_1,\ldots,g_q}
\newcommand \yq {y_1,\ldots,y_q}
\newcommand \xq {x_1,\ldots,x_q}
\newcommand \xzk {x_0,\ldots,x_k}
\newcommand \xzn {x_0,\ldots,x_n}
\newcommand \xhn {x_0:\ldots:x_n}
\newcommand \Xn {X_1,\ldots,X_n}
\newcommand \Xzn {X_0,\ldots,X_n}
\newcommand \Xm {X_1,\ldots,X_m}
\newcommand \Xr {X_1,\ldots,X_r}
\newcommand \xr {x_1,\ldots,x_r}
\newcommand \Yr {Y_1,\ldots,Y_r}
\newcommand \Yn {Y_1,\ldots,Y_n}
\newcommand \ym {y_1,\ldots,y_m}
\newcommand \Ym {Y_1,\ldots,Y_m}
\newcommand \yk {y_1,\ldots,y_k}
\newcommand \yr {y_1,\ldots,y_r}
\newcommand \yn {y_1,\ldots,y_n}
\newcommand \zn {z_1,\ldots,z_n}
\newcommand \Zn {Z_1,\ldots,Z_n}

\newcommand \xpn {x'_1,\ldots,x'_n}
\newcommand \uxp  {{\underline{x'}}}
\newcommand \ypm {y'_1,\ldots,y'_m}
\newcommand \uyp  {{\underline{y'}}}

\newcommand \aln {\alpha_1,\ldots,\alpha_n}
\newcommand \gan {\gamma_1,\ldots,\gamma_n}
\newcommand \xin {\xi_1,\ldots,\xi_n}
\newcommand \xihn {\xi_0:\ldots:\xi_n}

\newcommand \lfs {f_1,\ldots,f_s}

\newcommand\eme{$^{\mathrm{\grave eme}}$ }

%:  anneaux de polynomes
\newcommand \AT {{\gA[T]}}
\newcommand \AX {{\gA[X]}}
\newcommand \Ax {{\gA[x]}}
\newcommand \AXn {{\gA[\Xn]}}
\newcommand \Axn {{\gA[\xn]}}
\newcommand \KXn {{\gK[\Xn]}}
\newcommand \ZZXn {{\ZZ[\Xn]}}
\newcommand \AY {{\gA[Y]}}
\newcommand \BX {{\gB[X]}}
\newcommand \BY {{\gB[Y]}}
\newcommand \kX {{\gk[X]}}
\newcommand \kT {{\gk[T]}}
\newcommand \KT {{\gK[T]}}
\newcommand \KX {{\gK[X]}}
\newcommand \Kfi {{\gK[\varphi]}}
\newcommand \Kx {{\gK[x]}}
\newcommand \KY {{\gK[Y]}}
\newcommand \QQX {{\QQ[X]}}
\newcommand \VX {{\gV[X]}}
\newcommand \ZX {{\gZ[X]}}
\newcommand \Zx {{\gZ[x]}}
\newcommand \ZZX {{\ZZ[X]}}
\newcommand \ZZx {{\ZZ[x]}}
\newcommand \ZZxi {{\ZZ[\xi]}}

\newcommand \Aux {{\gA[\ux]}}
\newcommand \Auy {{\gA[\uy]}}
\newcommand \Bux {{\gB[\ux]}}
\newcommand \Kuy {{\gK[\uy]}}
\newcommand \Kuu {{\gK[\uu]}}
\newcommand \Kux {{\gK[\ux]}}

\newcommand \AuS {{\gA[\uS]}}
\newcommand \Ausi {{\gA[\usi]}}

\newcommand \AuX {{\gA[\uX]}}
\newcommand \BuX {{\gB[\uX]}}
\newcommand \KuX {{\gK[\uX]}}

\newcommand \AuY {{\gA[\uY]}}
\newcommand \BuY {{\gB[\uY]}}
\newcommand \KuY {{\gK[\uY]}}

%:   commentaires ...
\newcommand\comm{\noi\rdb{\it Commentaire. }}
\newcommand\comms{\noi\rdb{\it Commentaires. }}
\newcommand\REM[1]{\noi\rdb{\it Remarque#1 }}
\newcommand\rem{\noi\rdb{\it Remarque. }}
\newcommand\rep{\noi\rdb{\bf R\'eponse. }}
\newcommand\rems{\noi\rdb{\it Remarques. }}
\newcommand\exl{\noi\rdb{\bf Exemple. }}
\newcommand\EXL[1]{\noi\rdb{\bf Exemple #1}}
\newcommand\exls{\noi\rdb{\bf Exemples. }}

\newcommand\junk[1]{}

\newcommand\hum[1] {\sibrouillon{
{\begin{flushleft}

\texttt{\textbf{hum\,:}}
{\smalll\tt  #1\par}
\vspace{-.8mm}
\end{flushleft}}}
}

%:   mots mathematiques roman

\newcommand \Lin {\scL}

\newcommand \Id {\mathrm{Id}}
\newcommand \In {{\rI_n}}
\newcommand \I  {\mathrm{I}}

\newcommand\MA[1]{\mathop{#1}\nolimits}

\newcommand \Adj {\MA{\mathrm{Adj}}}
\newcommand \adj {\MA{\mathrm{adj}}}
\newcommand \Adu {\MA{\mathrm{Adu}}}
\newcommand \Alt {\MA{\mathrm{Alt}}}
\newcommand \Ann {\mathrm{Ann}}
\newcommand \Aut {\MA{\mathrm{Aut}}}
\newcommand \car {\MA{\mathrm{car}}}
\newcommand \Cl {\MA{\mathrm{Cl}}}
\newcommand \Coker {\MA{\mathrm{Coker}}}
\newcommand \Com {\MA{\mathrm{Com}}}
\renewcommand \det {\MA{\mathrm{det}}}
\renewcommand \deg {\MA{\mathrm{deg}}}
\newcommand \Diag {\MA{\mathrm{Diag}}}
\newcommand \di {\MA{\mathrm{di}}}
\newcommand \disc {\MA{\mathrm{disc}}}
\newcommand \Disc {\MA{\mathrm{Disc}}}
\newcommand \Div {\MA{\mathrm{Div}}}
\newcommand \dv {\MA{\mathrm{div}}}
\newcommand \ev {{\mathrm{ev}}}
\newcommand \End {\MA{\mathrm{End}}}
\newcommand \Fix {\MA{\mathrm{Fix}}}
\newcommand \Frac {\MA{\mathrm{Frac}}}
\newcommand \Gal {\MA{\mathrm{Gal}}}
\newcommand \Gfr {\MA{\mathrm{Gfr}}}
\newcommand \Gram {\MA{\mathrm{Gram}}}
\newcommand \gram {\MA{\mathrm{gram}}}
\newcommand \hauteur {\mathrm{hauteur}}
\newcommand \Iff {\MA{\mathrm{Iff}}}
\newcommand \Jac {\MA{\mathrm{Jac}}}
\newcommand \JAC {\MA{\mathrm{JAC}}}
\newcommand \Hom {\MA{\mathrm{Hom}}}
\newcommand \Ifr {\MA{\mathrm{Ifr}}}
\newcommand \Icl {\MA{\mathrm{Icl}}}
\renewcommand \Im {\MA{\mathrm{Im}}}
\newcommand \LIN {\mathrm{Lin}}
\newcommand \Mat {\MA{\mathrm{Mat}}}
\newcommand \Mip {\mathrm{Min}}
\newcommand \md {\mathrm{md}}
\newcommand \mod {\;\mathrm{mod}\;}
\newcommand \Mor {\MA{\mathrm{Mor}}}
\newcommand \poles {\hbox {\rm p\^oles}}
\newcommand \pgcd {\MA{\mathrm{pgcd}}}
\newcommand \ppcm {\MA{\mathrm{ppcm}}}
\newcommand \Rad {\MA{\mathrm{Rad}}}
\newcommand \Reg {\MA{\mathrm{Reg}}}
\newcommand \rg{\MA{\mathrm{rg}}}
\newcommand \rgd{\MA{\mathrm{rgd}}}
\newcommand \Res {\mathrm{Res}}
\newcommand \Rs {\MA{\mathrm{Rs}}}
\newcommand \rPr{\MA{\mathrm{Pr}}}
\newcommand \Rv {\mathrm{Rv}}
\newcommand \Sat {\MA{\mathrm{Sat}}}
\newcommand \Syl {\mathrm{Syl}}
\newcommand \Stp {\MA{\mathrm{Stp}}}
\newcommand \St {\mathrm{St}}
\newcommand \Tri {\MA{\mathrm{Tri}}}
\newcommand \Tor {\MA{\mathrm{Tor}}}
\newcommand \tr {\MA{\mathrm{tr}}}
\newcommand \Tr {\MA{\mathrm{Tr}}}
\newcommand \Tsc {\MA{\mathrm{Tsch}}}
\newcommand \Um {\MA{\mathrm{Um}}}
\newcommand \val {\MA{\mathrm{val}}}

\newcommand\add{\mathrm{add}}
\newcommand\mul{\mathrm{mul}}

%  symboles maths en francais
\newcommand\ilex{\exists}  
\newcommand\tout{\forall}  
\newcommand\Vers{\longrightarrow}  % fleche simple longue
\newcommand\donne{\mapsto}         % fleche   |-->
\newcommand\Donne{\longmapsto}     % fleche   |---->
\newcommand\imp{\Rightarrow}       % signe d'implication
\newcommand\equi{\Leftrightarrow}  % signe d'equivalence
\newcommand\Imp{\Longrightarrow}   % signe d'implication long
\newcommand\Equi{\Longleftrightarrow}   % signe d'equivalence long
\newcommand\souli[1]{\underline{#1}}
\newcommand\surli[1]{\overline{#1}}
\newcommand\fois{\,\times\,}

%%%%%%%%%%%%%%%%%%%%%%%%%%%%%%%%%%%%%%%%%%%%%%%%%%%%%%%%%%%%%%%%%%
%%%%%%%%%%%%%%%%%%%%%%%%%%%%%%%%% expressions maths usuelles
\newcommand \resp{resp. }
\newcommand \cf{cf. }
\newcommand \cad{c'est-\`a-dire }
\newcommand \cadz{c'est-\`a-dire}
\newcommand \cade{c'est-\`a-dire en\-co\-re }
\newcommand \Cad{C'est-\`a-dire }
\newcommand \cnes {con\-di\-tion n\'eces\-saire et suffi\-sante }
\newcommand \hdr {hypo\-th\`ese de \recu}
\newcommand \hdrz {hypo\-th\`ese de \recuz}
\newcommand \ssi {si, et seu\-le\-ment si, }
\newcommand \ssiz {si, et seu\-le\-ment si,}
\newcommand \spdg {sans per\-te de g\'en\'e\-ra\-lit\'e }
\newcommand \spdgz {sans per\-te de g\'en\'e\-ra\-lit\'e}
\newcommand \Spdg {Sans per\-te de g\'en\'e\-ra\-lit\'e }

\newcommand \Propeq {Les pro\-pri\-\'e\-t\'es sui\-van\-tes sont
\'equi\-val\-en\-tes.}
\newcommand \propeq {les pro\-pri\-\'e\-t\'es sui\-van\-tes sont
\'equi\-val\-en\-tes.}

\newcommand\qed{\hfill$\sqcap\kern-8.0pt\hbox{$\sqcup$}$}
\newcommand\cqvd{\noindent{ce qu'on vou\-lait d\'e\-mon\-trer. }}

%%%%%%%%%%%%%%%%%%%%%%%%%%%%%%%%%%%%%%%%%%%%%%%%%%%%%%%%%%%%%%%%%%
%: mots mathematiques

%:  Kev Alg etc...
\newcommand \Kev {$\gK$-\evc }
\newcommand \Kevs {$\gK$-\evcs }
\newcommand \Kevz {$\gK$-\evcz}
\newcommand \Kevsz {$\gK$-\evcsz}

\newcommand \Lev {$\gL$-\evc }
\newcommand \Levs {$\gL$-\evcs }
\newcommand \Levz {$\gL$-\evcz}
\newcommand \Levsz {$\gL$-\evcsz}

\newcommand \QQev {$\QQ$-\evc }
\newcommand \QQevs {$\QQ$-\evcs }
\newcommand \QQevz {$\QQ$-\evcz}
\newcommand \QQevsz {$\QQ$-\evcsz}

\newcommand \RRev {$\RR$-\evc }
\newcommand \RRevs {$\RR$-\evcs }
\newcommand \RRevz {$\RR$-\evcz}
\newcommand \RRevsz {$\RR$-\evcsz}

\newcommand \CCev {$\CC$-\evc }
\newcommand \CCevs {$\CC$-\evcs }
\newcommand \CCevz {$\CC$-\evcz}
\newcommand \CCevsz {$\CC$-\evcsz}

\newcommand \kev {$\gk$-\evc }
\newcommand \kevs {$\gk$-\evcs }
\newcommand \kevz {$\gk$-\evcz}
\newcommand \kevsz {$\gk$-\evcsz}

\newcommand \Alg {$\gA$-\alg}
\newcommand \Algs {$\gA$-\algs}
\newcommand \Algz {$\gA$-\algz}
\newcommand \Algsz {$\gA$-\algsz}

\newcommand \Blg {$\gB$-\alg}
\newcommand \Blgs {$\gB$-\algs}
\newcommand \Blgz {$\gB$-\algz}
\newcommand \Blgsz {$\gB$-\algsz}

\newcommand \Clg {$\gC$-\alg}
\newcommand \Clgs {$\gC$-\algs}
\newcommand \Clgz {$\gC$-\algz}
\newcommand \Clgsz {$\gC$-\algsz}

\newcommand \klg {$\gk$-\alg}
\newcommand \klgs {$\gk$-\algs}
\newcommand \klgz {$\gk$-\algz}
\newcommand \klgsz {$\gk$-\algsz}

\newcommand \Klg {$\gK$-\alg}
\newcommand \Klgs {$\gK$-\algs}
\newcommand \Klgz {$\gK$-\algz}
\newcommand \Klgsz {$\gK$-\algsz}

\newcommand \ZZlg {$\ZZ$-\alg}
\newcommand \ZZlgs {$\ZZ$-\algs}
\newcommand \ZZlgz {$\ZZ$-\algz}
\newcommand \ZZlgsz {$\ZZ$-\algsz}

\newcommand \Ali {appli\-cation~$\gA$-\lin }
\newcommand \Alis {appli\-cations~$\gA$-\lins }
\newcommand \Aliz {appli\-cation~$\gA$-\linz}
\newcommand \Alisz {appli\-cations~$\gA$-\linsz}

\newcommand \kli {appli\-cation $\gk$-\lin }
\newcommand \klis {appli\-cations $\gk$-\lins }
\newcommand \kliz {appli\-cation $\gk$-\linz}
\newcommand \klisz {appli\-cations $\gk$-\linsz}

\newcommand \Kli {appli\-cation $\gK$-\lin }
\newcommand \Klis {appli\-cations $\gK$-\lins }
\newcommand \Kliz {appli\-cation $\gK$-\linz}
\newcommand \Klisz {appli\-cations $\gK$-\linsz}

\newcommand \Bli {appli\-cation~$\gB$-\lin }
\newcommand \Blis {appli\-cations~$\gB$-\lins }
\newcommand \Bliz {appli\-cation~$\gB$-\linz}
\newcommand \Blisz {appli\-cations~$\gB$-\linsz}

\newcommand \Cli {appli\-cation $\gC$-\lin }
\newcommand \Clis {appli\-cations $\gC$-\lins }
\newcommand \Cliz {appli\-cation $\gC$-\linz}
\newcommand \Clisz {appli\-cations $\gC$-\linsz}

\newcommand \Zli {appli\-cation $\gZ$-\lin }
\newcommand \Zlis {appli\-cations $\gZ$-\lins }
\newcommand \Zliz {appli\-cation $\gZ$-\linz}
\newcommand \Zlisz {appli\-cations $\gZ$-\linsz}

\newcommand \QQli {appli\-cation $\QQ$-\lin }
\newcommand \QQlis {appli\-cations $\QQ$-\lins }
\newcommand \QQliz {appli\-cation $\QQ$-\linz}
\newcommand \QQlisz {appli\-cations $\QQ$-\linsz}

\newcommand \ZZli {appli\-cation $\ZZ$-\lin }
\newcommand \ZZlis {appli\-cations $\ZZ$-\lins }
\newcommand \ZZliz {appli\-cation $\ZZ$-\linz}
\newcommand \ZZlisz {appli\-cations $\ZZ$-\linsz}

\newcommand \Amo {$\gA$-module }
\newcommand \Amos {$\gA$-modules }
\newcommand \Amoz {$\gA$-module}
\newcommand \Amosz {$\gA$-modules}

\newcommand \Bmo {$\gB$-module }
\newcommand \Bmos {$\gB$-modules }
\newcommand \Bmoz {$\gB$-module}
\newcommand \Bmosz {$\gB$-modules}

\newcommand \Cmo {$\gC$-module }
\newcommand \Cmos {$\gC$-modules }
\newcommand \Cmoz {$\gC$-module}
\newcommand \Cmosz {$\gC$-modules}

\newcommand \Kmo {$\gK$-module }
\newcommand \Kmos {$\gK$-modules }
\newcommand \Kmoz {$\gK$-module}
\newcommand \Kmosz {$\gK$-modules}

\newcommand \Zmo {$\gZ$-module }
\newcommand \Zmos {$\gZ$-modules }
\newcommand \Zmoz {$\gZ$-module}
\newcommand \Zmosz {$\gZ$-modules}

\newcommand \ZZmo {$\ZZ$-module }
\newcommand \ZZmos {$\ZZ$-modules }
\newcommand \ZZmoz {$\ZZ$-module}
\newcommand \ZZmosz {$\ZZ$-modules}

%:  a
\newcommand \abi {\apl \bil}
\newcommand \abis {\apls \bils}
\newcommand \abiz {\apl \bilz}
\newcommand \abisz {\apls \bilsz}

\newcommand \adk {an\-neau de Dedekind }
\newcommand \adks {an\-neaux de Dedekind }
\newcommand \adkz {an\-neau de Dedekind}
\newcommand \adksz {an\-neaux de Dedekind}

\newcommand \adp {an\-neau de Pr\"ufer }
\newcommand \adps {an\-neaux de Pr\"ufer }
\newcommand \adpsz {an\-neaux de Pr\"ufer}
\newcommand \adpz {an\-neau de Pr\"ufer}

\newcommand \adu {\alg de \dcn \uvle }
\newcommand \adus {\algs de \dcn \uvle }
\newcommand \aduz {\alg de \dcn \uvlez}
\newcommand \adusz {\algs de \dcn \uvlez}

\newcommand \agB {\alg de Boole }
\newcommand \agBs {\algs de Boole }
\newcommand \agBz {\alg de Boole}
\newcommand \agBsz {\algs de Boole}

\newcommand \agq{al\-g\'e\-bri\-que }
\newcommand \agqs{al\-g\'e\-bri\-ques }
\newcommand \agqz{al\-g\'e\-bri\-que}
\newcommand \agqsz{al\-g\'e\-bri\-ques}

\newcommand \agqt{al\-g\'e\-bri\-que\-ment }
\newcommand \agqtz{al\-g\'e\-bri\-que\-ment}

\newcommand\alg{al\-g\`e\-bre }
\newcommand\algs{al\-g\`e\-bres }
\newcommand\algz{al\-g\`e\-bre}
\newcommand\algsz{al\-g\`e\-bres}

\newcommand \algo{al\-go\-rith\-me }
\newcommand \algos{al\-go\-rith\-mes }
\newcommand \algoz{al\-go\-rith\-me}
\newcommand \algosz{al\-go\-rith\-mes}

\newcommand \algq{al\-go\-rith\-mi\-que }
\newcommand \algqs{al\-go\-rith\-mi\-ques }
\newcommand \algqz{al\-go\-rith\-mi\-que}
\newcommand \algqsz{al\-go\-rith\-mi\-ques}

\newcommand \ali {\apl \lin }
\newcommand \alis {\apls \lins }
\newcommand \aliz {\apl \linz}
\newcommand \alisz {\apls \linsz}

\newcommand \anar {anneau \ari}
\newcommand \anars {anneaux \aris}
\newcommand \anarsz {anneaux \arisz}
\newcommand \anarz {anneau \ariz}

\newcommand \apl {appli\-cation }
\newcommand \apls {appli\-cations }
\newcommand \aplz {appli\-cation}
\newcommand \aplsz {appli\-cations}

\newcommand \ari{arith\-m\'e\-ti\-que }
\newcommand \ariz{arith\-m\'e\-ti\-que}
\newcommand \aris{arith\-m\'e\-ti\-ques }
\newcommand \arisz{arith\-m\'e\-ti\-ques}

\newcommand \auto {au\-to\-mor\-phis\-me }
\newcommand \autos {au\-to\-mor\-phis\-mes }
\newcommand \autoz {au\-to\-mor\-phis\-me}
\newcommand \autosz {au\-to\-mor\-phis\-mes}

%:  b

\newcommand \bdf {base de \fap }
\newcommand \bdfs {bases de \fap }
\newcommand \bdfz {base de \fapz}

\newcommand \bdg {base de Gr\"obner }
\newcommand \bdgs {bases de Gr\"obner }
\newcommand \bdgz {base de Gr\"obner}
\newcommand \bdgsz {bases de Gr\"obner}

\newcommand \bif {borne inf\'e\-rieure } %
\newcommand \bifz {borne inf\'e\-rieure} %
\newcommand \bifs {bornes inf\'e\-rieures } %
\newcommand \bifsz {bornes inf\'e\-rieures} %

\newcommand \bil {bili\-n\'e\-aire }
\newcommand \bils {bili\-n\'e\-aires }
\newcommand \bilz {bili\-n\'e\-aire}
\newcommand \bilsz {bili\-n\'e\-aires}

\newcommand \bol{boo\-l\'een }
\newcommand \bole{boo\-l\'een\-ne }
\newcommand \bols{boo\-l\'eens }
\newcommand \boles{boo\-l\'een\-nes }
\newcommand \bolz{boo\-l\'een}
\newcommand \bolez{boo\-l\'een\-ne}
\newcommand \bolsz{boo\-l\'eens}
\newcommand \bolesz{boo\-l\'een\-nes}

\newcommand \bsp {borne sup\'e\-rieure } %
\newcommand \bsps {borne sup\'e\-rieures } %
\newcommand \bspz {borne sup\'e\-rieure} %
\newcommand \bspsz {borne sup\'e\-rieures} %

%:  c

\newcommand \cara{ca\-rac\-t\'e\-ris\-ti\-que }
\newcommand \caras{ca\-rac\-t\'e\-ris\-ti\-ques }
\newcommand \caraz{ca\-rac\-t\'e\-ris\-ti\-que}
\newcommand \carasz{ca\-rac\-t\'e\-ris\-ti\-ques}

\newcommand \care{carac\-t\'e\-ris\'e }
\newcommand \caree{carac\-t\'e\-ris\'ee }
\newcommand \cares{carac\-t\'e\-ris\'es }
\newcommand \carees{carac\-t\'e\-ris\'ees }

\newcommand \carn{carac\-t\'e\-ri\-sation }
\newcommand \carns{carac\-t\'e\-ri\-sations }

\newcommand \carar{carac\-t\'e\-riser }

\newcommand \cba{chan\-ge\-ment de ba\-se }
\newcommand \cbas{chan\-ge\-ments de ba\-se }
\newcommand \cbaz{chan\-ge\-ment de ba\-se}
\newcommand \cbasz{chan\-ge\-ments de ba\-se}

\newcommand \cdr{corps de racines }
\newcommand \cdrz{corps de racines}

\newcommand \coe{coef\-fi\-cient }
\newcommand \coes{coef\-fi\-cients }
\newcommand \coez{coef\-fi\-cient}
\newcommand \coesz{coef\-fi\-cients}

\newcommand \coh {co\-h\'e\-rent }
\newcommand \cohs {co\-h\'e\-rents }
\newcommand \cohz {co\-h\'e\-rent}
\newcommand \cohsz {co\-h\'e\-rents}

\newcommand \cohc {co\-h\'e\-rence }
\newcommand \cohcz {co\-h\'e\-rence}

\newcommand \coli {com\-bi\-nai\-son \lin }
\newcommand \colis {com\-bi\-nai\-sons \lins }
\newcommand \coliz {com\-bi\-nai\-son \linz}
\newcommand \colisz {com\-bi\-nai\-sons \linsz}

\newcommand \com {co\-maxi\-maux }
\newcommand \comz {co\-maxi\-maux}
\newcommand \come {co\-maxi\-ma\-les }
\newcommand \comez {co\-maxi\-ma\-les}

\newcommand \coo {coor\-donn\'ee }
\newcommand \coos {coor\-donn\'ees }
\newcommand \cooz {coor\-donn\'ee}
\newcommand \coosz {coor\-donn\'ees}

\newcommand \cop {compl\'e\-men\-taire }
\newcommand \cops {compl\'e\-men\-taires }
\newcommand \copz {compl\'e\-men\-taire}
\newcommand \copsz {compl\'e\-men\-taires}

%:  d
\newcommand \dcn {d\'e\-com\-po\-sition }
\newcommand \dcns {d\'e\-com\-po\-sitions }
\newcommand \dcnz {d\'e\-com\-po\-sition}
\newcommand \dcnsz {d\'e\-com\-po\-sitions}

\newcommand \ddk {dimension de~Krull }
\newcommand \ddkz {dimension de~Krull}

\newcommand \ddp {do\-maine de~Pr\"ufer }
\newcommand \ddps {do\-maines de~Pr\"ufer }
\newcommand \ddpsz {do\-maines de~Pr\"ufer}
\newcommand \ddpz {do\-maine de~Pr\"ufer}

\newcommand \Demo{D\'emonstration }     
\newcommand \Demoz{D\'emonstration}     

\newcommand \dem{d\'emons\-tra\-tion }     
\newcommand \demz{d\'emons\-tra\-tion}     
\newcommand \dems{d\'emons\-tra\-tions }     
\newcommand \demsz{d\'emons\-tra\-tions}

\newcommand \deno{d\'e\-no\-mi\-na\-teur }
\newcommand \denos{d\'e\-no\-mi\-na\-teurs }
\newcommand \denoz{d\'e\-no\-mi\-na\-teur}
\newcommand \denosz{d\'e\-no\-mi\-na\-teurs}     

\newcommand \deter {d\'eter\-minant }
\newcommand \deters {d\'eter\-minants }
\newcommand \deterz {d\'eter\-minant}
\newcommand \detersz {d\'eter\-minants}

\newcommand \dfc {d\'evelop\-pement en fraction conti\-nue }
\newcommand \dfcz {d\'evelop\-pement en fraction conti\-nue}
\newcommand \dfcs {d\'evelop\-pements en fraction conti\-nue }
\newcommand \dfcsz {d\'evelop\-pements en fraction conti\-nue}

\newcommand \dfn{d\'efi\-nition }
\newcommand \dfns{d\'efi\-nitions }
\newcommand \dfnz{d\'efi\-nition}
\newcommand \dfnsz{d\'efi\-nitions}

\newcommand \discri{discri\-minant }
\newcommand \discris{discri\-minants }
\newcommand \discriz{discri\-minant}
\newcommand \discrisz{discri\-minants}

\newcommand \dit{distri\-bu\-ti\-vit\'e }
\newcommand \ditz{distri\-bu\-ti\-vit\'e}
 
\newcommand \dok {domaine de~Dedekind }
\newcommand \doks {domaines de~Dedekind }
\newcommand \dokz {domaine de~Dedekind}
\newcommand \doksz {domaines de~Dedekind}

\newcommand \dve {divi\-si\-bi\-lit\'e }
\newcommand \dvez {divi\-si\-bi\-lit\'e}

\newcommand \dvz {di\-viseur de z\'ero }
\newcommand \dvzs {di\-viseurs de z\'ero }
\newcommand \dvzz {di\-viseur de z\'ero}
\newcommand \dvzsz {di\-viseurs de z\'ero}

%:  e
\newcommand \eco {\elts \com}
\newcommand \ecoz {\elts \comz}

\newcommand \eds {exten\-sion des sca\-laires }
\newcommand \edsz {exten\-sion des sca\-laires}

\newcommand \egmt {\'ega\-lement }
\newcommand \egmtz {\'ega\-lement}

\newcommand \egt {\'ega\-lit\'e }
\newcommand \egts {\'ega\-lit\'es }
\newcommand \egtz {\'ega\-lit\'e}
\newcommand \egtsz {\'ega\-lit\'es}

\newcommand \elr{\'el\'e\-men\-taire }
\newcommand \elrs{\'el\'e\-men\-taires }
\newcommand \elrz{\'el\'e\-men\-taire}
\newcommand \elrsz{\'el\'e\-men\-taires}

\newcommand \elrt{\'el\'e\-men\-tai\-rement }

\newcommand \elt{\'el\'e\-ment }
\newcommand \elts{\'el\'e\-ments }
\newcommand \eltz{\'el\'e\-ment}
\newcommand \eltsz{\'el\'e\-ments}

\def \endo {en\-do\-mor\-phis\-me }
\def \endos {en\-do\-mor\-phis\-mes }
\def \endoz {en\-do\-mor\-phis\-me}
\def \endosz {en\-do\-mor\-phis\-mes}

\newcommand \eqn{\'equa\-tion }
\newcommand \eqns{\'equa\-tions }
\newcommand \eqnz{\'equa\-tion}
\newcommand \eqnsz{\'equa\-tions}

\newcommand \eqv {\'equi\-valent }
\newcommand \eqve {\'equi\-valente }
\newcommand \eqvs {\'equi\-valents }
\newcommand \eqves {\'equi\-valentes }
\newcommand \eqvz {\'equi\-valent}
\newcommand \eqvez {\'equi\-valente}
\newcommand \eqvsz {\'equi\-valents}
\newcommand \eqvesz {\'equi\-valentes}

\newcommand \eqvc {\'equi\-va\-lence }
\newcommand \eqvcs {\'equi\-va\-lences }
\newcommand \eqvcz {\'equi\-va\-lence}
\newcommand \eqvcsz {\'equi\-va\-lences}

\newcommand \etp{en temps po\-ly\-no\-mial }
\newcommand \etpz{en temps po\-ly\-no\-mial}

\newcommand\evc{espa\-ce vec\-to\-riel }
\newcommand\evcs{espa\-ces vec\-to\-riels }
\newcommand\evcz{espa\-ce vec\-to\-riel}
\newcommand\evcsz{espa\-ces vec\-to\-riels}

\newcommand\evn{\'eva\-lua\-tion }
\newcommand\evnz{\'eva\-lua\-tion}
\newcommand\evns{\'eva\-lua\-tions }
\newcommand\evnsz{\'eva\-lua\-tions}

%:  f,g

\newcommand \fap {\fcn par\-tiel\-le }
\newcommand \faps {\fcns par\-tiel\-les }
\newcommand \fapz {\fcn par\-tiel\-le}
\newcommand \fapsz {\fcns par\-tiel\-les}

\newcommand \fcn {facto\-ri\-sation }
\newcommand \fcns {facto\-ri\-sations }
\newcommand \fcnz {facto\-ri\-sation}
\newcommand \fcnsz {facto\-ri\-sations}

\newcommand \fdi {for\-te\-ment dis\-cret }
\newcommand \fdis {for\-te\-ment dis\-crets }
\newcommand \fdisz {for\-te\-ment dis\-crets}
\newcommand \fdiz {for\-te\-ment dis\-cret}

\newcommand\gmt{g\'eom\'e\-trie }
\newcommand\gmts{g\'eom\'e\-tries }
\newcommand\gmtz{g\'eom\'e\-trie}
\newcommand\gmtsz{g\'eom\'e\-tries}

\newcommand\gmq{g\'eom\'e\-trique }
\newcommand\gmqs{g\'eom\'e\-triques }
\newcommand\gmqz{g\'eom\'e\-trique}
\newcommand\gmqsz{g\'eom\'e\-triques}

\newcommand\gne{g\'en\'e\-ra\-lis\'e }
\newcommand\gnee{g\'en\'e\-ra\-lis\'ee }
\newcommand\gnes{g\'en\'e\-ra\-lis\'es }
\newcommand\gnees{g\'en\'e\-ra\-lis\'ees }
\newcommand\gnez{g\'en\'e\-ra\-lis\'e}
\newcommand\gneez{g\'en\'e\-ra\-lis\'ee}
\newcommand\gnesz{g\'en\'e\-ra\-lis\'es}
\newcommand\gneesz{g\'en\'e\-ra\-lis\'ees}

\newcommand\gnl{g\'en\'e\-ral }
\newcommand\gnle{g\'en\'e\-rale }
\newcommand\gnls{g\'en\'e\-raux }
\newcommand\gnles{g\'en\'e\-rales }
\newcommand\gnlz{g\'en\'e\-ral}
\newcommand\gnlez{g\'en\'e\-rale}
\newcommand\gnlsz{g\'en\'e\-raux}
\newcommand\gnlesz{g\'en\'e\-rales}

\newcommand\gnlt{g\'en\'e\-ra\-lement }
\newcommand\gnltz{g\'en\'e\-ra\-lement}

\newcommand\gnn{g\'en\'e\-ra\-li\-sa\-tion }
\newcommand\gnns{g\'en\'e\-ra\-li\-sa\-tions }
\newcommand\gnnz{g\'en\'e\-ra\-li\-sa\-tion}
\newcommand\gnnsz{g\'en\'e\-ra\-li\-sa\-tions}

\newcommand\gnq{g\'en\'e\-rique }
\newcommand\gnqs{g\'en\'e\-riques }
\newcommand\gnqz{g\'en\'e\-rique}
\newcommand\gnqsz{g\'en\'e\-riques}

\newcommand\gnr{g\'en\'e\-ra\-li\-ser }

\newcommand\gns{g\'en\'e\-ra\-li\-se }

\newcommand\gnt{g\'en\'e\-ra\-lit\'e }
\newcommand\gnts{g\'en\'e\-ra\-lit\'es }
\newcommand\gntz{g\'en\'e\-ra\-lit\'e}
\newcommand\gntsz{g\'en\'e\-ra\-lit\'es}

\newcommand\gtr{g\'en\'e\-rateur }
\newcommand\gtrs{g\'en\'e\-rateurs }
\newcommand\gtrz{g\'en\'e\-rateur}
\newcommand\gtrsz{g\'en\'e\-rateurs}

%:  h

\newcommand \hbq {hyper\-bo\-lique }
\newcommand \hbqs {hyper\-bo\-liques }
\newcommand \hbqz {hyper\-bo\-lique}
\newcommand \hbqzs {hyper\-bo\-liques}

\newcommand \homo {ho\-mo\-mor\-phis\-me }
\newcommand \homoz {ho\-mo\-mor\-phis\-me}
\newcommand \homos {ho\-mo\-mor\-phis\-mes }
\newcommand \homosz {ho\-mo\-mor\-phis\-mes}

\newcommand \hmg {homo\-g\`ene }
\newcommand \hmgs {homo\-g\`enes }
\newcommand \hmgz {homo\-g\`ene}
\newcommand \hmgsz {homo\-g\`enes}

\newcommand \hmt {ho\-mo\-th\'e\-tie }
\newcommand \hmts {ho\-mo\-th\'e\-ties }
\newcommand \hmtz {ho\-mo\-th\'e\-tie}
\newcommand \hmtsz {ho\-mo\-th\'e\-ties}

%:  i

\newcommand \ica {injection canonique }
\newcommand \icas {injections canoniques }
\newcommand \icaz {injection canonique}
\newcommand \icasz {injections canoniques}

\newcommand \icl {int\'e\-gra\-le\-ment clos }
\newcommand \iclz {int\'e\-gra\-le\-ment clos}
\newcommand \icle {int\'e\-gra\-le\-ment close }
\newcommand \iclez {int\'e\-gra\-le\-ment close}
\newcommand \icles {int\'e\-gra\-le\-ment closes }
\newcommand \iclesz {int\'e\-gra\-le\-ment closes}

\newcommand \id {id\'eal }
\newcommand \ids {id\'eaux }
\newcommand \idz {id\'eal}
\newcommand \idsz {id\'eaux}

\newcommand \ida {\idt \agq }
\newcommand \idas {\idts \agqs }
\newcommand \idasz {\idts \agqsz}
\newcommand \idaz {\idt \agqz}

\newcommand \dtl {d\'eter\-mi\-nantiel }
\newcommand \dtlz {d\'eter\-mi\-nantiel}
\newcommand \dtls {d\'eter\-mi\-nantiels }
\newcommand \dtlsz {d\'eter\-mi\-nantiels}

\newcommand \Idd {Id\'eal \dtl}
\newcommand \Idds {Id\'eaux \dtls }

\newcommand \idd {id\'eal \dtl }
\newcommand \idds {id\'eaux \dtls }
\newcommand \iddz {id\'eal \dtlz}
\newcommand \iddsz {id\'eaux \dtlsz}

\newcommand \idema {\id maxi\-mal }
\newcommand \idemas {\ids maxi\-maux }
\newcommand \idemaz {\id maxi\-mal}
\newcommand \idemasz {\ids maxi\-maux}

\newcommand \idep {\id pre\-mier }
\newcommand \idepz {\id pre\-mier}
\newcommand \ideps {\ids pre\-miers }
\newcommand \idepsz {\ids pre\-miers}

\newcommand \idf {id\'eal de Fitting }
\newcommand \idfs {id\'eaux de Fitting }
\newcommand \idfz {id\'eal de Fitting}
\newcommand \idfsz {id\'eaux de Fitting}

\newcommand \idm {idem\-po\-tent }
\newcommand \idms {idem\-po\-tents }
\newcommand \idmz {idem\-po\-tent}
\newcommand \idmsz {idem\-po\-tents}

\newcommand \idme {idem\-po\-tente }
\newcommand \idmes {idem\-po\-tentes }
\newcommand \idmez {idem\-po\-tente}
\newcommand \idmesz {idem\-po\-tentes}

\newcommand \idp {id\'e\-al prin\-ci\-pal }
\newcommand \idps {id\'e\-aux prin\-ci\-paux }
\newcommand \idpsz {id\'e\-aux prin\-ci\-paux}
\newcommand \idpz {id\'e\-al prin\-ci\-pal}

\newcommand \idt {iden\-ti\-t\'e }
\newcommand \idts {iden\-ti\-t\'es }
\newcommand \idtz {iden\-ti\-t\'e}
\newcommand \idtsz {iden\-ti\-t\'es}

\newcommand \idtr {in\-d\'e\-ter\-mi\-n\'ee }
\newcommand \idtrs {in\-d\'e\-ter\-mi\-n\'ees }
\newcommand \idtrz {in\-d\'e\-ter\-mi\-n\'ee}
\newcommand \idtrsz {in\-d\'e\-ter\-mi\-n\'ees}

\newcommand \ifo {in\-for\-ma\-ti\-que }
\newcommand \ifos {in\-for\-ma\-ti\-ques }
\newcommand \ifoz {in\-for\-ma\-ti\-que}
\newcommand \ifosz {in\-for\-ma\-ti\-ques}

\newcommand \ifr {id\'eal frac\-tion\-naire }
\newcommand \ifrs {id\'eaux frac\-tion\-naires }
\newcommand \ifrz {id\'eal frac\-tion\-naire}
\newcommand \ifrsz {id\'eaux frac\-tion\-naires}

\newcommand \imd {imm\'e\-diat }
\newcommand \imde {imm\'e\-diate }
\newcommand \imds {imm\'e\-diats }
\newcommand \imdes {imm\'e\-diates }
\newcommand \imdz {imm\'e\-diat}
\newcommand \imdez {imm\'e\-diate}
\newcommand \imdsz {imm\'e\-diats}
\newcommand \imdesz {imm\'e\-diates}

\newcommand \imdt {imm\'e\-dia\-te\-ment }
\newcommand \imdtz {imm\'e\-dia\-te\-ment}

\newcommand \iMR {\id de MacRae }
\newcommand \iMRs {\ids de MacRae }
\newcommand \iMRz {\id de MacRae}
\newcommand \iMRsz {\ids de MacRae}

\newcommand \ing {in\-ver\-se \gne }
\newcommand \ings {in\-ver\-ses \gnes }
\newcommand \ingz {in\-ver\-se \gnez}
\newcommand \ingsz {in\-ver\-ses \gnesz}

\newcommand \ird {irr\'e\-duc\-tible }
\newcommand \irds {irr\'e\-duc\-tibles }
\newcommand \irdz {irr\'e\-duc\-tible}
\newcommand \irdsz {irr\'e\-duc\-tibles}

\newcommand \isi {inva\-riant de simi\-li\-tude }
\newcommand \isis {inva\-riants de simi\-li\-tude }
\newcommand \isiz {inva\-riant de simi\-li\-tude}
\newcommand \isisz {inva\-riants de simi\-li\-tude}

\newcommand \iso {isomor\-phis\-me }
\newcommand \isos {isomor\-phis\-mes }
\newcommand \isosz {isomor\-phis\-mes}
\newcommand \isoz {isomor\-phis\-me}

\newcommand \itf {id\'e\-al \tf}
\newcommand \itfs {id\'e\-aux \tf}
\newcommand \itfz {id\'e\-al \tfz}
\newcommand \itfsz {id\'e\-aux \tfz}

\newcommand \ist {isom\'e\-trie }
\newcommand \ists {isom\'e\-tries }
\newcommand \istz {isom\'e\-trie}
\newcommand \istsz {isom\'e\-tries}

\newcommand \ivt {inver\-si\-bi\-lit\'e }
\newcommand \ivtz {inver\-si\-bi\-lit\'e}

\newcommand \iv {inver\-sible }
\newcommand \ivs {inver\-sibles }
\newcommand \ivz {inver\-sible}
\newcommand \ivsz {inver\-sibles}

%:  j, k, l

\newcommand \KRA {\index{Kronecker!astuce de ---}Kronecker }
\newcommand \KRO {\index{Kronecker!\tho de ---}Kronecker }
\newcommand \KRAz {\index{Kronecker!astuce de ---}Kronecker}
\newcommand \KROz {\index{Kronecker!\tho de ---}Kronecker}

\newcommand \lgb {local-global }
\newcommand \lgbe {locale-globale }
\newcommand \lgbes {locales-globales }
\newcommand \lgbs {local-globals }
\newcommand \lgbz {local-global}
\newcommand \lgbez {locale-globale}
\newcommand \lgbsz {local-globals}

\newcommand \lin {lin\'e\-aire }
\newcommand \lins {lin\'e\-aires }
\newcommand \linz {lin\'e\-aire}
\newcommand \linsz {lin\'e\-aires}

\newcommand \lint {lin\'e\-ai\-rement }

\newcommand \lnl {\lot simple }
\newcommand \lnls {\lot simples }
\newcommand \lnlz {\lot simple}
\newcommand \lnlsz {\lot simples}

\newcommand \lon {loca\-li\-sa\-tion }
\newcommand \lons {loca\-li\-sa\-tions }
\newcommand \lonz {loca\-li\-sa\-tion}
\newcommand \lonsz {loca\-li\-sa\-tions}

\newcommand \lop {\lot prin\-ci\-pal }
\newcommand \lops {\lot prin\-ci\-paux }
\newcommand \lopsz {\lot prin\-ci\-paux}
\newcommand \lopz {\lot prin\-ci\-pal}

\newcommand \lot {loca\-lement }
\newcommand \lotz {loca\-lement}

%:  m

\newcommand \maca {ma\-tri\-ce car\-r\'ee }
\newcommand \macas {ma\-tri\-ces car\-r\'ees }
\newcommand \macaz {ma\-tri\-ce car\-r\'ee}
\newcommand \macasz {ma\-tri\-ces car\-r\'ees}

\newcommand \maths {ma\-th\'e\-ma\-ti\-ques }
\newcommand \mathsz {ma\-th\'e\-ma\-ti\-ques}
\newcommand \mathe {ma\-th\'e\-ma\-ti\-que }
\newcommand \mathz {ma\-th\'e\-ma\-ti\-que}

\newcommand \mec{m\'eca\-ni\-quement calcu\-lable }
\newcommand \mecs{m\'eca\-ni\-quement calcu\-lables }
\newcommand \mecz{m\'eca\-ni\-quement calcu\-lable}
\newcommand \mecsz{m\'eca\-ni\-quement calcu\-lables}

\newcommand \met{m\'e\-thode }
\newcommand \mets{m\'e\-thodes }
\newcommand \metz{m\'e\-thode}
\newcommand \metsz{m\'e\-thodes}

\newcommand \mlp {ma\-tri\-ce de \lon princi\-pale }
\newcommand \mlpz {ma\-tri\-ce de \lon princi\-pale}
\newcommand \mlps {ma\-tri\-ces de \lon princi\-pale }
\newcommand \mlpsz {ma\-tri\-ces de \lon princi\-pale}

\newcommand \mnp {mani\-pu\-lation }
\newcommand \mnps {mani\-pu\-lations }
\newcommand \mnpz {mani\-pu\-lation}
\newcommand \mnpsz {mani\-pu\-lations}

\newcommand \mnr {\mnp \elr }
\newcommand \mnrs {\mnps \elrs }
\newcommand \mnrz {\mnp \elrz}
\newcommand \mnrsz {\mnps \elrsz}

\newcommand \mo {mo\-no\"{\i}de }
\newcommand \mos {mo\-no\"{\i}des }
\newcommand \mosz {mo\-no\"{\i}des}
\newcommand \moz {mo\-no\"{\i}de}

\newcommand \mom {mo\-n\^o\-me }
\newcommand \moms {mo\-n\^o\-mes }
\newcommand \momz {mo\-n\^o\-me}
\newcommand \momsz {mo\-n\^o\-mes}

\newcommand \mor{mor\-phis\-me }
\newcommand \mors{mor\-phis\-mes }
\newcommand \morz{mor\-phis\-me}
\newcommand \morsz{mor\-phis\-mes}

\newcommand \mpf {mo\-dule \pf}
\newcommand \mpfs {mo\-dules \pf}
\newcommand \mpfz {mo\-dule \pfz}
\newcommand \mpfsz {mo\-dules \pfz}

\newcommand \mpn {ma\-trice de \pn }
\newcommand \mpns {ma\-trices de \pn }
\newcommand \mpnz {ma\-trice de \pnz}
\newcommand \mpnsz {ma\-trices de \pnz}

\newcommand \mprn {ma\-trice de \prn }
\newcommand \mprns {ma\-trices de \prn }
\newcommand \mprnz {ma\-trice de \prnz}
\newcommand \mprnsz {ma\-trices de \prnz}

\newcommand \mptf {mo\-dule \ptf}
\newcommand \mptfs {mo\-dules \ptfs}
\newcommand \mptfz {mo\-dule \ptfz}
\newcommand \mptfsz {mo\-dules \ptfsz}

\newcommand \mrc {mo\-dule \prc }
\newcommand \mrcz {mo\-dule \prcz}
\newcommand \mrcs {mo\-dules \prcs }
\newcommand \mrcsz {mo\-dules \prcsz}

\newcommand \mtf {mo\-dule \tf}
\newcommand \mtfs {mo\-dules \tf}
\newcommand \mtfz {mo\-dule \tfz}
\newcommand \mtfsz {mo\-dules \tfz}

\newcommand\MT{Ma\-chi\-ne de Tu\-ring }
\newcommand\MTz{Ma\-chi\-ne de Tu\-ring}
\newcommand\MTs{Ma\-chi\-nes de Tu\-ring }
\newcommand\MTsz{Ma\-chi\-nes de Tu\-ring}

\newcommand\MTu{\MT \uvle }
\newcommand\MTuz{\MT \uvlez}

%:  n

\newcommand \ncr{n\'e\-ces\-sai\-re }
\newcommand \ncrs{n\'e\-ces\-sai\-res }
\newcommand \ncrz{n\'e\-ces\-sai\-re}
\newcommand \ncrsz{n\'e\-ces\-sai\-res}

\newcommand \ncrt{n\'e\-ces\-sai\-re\-ment }
\newcommand \ncrtz{n\'e\-ces\-sai\-re\-ment}

\newcommand \ndz {r\'egu\-lier }
\newcommand \ndzs {r\'egu\-liers }
\newcommand \ndzz {r\'egu\-lier}
\newcommand \ndzsz {r\'egu\-liers}

\newcommand \noe {noeth\'e\-rien }
\newcommand \noes {noeth\'e\-riens }
\newcommand \noee {noeth\'e\-rienne }
\newcommand \noees {noeth\'e\-riennes }
\newcommand \noez {noeth\'e\-rien}
\newcommand \noesz {noeth\'e\-riens}
\newcommand \noeez {noeth\'e\-rienne}
\newcommand \noeesz {noeth\'e\-riennes}

\newcommand \noet {noeth\'e\-ria\-nit\'e }
\newcommand \noetz {noeth\'e\-ria\-nit\'e}

%:  o

\newcommand\op{op\'e\-ra\-tion }
\newcommand\ops{op\'e\-ra\-tions }
\newcommand\opz{op\'e\-ra\-tion}
\newcommand\opsz{op\'e\-ra\-tions}
\newcommand\opari{\op\arith}
\newcommand\oparis{\ops\ariths}
\newcommand\opariz{\op\arithz}
\newcommand\oparisz{\ops\arithsz}
\newcommand\oparisv{\ops\arithsv}

\newcommand \ort{or\-tho\-go\-nal }
\newcommand \orte{or\-tho\-go\-na\-le }
\newcommand \orts{or\-tho\-go\-naux }
\newcommand \ortes{or\-tho\-go\-na\-les }
\newcommand \ortz{or\-tho\-go\-nal}
\newcommand \ortez{or\-tho\-go\-na\-le}
\newcommand \ortsz{or\-tho\-go\-naux}
\newcommand \ortesz{or\-tho\-go\-na\-les}

%:  p
\newcommand \paral{pa\-ral\-l\`e\-le }
\newcommand \parals{pa\-ral\-l\`e\-les }
\newcommand \paralz{pa\-ral\-l\`e\-le}
\newcommand \paralsz{pa\-ral\-l\`e\-les}

\newcommand \paralm{pa\-ral\-l\`e\-le\-ment }

\newcommand\pb{pro\-bl\`e\-me }
\newcommand\pbs{pro\-bl\`e\-mes }
\newcommand\pbz{pro\-bl\`e\-me}
\newcommand\pbsz{pro\-bl\`e\-mes}

\newcommand \pca {projection canonique }
\newcommand \pcas {projections canoniques }
\newcommand \pcaz {projection canonique}
\newcommand \pcasz {projections canoniques}

\newcommand \pf {de \pn finie }
\newcommand \pfz {de \pn finie}

\newcommand \plg {prin\-cipe \lgb }
\newcommand \plgs {prin\-cipes \lgbs }
\newcommand \plgz {prin\-cipe \lgbz}
\newcommand \plgsz {prin\-cipes \lgbsz}

\newcommand \pn {pr\'esen\-ta\-tion }
\newcommand \pns {pr\'esen\-ta\-tions }
\newcommand \pnz {pr\'esen\-ta\-tion}
\newcommand \pnsz {pr\'esen\-ta\-tions}

\newcommand \pog {\pol \hmg }
\newcommand \pogs {\pols \hmgs }
\newcommand \pogz {\pol \hmgz}
\newcommand \pogsz {\pols \hmgsz}

\newcommand\pol{poly\-n\^ome }
\newcommand\polz{poly\-n\^ome}
\newcommand\pols{poly\-n\^omes }
\newcommand\polsz{poly\-n\^omes}

\newcommand\Pol{Poly\-n\^ome }
\newcommand\Pols{Poly\-n\^omes }

\newcommand\Polcar{\Pol \cara}
\newcommand\polcar{\pol \cara}
\newcommand\polcarz{\pol \caraz}
\newcommand\polcars{\pols \caras }
\newcommand\polcarsz{\pols \carasz}

\newcommand\poll{poly\-no\-mial }
\newcommand\pollz{poly\-no\-mial}
\newcommand\polle{poly\-no\-mia\-le }
\newcommand\polles{poly\-no\-mia\-les }
\newcommand\pollesz{poly\-no\-mia\-les}
\newcommand\pollez{poly\-no\-mia\-le}

\newcommand\polt{poly\-no\-mia\-le\-ment }
\newcommand\poltz{poly\-no\-mia\-le\-ment}

\newcommand \polmin {\pol mini\-mal }
\newcommand \polmins {\pols mini\-maux }
\newcommand \polminsz {\pols mini\-maux}
\newcommand \polminz {\pol mini\-mal}

\newcommand \prg{pro\-gram\-me }
\newcommand \prgs{pro\-gram\-mes }
\newcommand \prgz{pro\-gram\-me}
\newcommand \prgsz{pro\-gram\-mes}

\newcommand \prmt {pr\'eci\-s\'e\-ment }
\newcommand \prmtz {pr\'eci\-s\'e\-ment}
\newcommand \Prmt {Pr\'eci\-s\'e\-ment }
\newcommand \Prmtz {Pr\'eci\-s\'e\-ment}

\newcommand \prc {\pro de rang constant }
\newcommand \prcs {\pros de rang constant }
\newcommand \prcz {\pro de rang constant}
\newcommand \prcsz {\pros de rang constant}

\newcommand \prn {pro\-jec\-tion }
\newcommand \prns {pro\-jec\-tions }
\newcommand \prnz {pro\-jec\-tion}
\newcommand \prnsz {pro\-jec\-tions}

\newcommand \pro {pro\-jec\-tif }
\newcommand \pros {pro\-jec\-tifs }
\newcommand \proz {pro\-jec\-tif}
\newcommand \prosz {pro\-jec\-tifs}

\newcommand \prr {pro\-jec\-teur }
\newcommand \prrs {pro\-jec\-teurs }
\newcommand \prrz {pro\-jec\-teur}
\newcommand \prrsz {pro\-jec\-teurs}

\newcommand \prt {pro\-pri\'et\'e }
\newcommand \prts {pro\-pri\'et\'es }
\newcommand \prtz {pro\-pri\'et\'e}
\newcommand \prtsz {pro\-pri\'et\'es}

\newcommand \ptf {\pro \tf }
\newcommand \ptfz {\pro \tfz}
\newcommand \ptfs {\pros \tf }
\newcommand \ptfsz {\pros \tfz}

%:  r

\newcommand \rde {rela\-tion de d\'epen\-dance }
\newcommand \rdes {rela\-tions de d\'epen\-dance }
\newcommand \rdesz {rela\-tions de d\'epen\-dance}
\newcommand \rdez {rela\-tion de d\'epen\-dance}

\newcommand \rdi {\rde int\'e\-grale }
\newcommand \rdis {\rdes int\'e\-grale }
\newcommand \rdiz {\rde int\'e\-grale}
\newcommand \rdisz {\rdes int\'e\-grale}

\newcommand \rdl {\rde \lin }
\newcommand \rdls {\rdes \lins }
\newcommand \rdlsz {\rdes \linsz}
\newcommand \rdlz {\rde \linz}

\newcommand\recu{r\'e\-cur\-ren\-ce }
\newcommand\recuz{r\'e\-cur\-ren\-ce}

\newcommand\reg{r\'e\-gu\-li\`e\-re }
\newcommand\regs{r\'e\-gu\-li\`e\-res }
\newcommand\regz{r\'e\-gu\-li\`e\-re}
\newcommand\regsz{r\'e\-gu\-li\`e\-res}

%:  s

\newcommand \sca {surjection canonique }
\newcommand \scas {surjections canoniques }
\newcommand \scaz {surjection canonique}
\newcommand \scasz {surjections canoniques}

\newcommand \seco {\sex courte }
\newcommand \secos {\sexs courtes }
\newcommand \secoz {\sex courte}
\newcommand \secosz {\sexs courtes}

\newcommand \sex {suite exacte }
\newcommand \sexs {suites exactes }
\newcommand \sexz {suite exacte}
\newcommand \sexsz {suites exactes}

\newcommand \sfio {\sys fonda\-mental d'\idms \orts }
\newcommand \sfios {\syss fonda\-mentaux d'\idms \orts }
\newcommand \sfioz {\sys fonda\-mental d'\idms \ortsz}
\newcommand \sfiosz {\syss fonda\-mentaux d'\idms \ortsz}

\newcommand \sgr {\sys \gtr }
\newcommand \sgrs {\syss \gtrs }
\newcommand \sgrz {\sys \gtrz}
\newcommand \sgrsz {\syss \gtrsz}

\newcommand \sli {\sys \lin }
\newcommand \slis {\syss \lins }
\newcommand \slisz {\syss \linsz}
\newcommand \sliz {\sys \linz}

\newcommand \smq {sym\'e\-trique }
\newcommand \smqs {sym\'e\-triques }
\newcommand \smqz {sym\'e\-trique}
\newcommand \smqsz {sym\'e\-triques}

\newcommand \spl {s\'epa\-rable }  % polynomes
\newcommand \spls {s\'epa\-rables }
\newcommand \splz {s\'epa\-rable}
\newcommand \splsz {s\'epa\-rables}

\newcommand \sul {suppl\'e\-men\-taire }
\newcommand \suls {suppl\'e\-men\-taires }
\newcommand \sulz {suppl\'e\-men\-taire}
\newcommand \sulsz {suppl\'e\-men\-taires}

\newcommand\supt{sup\-pl\'e\-ment\-aire }
\newcommand\supts{sup\-pl\'e\-ment\-aires }
\newcommand\suptz{sup\-pl\'e\-ment\-aire}
\newcommand\suptsz{sup\-pl\'e\-ment\-aires}

\newcommand \syc {\sys de \coos }
\newcommand \sycs {\syss de \coos }
\newcommand \sycz {\sys de \coosz}
\newcommand \sycsz {\syss de \coosz}

\newcommand \sys {sys\-t\`e\-me }
\newcommand \syss {sys\-t\`e\-mes }
\newcommand \sysz {sys\-t\`e\-me}
\newcommand \syssz {sys\-t\`e\-mes}

\newcommand \syse {\sys d'\'equa\-tions }
\newcommand \syses {\syss d'\'equa\-tions }
\newcommand \sysez {\sys d'\'equa\-tions}
\newcommand \sysesz {\syss d'\'equa\-tions}

\newcommand \Tfa {\Tho fondamental de l'\alg }
\newcommand \Tfaz {\Tho fondamental de l'\algz}
\newcommand \tfa {\tho fondamental de l'\alg }
\newcommand \tfaz {\tho fondamental de l'\algz}

%:  t,u ...

\newcommand \tf {de ty\-pe fi\-ni }
\newcommand \tfz {de ty\-pe fi\-ni}

\newcommand \Tho {Th\'eo\-r\`eme }
\newcommand \tho {th\'eo\-r\`eme }
\newcommand \thos {th\'eo\-r\`emes }
\newcommand \thoz {th\'eo\-r\`eme}
\newcommand \thosz {th\'eo\-r\`emes}

\newcommand \tvi {\tho des valeurs inter\-m\'e\-diaires }
\newcommand \tviz {\tho des valeurs inter\-m\'e\-diaires}
\newcommand \Tvi {\Tho des valeurs inter\-m\'e\-diaires }
\newcommand \Tviz {\Tho des valeurs inter\-m\'e\-diaires}

\newcommand \umd {unimo\-du\-laire }
\newcommand \umds {unimo\-du\-laires }
\newcommand \umdz {unimo\-du\-laire}
\newcommand \umdsz {unimo\-du\-laires}

\newcommand \unt {uni\-taire }
\newcommand \unts {uni\-taires }
\newcommand \untz {uni\-taire}
\newcommand \untsz {uni\-taires}

\newcommand\usl{usu\-el }
\newcommand\usle{usu\-elle }
\newcommand\usls{usu\-els }
\newcommand\usles{usu\-elles }
\newcommand\uslz{usu\-el}
\newcommand\uslez{usu\-elle}
\newcommand\uslsz{usu\-els}
\newcommand\uslesz{usu\-elles}

\newcommand\uvl{uni\-versel }
\newcommand\uvle{uni\-verselle }
\newcommand\uvls{uni\-versels }
\newcommand\uvles{uni\-verselles }
\newcommand\uvlz{uni\-versel}
\newcommand\uvlez{uni\-verselle}
\newcommand\uvlsz{uni\-versels}
\newcommand\uvlesz{uni\-verselles}

\newcommand \vmd {vecteur \umd}
\newcommand \vmds {vecteur \umds}
\newcommand \vmdz {vecteur \umdz}
\newcommand \vmdsz {vecteur \umdsz}

\newcommand \zed {z\'{e}\-ro-di\-men\-sion\-nel }
\newcommand \zedz {z\'{e}\-ro-di\-men\-sion\-nel}
\newcommand \zede {z\'{e}\-ro-di\-men\-sion\-nel\-le }
\newcommand \zedez {z\'{e}\-ro-di\-men\-sion\-nel\-le}
\newcommand \zeds {z\'{e}\-ro-di\-men\-sion\-nels }
\newcommand \zedsz {z\'{e}\-ro-di\-men\-sion\-nels}
\newcommand \zedes {z\'{e}\-ro-di\-men\-sion\-nel\-les }
\newcommand \zedesz {z\'{e}\-ro-di\-men\-sion\-nel\-les}

\newcommand \zedr {\zed r\'eduit }
\newcommand \zedrs {\zeds r\'eduits }
\newcommand \zedrz {\zed r\'eduit}
\newcommand \zedrsz {\zeds r\'eduits}

%:  constructifs
\newcommand \cof {cons\-truc\-tif }
\newcommand \cofs {cons\-truc\-tifs }
\newcommand \cofz {cons\-truc\-tif}
\newcommand \cofsz {cons\-truc\-tifs}

\newcommand \cov {cons\-truc\-ti\-ve }
\newcommand \covz {cons\-truc\-ti\-ve}
\newcommand \covsz {cons\-truc\-ti\-ves}
\newcommand \covs {cons\-truc\-ti\-ves }

\newcommand \coma {\maths\covs}
\newcommand \comaz {\maths\covsz}
\newcommand \clama {\maths clas\-si\-ques }
\newcommand \clamaz {\maths clas\-si\-ques}

\renewcommand \cot {cons\-truc\-ti\-ve\-ment }
\newcommand \cotz {cons\-truc\-ti\-ve\-ment}

\newcommand \prco {preuve \cov}
\newcommand \prcos {preuves \covs}
\newcommand \prcoz {preuve \covz}
\newcommand \prcosz {preuve \covsz}

%: Macros compliquees
\newcommand \cm{cm}
\makeatletter

%:  \\ retabli
\DeclareRobustCommand\\{%
  \let \reserved@e \relax
  \let \reserved@f \relax
  \@ifstar{\let \reserved@e \vadjust \let \reserved@f \nobreak
             \@xnewline}%
          \@xnewline}
\makeatother

%:             blocs  (2 X 2 blocs)
\newcommand{\blocs}[8]{%
{\setlength{\unitlength}{.0833\textwidth}
\tabcolsep0pt\renewcommand{\arraystretch}{0}%
\begin{tabular}{|c|c|}
\hline
\parbox[t][#3\cm][c]{#1\cm}{\begin{minipage}[c]{#1\cm}
\centering#5
\end{minipage}}&
\parbox[t][#3\cm][c]{#2\cm}{\begin{minipage}[c]{#2\cm}
\centering#6
\end{minipage}}\\
\hline
\parbox[t][#4\cm][c]{#1\cm}{\begin{minipage}[c]{#1\cm}
\centering#7
\end{minipage}}&
\parbox[t][#4\cm][c]{#2\cm}{\begin{minipage}[c]{#2\cm}
\centering#8
\end{minipage}}\\
\hline
\end{tabular}
}}

% exemple
% \blocs{.8}{.6}{.8}{.6}{$\I_k$}{$0$}{$0$}{$0$}

%:             tri (triangle)
\newcommand \tri[7]{
$$\quad\quad\quad\quad
\vcenter{\xymatrix@C=1.5cm
{
#1 \ar[d]_{#2} \ar[dr]^{#3} \\
{#4} \ar[r]_{{#5}}   & {#6} \\
}}
\quad\quad \vcenter{\hbox{\smalll {#7}}\hbox{~\\[1mm] ~ }}
$$
}

%:             carre (carr\'e)

\newcommand \carre[8]{
$$
\xymatrix @C=1.2cm{
#1\,\ar[d]^{#4}\ar[r]^{#2}   & \,#3\ar[d]^{#5}   \\
#6\,\ar[r]    ^{#7}    & \,#8  \\
}
$$
}

%: dcan
\newcommand\dcan[8]{
\xymatrix @C=1.2cm{
#1\,\ar[d]_{#4}\ar[r]^{#2}   & \,#3   \\
#6\,\ar[r]_{#7}    & #8\ar[u]_{#5}  %\\
}
}

%:             pun  fragile
\newcommand \pun[7]{
$$\quad\quad\quad\quad
\vcenter{\xymatrix@C=1.5cm
{
#1 \ar[d]_{#2} \ar[dr]^{#3} \\
{#4} \ar@{-->}[r]_{{#5}\,!}   & {#6} \\
}}
\quad\quad \vcenter{\hbox{\smalll {#7}}\hbox{~\\[1mm] ~ }}
$$
}

\newcommand \puN[8]{
$$\hspace{#8}
\vcenter{\xymatrix@C=1.5cm
{
#1 \ar[d]_{#2} \ar[dr]^{#3} \\
{#4} \ar@{-->}[r]_{{#5}\,!}   & {#6} \\
}}
\quad\quad \vcenter{\hbox{\smalll {#7}}\hbox{~\\[1mm] ~ }}
$$
}

%:              Pun et autres
\newcommand \Pun[8]{
$$\quad\quad\quad\quad
\vcenter{\xymatrix@C=1.5cm
{
#1 \ar[d]_{#2} \ar[dr]^{#3} \\
{#4} \ar@{-->}[r]_{{#5}\,!}   & {#6} \\
}}
\quad\quad
\vcenter{\hbox{\smalll {#7}}
\hbox{~\\[3.5mm] ~ }
\hbox{\smalll {#8}}
\hbox{~\\[-3.5mm] ~ }}
$$
}

%\def\PUn{\PUN}

%%%%%%%%%%%%%%%%%%%%%%%%%%%%%%%%%%%%%%%%%
\newcommand \PUN[9]{
$$\quad\quad
\vcenter{\xymatrix@C=1.5cm
{
#1 \ar[d]_{#2} \ar[dr]^{#3} \\
{#4} \ar@{-->}[r]_{{#5}\,!}   & {#6} \\
}}
\quad\quad
\vcenter{
\hbox{\smalll {#7}}
\hbox{~\\[-3mm] ~}
\hbox{\smalll {#8}}
\hbox{~\\[-3mm] ~}
\hbox{\smalll {#9}}
\hbox{~\\[-3.5mm] ~ }}
$$
}

%%%%%%%%%%%%%%%%%%%%%%%%%%%%%%%%%%%%%%%%%
\newcommand \Pnv[9]{
$$\quad\quad\quad\quad
\vcenter{\xymatrix@C=1.5cm
{
#1 \ar[d]_{#2} \ar[dr]^{#3} \\
{#4} \ar@{-->}[r]_{{#5}\,!}   & {#6} \\
}}
\quad\quad
\vcenter{
\hbox{\smalll {#7}}
\hbox{~\\[1mm] ~}
\hbox{\smalll {#8}}
\hbox{~\\[-1mm] ~}
\hbox{\smalll {#9}}
\hbox{~\\[0mm] ~ }}
$$
}

%%%%%%%%%%%%%%%%%%%%%%%%%%%%%%%%%%%%%%%%%
\newcommand \pnv[9]{
$$%\hspace{3pt}
\vcenter{\xymatrix@C=1.5cm
{
#1 \ar[d]_{#2} \ar[dr]^{#3} \\
{#4} \ar@{-->}[r]_{{#5}\,!}   & {#6} \\
}}
\quad\quad
\vcenter{
\hbox{\smalll {#7}}
\hbox{~\\[1mm] ~}
\hbox{\smalll {#8}}
\hbox{~\\[-1mm] ~}
\hbox{\smalll {#9}}
\hbox{~\\[0mm] ~ }}
$$
}

%%%%%%%%%%%%%%%%%%%%%%%%%%%%%%%%%%%%%%%%%
\newcommand \PNV[9]{
$$\quad\quad\quad\quad
\vcenter{\xymatrix@C=1.5cm
{
#1 \ar[d]_{#2} \ar[dr]^{#3} \\
{#4} \ar@{-->}[r]_{{#5}\,!}   & {#6} \\
}}
\quad\quad
\vcenter{
\vspace{4mm}
\hbox{\smalll {#7}}
\hbox{~\\[-1.7mm] ~}
\hbox{\smalll {#8}}
\hbox{~\\[-1.7mm] ~}
\hbox{\smalll {#9}}
\hbox{~\\[2mm] ~ }
}
$$
}

%%%%%%%%%%%%%%%%%%%%%%%%%%%%%%%%%%%%%%%%%

%:     SCO  (suites complementaires)
\newdimen\xyrowsp
\xyrowsp=3pt
\newcommand{\SCO}[6]{
\xymatrix @R = \xyrowsp {
                                  &1 \ar@{-}[dl] \ar@{-}[dr] \\
#3 \ar@{-}[ddr]                   &   & #6 \ar@{-}[ddl] \\
                                  &\bullet\ar@{-}[d] \\
                                  &\bullet   \\
#2 \ar@{-}[ddr] \ar@{-}[uur]      &   & #5 \ar@{-}[ddl] \ar@{-}[uul] \\
                                  &\bullet \ar@{-}[d] \\
                                  &\bullet  \\
#1 \ar@{-}[uur]                   &   & #4 \ar@{-}[uul] \\
                                  & 0 \ar@{-}[ul] \ar@{-}[ur] \\
}
}

%%%%%%%%%%%%%%%%%%%%%%%%%%%%%%%%%%%%%%%%%%%%%%%%%%%%%%%%%%%%%%%%%%%%%%%%%%%

\makeatletter
\newif\if@borderstar
\def\bordercmatrix{\@ifnextchar*{%
  \@borderstartrue\@bordercmatrix@i}{\@borderstarfalse\@bordercmatrix@i*}%
}
\def\@bordercmatrix@i*{\@ifnextchar[{%
  \@bordercmatrix@ii}{\@bordercmatrix@ii[()]}
}
\def\@bordercmatrix@ii[#1]#2{%
  \begingroup
    \m@th\@tempdima.875em\setbox\z@\vbox{%
      \def\cr{\crcr\noalign{\kern 2\p@\global\let\cr\endline}}%
      \ialign {$##$\hfil\kern.2em\kern\@tempdima&\thinspace%
      \hfil$##$\hfil&&\quad\hfil$##$\hfil\crcr\omit\strut%
      \hfil\crcr\noalign{\kern-\baselineskip}#2\crcr\omit%
      \strut\cr}}%
    \setbox\tw@\vbox{\unvcopy\z@\global\setbox\@ne\lastbox}%
    \setbox\tw@\hbox{\unhbox\@ne\unskip\global\setbox\@ne\lastbox}%
    \setbox\tw@\hbox{%
      $\kern\wd\@ne\kern-\@tempdima\left\@firstoftwo#1%
        \if@borderstar\kern.2em\else\kern -\wd\@ne\fi%
      \global\setbox\@ne\vbox{\box\@ne\if@borderstar\else\kern.2em\fi}%
      \vcenter{\if@borderstar\else\kern-\ht\@ne\fi%
        \unvbox\z@\kern-\if@borderstar2\fi\baselineskip}%
\if@borderstar\kern-2\@tempdima\kern.4em\else\,\fi\right\@secondoftwo#1 $%
    }\null\;\vbox{\kern\ht\@ne\box\tw@}%
  \endgroup
}
\makeatother

\newcommand{\incrementeexosetprob}{}
\newcommand{\finincrementeexosetprob}{}

\let\oldshowchapter\showchapter
\let\oldshowsection\showsection
\def\showchapter#1{\edef\temp{\thechapter}\def\ttemp{#1}\ifx\ttemp\temp\relax\def\showsection##1{##1}\else\let\showsection\oldshowsection\oldshowchapter{#1}\fi}
%%%%%%%%%%%%%%%%%%%%%%%%%%%%%%

 % avec cadre
%\input{Header.tex} % sans cadre
\newcommand\sibrouillon[1]{}
\newcommand\Today{\sibrouillon{, \today}}

\newcommand\pba{pseudo-base }
\newcommand\pbas{pseudo-bases }
\newcommand\pbaz{pseudo-base}
\newcommand\pbasz{pseudo-bases}

\newcommand\pma{pseudo-matrice }
\newcommand\pmas{pseudo-matrices }
\newcommand\pmaz{pseudo-matrice}
\newcommand\pmasz{pseudo-matrices}

\newcommand\cpb{changement de \pba}
\newcommand\cpbs{changement de \pbas}
\newcommand\cpbz{changement de \pbaz}
\newcommand\cpbsz{changement de \pbasz}

\newcommand\fde{\MA{\mathfrak{det}}}

\setcounter{chapter}{1}

\pagestyle{headings}
\markboth{Calcul matriciel généralisé}{Calcul matriciel généralisé}
\thispagestyle{CMcadreseul}

\vspace{2em}
%:     center
\begin{center} \label{center}

{\bf\Large Calcul matriciel généralisé 
sur les 

\smallskip \ddps}

\medskip {Gema M. D\'{\i}az-Toca, Henri Lombardi, 27 juillet 2015}

\end{center}

\vspace{2em}
\begin{quotation}{\small
\centerline{\bf Résumé}

\medskip Dans cet article nous présentons un algorithme pour calculer la forme de Hermite d'une \pma sur un \ddpz. Ceci nous donne des \prcos
des principaux résultats théoriques concernant les \mpfs et \ptfs sur le \ddpsz, et nous permet de discuter la résolution des \slis sur les \ddpsz. Nous généralisons ainsi la méthodologie développée par Henri Cohen pour les \doksz. Nous présentons \egmt des résultats concernant la réduction de Smith sur les \ddps de dimension $\leq 1$.

\bigskip \centerline{\bf Abstract}

\medskip 
In this paper, we first present an algorithm for computing the Hermite normal form of pseudo-matrices over Pr\"ufer domains. This algorithm allows us to provide constructive proofs of the main theoretical results on finitely presented modules over Pr\"ufer domains and  to discuss the resolution of linear systems. In some sense, we generalize the methodology developed by Henri Cohen for Dedekind domains. Finally, we present some results over Pr\"ufer domains of dimension one about the Smith normal form.
}
\end{quotation}

%\section{Introduction}
 
\Intro
\smallskip La résolution \algq des \slis sur les corps ou sur les anneaux principaux
est classique, et elle est \eqve au traitement des matrices par des \mnrs
(et des manipulations de Bezout dans le cas des anneaux principaux) de fa{\c c}on à les ramener à une forme réduite convenable.

On vise ici la \gnn de ce type de procédé à des \ddps arbitraires. 

\`A cet effet, on adapte pour un \ddp arbitraire le calcul matriciel \gne 
d'Henri Cohen [Cohen, Chapitre 1],
qu'il a introduit dans le cadre du traitement \algq des anneaux d'entiers de corps de nombres.

Notre but  est  de montrer qu'avec un calcul matriciel \gne
on peut obtenir pour un \ddp  l'analogue des réductions classiques 
 pour un domaine de Bezout.

Ceci donne un traitement \algq  systématique et \gui{agréable}, qui permet d'obtenir, comme conséquences du calcul matriciel considéré, les \thos \gnls sur les \ddpsz. 

En particulier cela donne un moyen pratique pour discuter les \slis
à \coes et inconnues dans le \ddp considéré.

Du point de vue du Calcul Formel, notre approche permet de traiter des situations plus générales que les approches usuelles. Par exemple, pour les anneaux d'entiers de corps de nombres, comme le véritable ingrédient de nos \algos se limite à l'inversion des \itfsz, notre approche permet de traiter le cas où la \fcn du \discri est impossible  en pratique, et plus \gnlt le cas où le calcul d'une base d'entiers sur $\ZZ$ est trop coûteux. 

Notons que pour les \doksz, et plus \gnlt pour les \ddps de dimension $\leq 1$,
on obtient des résultats plus précis, du type de la réduction de Smith des matrices sur les anneaux principaux. Cela sera discuté dans la section \ref{CasDim=1}. Les méthodes ici relèvent des calculs modulaires et sont les seules vraiment efficaces pour éviter l'explosion de la taille des objets dans les calculs intermédiaires.

%%%%%%%%%%%%%%%%%%%%%%%%%%%%%%%%%%%%%%%%%%%%%%%%%%%%%%%%%%%%%%%%%%%%%%%%%%%
\section{Quelques rappels}\label{secRappels}

Le premier paragraphe de cette section donne des \dfns classiques que nous donnons si \ncr  sous forme \covz. 
Le deuxième est consacré à des rappels concernant les calculs avec les \itfs sur les \ddpsz.  
Le troisième rappelle les principaux résultats concernant les calculs matriciels sur les domaines de Bezout. 

Notre référence pour les résultats classiques exposés dans cette section est le livre [Modules], dans lequel tous les résultats sont démontrés de manière \algqz.

%: subsec{Quelques \dfnsz}
\subsec{Quelques \dfnsz}

Un anneau $\gA$ est dit \emdf{\zedz} si l'axiome suivant est vérifié
$$
\forall a\in\gA,\ \exists n\in\NN\ \exists x\in\gA,\  x^n(1-ax)=0.
$$

\smallskip Un anneau intègre $\gA$ est dit \emdf{de dimension $\leq 1$} lorsque pour tout $b\neq 0$ dans $\gA$, l'anneau quotient $\aqo\gA b$ est \zedz.

\smallskip Sur un anneau $\gA$ arbitraire, un \itf $\fa=\gen{\an}$ est dit \emdf{\lopz}  s'il existe
 $s_1$, \dots, $s_n\in\gA$  tels
que \smash{$\sum_{i\in\lrbn}s_i=1$}  et $s_i\fa\subseteq \gen{a_i}$ pour chaque $s_i$.

\smallskip Un anneau $\gA$ arbitraire est dit \emdf{arithmétique}  si tout \itf est \lopz.

\smallskip Un \itf $\fa=\gen{\an}$ est dit \emdf{\ivz} s'il existe un \elt \ndz $c$ et un \itf 
$\fb$ tel que $\fa\fb=\gen{c}$. Il revient au même de dire que $\fa$ est \lop et contient un \elt \ndzz.

\smallskip Un \emdf{\ddpz} est un \anar intègre. De manière \eqvez, c'est un anneau intègre dans lequel tout \itf non nul est \ivz.
Dans ce cas, pour tout $a\neq 0$ dans l'\itf $\fa$, le transporteur $(\gen{a}:\fa)$  est \tfz, avec l'\egt $\fa\,(\gen{a}:\fa)=\gen{a}$.

\smallskip Une matrice $A=(a_{ij})$ sur un anneau arbitraire est dite \emdf{en forme de Smith} lorsque: 
\begin{itemize}
\item tous les \coes 
  $a_{ij}$ sont nuls sauf éventuellement des \gui{\coes diagonnaux} $a_i=a_{ii}$,
\item  et $a_i$ divise $a_{i+1}$ pour $i\geq 1$.
\end{itemize}

\smallskip 
 Lorsque la matrice $M$ admet une réduction de Smith, le module  $\Coker M$ est isomorphe à la somme directe des $\aqo\gA{a_i}$ et les \ids $\gen{a_i}$ sont uniquement déterminés en vertu du fait \ref{factMV-9.1}.

 Pour un anneau intègre, cela signifie que les $a_i$ sont eux-mêmes bien déterminés \gui{à une unité multiplicative près}.
 
 On appelle \emdf{anneau de Smith} un anneau $\gA$ sur lequel toute matrice~$M$ est équivalente à une matrice  en forme de Smith.

\smallskip L'\emdf{\idd d'ordre $k$} d'une matrice $M$ est l'\itf $\fD_k(M)$
engendré par les mineurs d'ordre $k$ de $M$.

\smallskip L'\emdf{\idf d'ordre $k$} d'un \mpf $P$, conoyau d'une matrice $M\in\MM_{n,m}(\gA)$, est défini par l'\egt  $\fF_k(P):=\fD_{n-k}(M)$. Il ne dépend pas de la \pn choisie pour $P$.

\smallskip On rappelle  les résultats suivants.
%:     fact factMV-9.1
\begin{fact} \label{factMV-9.1} \emph{([Modules, \thoz~V-9.1], [MRR, \thoz~V-2.4])}\\
Soient $ \fa_1\subseteq \cdots\subseteq \fa_n$ et $ \fb_1\subseteq
\cdots\subseteq \fb_m$ des id\'eaux de $\gA$ avec $n \leq m$. 
Si un module $P$ est isomorphe \`a la fois \`a
$\gA/ \fa_1\oplus\cdots\oplus \gA/ \fa_n$ et $\gA/ \fb_1\oplus\cdots\oplus
\gA/ \fb_m$, alors on a
%-----------------begin enum------------------
\begin{enumerate}
\item $ \fb_k=\gA$ pour $n < k \le m$;
\item  et  $ \fb_k= \fa_k$ pour $1 \le k \le n$.
\end{enumerate}
%-----------------end enum------------------
\end{fact}
%:     fact  factidfSmith
\begin{fact} \label{factidfSmith} \emph{([ACMC, exercice IV-16])}\\
Soit $P\simeq\gA/ \fa_1\oplus\cdots\oplus \gA/ \fa_n$ avec des \itfs
$\fa_1\supseteq \cdots\supseteq \fa_n$. Alors pour $k\in\lrbn$, on a 
$\fF_{n-k}(P)=\fa_1\cdots \fa_{k}$.
\end{fact}

%:     fact  factArizedSmith
\begin{fact} \label{factArizedSmith}  \emph{([Modules, exercice~XVI-9])}\\ Sur un \anar \zedz, toute matrice se ramène par \mnrs de lignes et de colonnes à une matrice en forme de Smith.
\end{fact}

%: subsubsec{Idéaux fractionnaires \tf}
\subsubsec{Idéaux fractionnaires \tf}
Il est souvent confortable de considérer, sur un pied d'\egt avec les 
\itfsz, les \emdf{\ifrs \tfz}. 

Pour un anneau intègre $\gA$ de corps de fractions $\gK$, nous noterons $\Iff(\gA)$ l'ensemble des \ifrs \tf de $\gA$, \cad l'ensemble des sous-\Amos \tf de $\gK$.  

Un \ifr \tf $\fb$ peut toujours s'écrire 
$$\preskip.4em \postskip.4em 
\ndsp\fb=\fraC1a \,\fa = \fraC1a\gen\an \;\hbox{avec}\;\an\in\gA, \hbox{ et }a\in\Atl.
$$ 
On dit qu'il est \lop lorsque $\fa$ est \lopz. Il admet alors la même \mlp (à \coes dans~$\gA$) que $\fa$, avec les mêmes \prtsz. 

Le produit de deux \ifrs est défini de façon naturelle et l'on a 
$x\fa\,y\fb=xy\,\fa\fb$ pour $x,y\in\gK$. 

Concernant les \emdf{transporteurs}, il faut distinguer entre 
l'\ifr 
$$\preskip.4em \postskip.4em 
 (\fa:\fb)_\gK=\sotq{x\in\gK}{x\fb\subseteq \fa}
$$ 
et l'\id (au sens usuel) $(\fa:\fb)_\gA=\sotq{x\in\gA}{x\fb\subseteq \fa}$.

Un \ifr \lop non nul $\fb$ admet  un inverse $\fb^{-1}$ au sens du produit défini dans $\Iff(\gA)$:
c'est l'\ifr $(\gA:\fb)_\gK$. 

Notons que l'application bi\lin $\theta:\fb\times \fb^{-1}\to \gA$, $(x,y)\mt xy$
est une dualité, \cad qu'elle permet d'identifier $\fb^{-1}$ au \Amo dual de $\fb$.

Un \ifr \tf contenu dans $\gA$ est un \itf usuel et il est appelé un \emdf{\id entier} si l'on est dans le contexte des \ifrsz.

%:subsec{Calculs sur les \itfs dans les \ddpsz}
\subsec{Calculs sur les \itfs dans les \ddpsz}

%:subsubsec{Matrice de \lon principale}
\subsubsec{Matrice de \lon principale}

Un anneau intègre $\gZ$ est un \ddp \emph{d'un point de vue \algqz} lorsque
l'on dispose d'un \algo \gnl qui permet d'inverser un \itf arbitraire.

Comme l'a remarqué Dedekind lorsqu'il insiste sur ce qu'il estime être la \prt la plus décisive dans les anneaux d'entiers de corps de nombres, on peut formuler l'inversibilité de l'\id $\fa=\gen{\an}\subseteq \gZ$ (les $a_i\neq 0$) de la manière suivante.

\smallskip 
\textit{1. Il existe $\gamma_1$, \ldots, $\gamma_n$ dans $\gK$ tels que 
$\sum_i\gamma_ia_i=1$ et chacun des $\gamma_ia_j$ est dans $\gZ$}.

\smallskip Deux formulations équivalentes sont les suivantes.

\smallskip \textit{2.  L'\id $\fa$ est \lopz, autrement dit il existe
 $s_1$, \dots, $s_n\in\gA$ de somme $1$    tels
que l'on \hbox{ait $s_i\fa\subseteq \gen{a_i}$} pour chaque $s_i$ non nul.}

\smallskip \textit{3. Il existe une \emdf{\mlp pour $(\an)$}
i.e., une matrice $C=(c_{ij})\in\Mn(\gZ)$ satisfaisant}
\begin{itemize}\itemsep0pt
\item $\Tr(C)=1$,
\item \textit{chaque ligne de $C$ est \emdf{proportionnelle} à $\lst{a_1\;\dots\;a_n}$, i.e.} 

\snic{\dmatrix{c_{i,\ell}&c_{i,j}\cr a_\ell&a_j}=0\hbox{ \emph{pour tous} }i, j, \ell.}

\end{itemize}

\medskip 
Pour passer de \emph{1} à \emph{2} on pose $s_i={\gamma_i}{a_i}$,
et dans l'autre sens $\gamma_i=\fraC{s_i}{a_i}$. 

\smallskip Pour passer de \emph{2} à \emph{3} on pose $c_{ij}=\fraC{s_ia_j}{a_i}$.

\smallskip La \mlp permet de nombreux calculs, comme indiqués dans 
[Modules, \thoz~\hbox{IX-2.3}, lemme\hbox{~IX-2.5}, et exercice~\hbox{XVI-8}]. 
En particulier pour un anneau intègre, on a $C^{2}=C$ (\mprn de rang $1$) et
$\In-C$ est une \mpn de l'\id $\fa$ pour le \sgr $(\an)$.
En tant que \Zmoz, $\fa$ est isomorphe à l'image de $C$ dans $ \gZ^n$.

\smallskip  Notons que, lorsque la \dve est explicite dans $\gZ$, une \mlp pour $(\an)$ est connue à partir de sa diagonale.

\fbox{Nous ferons cette hypothèse de \dve explicite dans la suite}.

Le 
 \sgr  $(\an)$ d'un \itf sera toujours accompagné d'une liste d'\elts $s_i$ de somme $1$ qui satisfont les inclusions $s_i\fa\subseteq \gen{a_i}$.
 
Sans l'hypothèse de \dve explicite, il faudrait travailler avec les \mlps elles-mêmes, qui sont des objets plus lourds à manipuler.

\smallskip   Pour qu'un anneau intègre $\gZ$  soit un \ddpz, il suffit que les \ids \emph{à deux \gtrsz} soient \ivsz, ce qui s'exprime comme suit. 

\smallskip 
\emph{Pour tous $a$, $b\in\gZ$, il existe $u$, $v$, $s$, $t\in\gZ$ satisfaisant}: 
%:     equation
\begin{equation} \label{eqmlp2} \preskip.4em \postskip.4em
s+t=1, \;sa=vb, \;tb=wa.
\end{equation}
Autrement dit, on a pour $(a,b)$ la \mlp $\cmatrix{t&w\\v&s}.$

Supposons $a\neq 0$. Alors $ts=vw$, $\gen{t,w}=\fraC t a \gen{a,b}$ et $\gen{t,w}\gen{t,v}=\gen t$. Donc, si $t,a\neq 0$,  $\gen{t,w}^{-1}=\fraC 1 t\gen{t,v}$ et $\gen{a,b}^{-1}=\fraC 1 a\gen{t,v}$ (notation des \ifrsz).

%%%%%%%%%%%%%%%%%%%%%%%%%%%%%%%%%%%%%%%%%%%%%%%%%%%%%%%%%%%%%%%%%%%%%%%%%%%
%: subsubsec{Opérations \elrs sur les \itfs}
\subsubsec{Opérations \elrs sur les \itfs}

Nous supposons travailler avec un \ddp explicite.
Autrement dit une boite noire\footnote{Dans les cas que l'on traite sur machine, cette boite noire cache un \algoz.} 
 donne pour toute suite finie non nulle $(\an)=(\ua)$ une suite finie $(\sn)=(\us)$ avec
$\sum_is_i=1$ et $s_ia_j\in\gen{a_i}$ pour \hbox{tous $i$, $j$}.

Nous supposons aussi que le \ddp est  à \dve explicite.

Tout ce qui est écrit par la suite fonctionne aussi bien avec des \ifrs \tf qu'avec des \itfs entiers.
Nous parlerons donc simplement d'\itfsz.

\paragraph{Signification de $s_k=0$.} 
\ Si $s_k=0$, en additionnant les \egts $s_i a_k=c_{i,k}a_i$ pour $i\neq k$,  on obtient l'\egt $a_k=\sum_{i\neq k}c_{i,k}a_i$, autrement dit le \gtr $a_k$ de l'\id $\fa$ est superflu. En supprimant la ligne et la colonne correspondante, on obtient une \mlp pour les $a_i$ restants.

Inversement, si l'on rajoute un \gtr \coli de ceux déjà donnés, on peut obtenir une  \mlp pour le nouveau \sgr
en rajoutant une ligne nulle pour le nouveau \gtrz, et la colonne correspondant au nouveau \gtr qui exprime la \coli comme indiqué ci-avant. 

%%%%%%%%%%%%%%%%%%%%%%%%%%%%%%%%%%%%%%%%%%%%%%%%%%%%%%%%%%%%%%%%%%%
\paragraph{Pour $\fa$ et $\fb$ \tfz, on a {\fbox{$1\in(\fa:\fb)_\gZ+(\fb:\fa)_\gZ$}}} 
(voir [Modules, exercice XVI-8]).\  
On écrit $\fa=\gen{\an}$ et $\fb=\gen{\bp}$. On a des $s_i$ de somme~$1$ tels que $s_i\fa\subseteq \gen{a_i}$, des
$t_j$ de somme $1$ tels que $t_j\fb\subseteq \gen{b_j}$, et pour chaque
couple $(i,j)$ deux \elts $u_{ij}$ et $v_{ij}$ de somme $1$ tels que

\snic{u_{ij}\!\gen{a_i,b_j}\subseteq \gen{a_i} \hbox{ et }v_{ij}\!\gen{a_i,b_j}\subseteq \gen{b_j}.}

\snii
On obtient donc 

\snic{u_{ij}s_it_j(\fa+\fb)\subseteq \gen{a_i}\subseteq \fa\hbox{ 
et  }v_{ij}s_it_j(\fa+\fb)\subseteq \gen{b_j}\subseteq  \fb.}

\snii
On prend \fbox{$s=\sum_{ij}u_{ij}s_it_j$}
et \fbox{$t=\sum_{ij}v_{ij}s_it_j$}. On a alors $s+t=1$ et

\snic{s(\fa+\fb)\subseteq \fa \hbox{ (i.e. }s\fb\subseteq \fa)\quad \hbox{et}\quad 
t(\fa+\fb)\subseteq \fb \hbox{  (i.e.  }t\fa\subseteq \fb).}

Notons que lorsque $p=1$ la \mlp pour $(b_1)$ est simplement $\I_1$.

%%%%%%%%%%%%%%%%%%%%%%%%%%%%%%%%%%%%%%%%%%%%%%%%%%%%%%%%%%%%%%%%%%%%
\paragraph{Somme (pgcd).}\ 
Avec les mêmes données qu'au point précédent, on voit que si l'on pose

\smallskip 
\centerline{\fbox{$u_i=\sum_ju_{ij}s_it_j$}\quad  et \quad  \fbox{$v_j=\sum_iv_{ij}s_it_j$},}

 on a $u_i(\fa+\fb)\subseteq \gen{a_i}$, $v_j(\fa+\fb)\subseteq \gen{b_j}$ et $\sum_iu_i+\sum_jv_j=1$.
\\
Ainsi la liste $(\un,v_1,\dots,v_p)$ convient pour le \sgr $(\an,\bp)$ de l'\id $\fa+\fb$.
\\
NB: cette procédure explique pourquoi on peut se ramener au cas où la boite noire
fonctionne uniquement pour les \ids à deux \gtrsz.

%%%%%%%%%%%%%%%%%%%%%%%%%%%%%%%%%%%%%%%%%%%%%%%%%%%%%%%%%%%%%%%%%%%
\paragraph{Intersection (ppcm).}\
Avec les mêmes données qu'aux points précédents, on a \fbox{$s\fb+t\fa=\fa\cap\fb$}. En effet, on a évidemment $s\fb\subseteq \fb\cap\fa$ et $t\fa\subseteq \fb\cap\fa$.  
Et pour l'autre inclusion $s\fb+t\fa\supseteq \fa\cap\fb$, si $m \in \fa\cap\fb$, alors 
${m = (s+t)m  = sm+tm \in s\fb+t\fa.}$

%%%%%%%%%%%%%%%%%%%%%%%%%%%%%%%%%%%%%%%%%%%%%%%%%%%%%%%%%%%%%%%%%%%
\paragraph{Inverse}\  
([Modules, lemme IX-2.5]).
 Soit $a=\sum_ia_ix_i=\lst{a_1\;\cdots \;a_n} X$, alors on a 

\snic{\gen{\ua}\gen{\ub}=\gen{a}$
avec  $\cmatrix{b_1\\ \vdots \\ b_n}=C\cmatrix{x_1\\ \vdots \\ x_n}=CX}

 (o\`u $C$ est une \mlp pour $(\ua)$).
En effet 

\snic{\dsp\som_i a_ib_i=\lst{a_1 \;\cdots \; a_n}\,C\,X=\lst{a_1 \;\cdots \; a_n}\,X=a,}

donc $a\in \fb\fa$. En outre 

\snic
{\dsp
b_i a_{j} %= \big(\som_{k} c_{ik} \,x_{k}\big) \,a_j 
= \som_{k} c_{ik} \, x_{k}\,a_j=
\som_{k}  c_{ij}\, x_{k} \,a_k =
c_{ij}\,\big(\som_{k}x_{k} a_k\big) = c_{ij}\, a,
}

donc $\fb\fa\subseteq \gen{a}$. Ainsi $\fb\fa=\gen{a}$.

%%%%%%%%%%%%%%%%%%%%%%%%%%%%%%%%%%%%%%%%%%%%%%%%%%%%%%%%%%%%%%%%%%%
\paragraph{Test d'appartenance.}\
Si $\fa=\gen{\an}$ et si $s_1,\dots,s_n\in\gZ$ vérifient $s_i\fa\subseteq \gen{a_i}$ et $\sum_is_i=1$, alors $x\in\fa$ \ssi chaque $s_ix\in\gen{a_i}$.
Comme on suppose que $\gZ$ est à \dve explicite on a un test d'appartenance aux \ifrsz, et donc aussi un test pour l'inclusion d'un \ifr dans un autre.

%%%%%%%%%%%%%%%%%%%%%%%%%%%%%%%%%%%%%%%%%%%%%%%%%%%%%%%%%%%%%%%%%%%%
%\paragraph{Changement de \sgrz.}\
%Si $\gen{\an}=\gen{b_1,\dots,b_m}$ et si l'on a $s_ia_j\in\gen{a_i}$ avec $\sum_is_i=1$,
%on peut calculer $(t_1,\dots,t_m)$ vérifiant $t_ib_j\in\gen{b_i}$ avec $\sum_it_i=1$.\\
%Pour cela on utilise l'implication \emph{1} $\Rightarrow$ \emph{2}
%dans les \prts \eqves du départ  puisque l'on dispose d'un inverse de \hbox{l'\id $\fa=\gen{b_1,\dots,b_m}$}.
%\footnote{
%On peut noter aussi qu'il suffit de savoir faire le calcul lorsque l'on rajoute ou lorsque l'on enlève un \gtr superflu.
%On a déjà traité le rajout d'un \gtr superflu.
%Lorsque l'on enlève un \gtr superflu, par exemple si $a_n=\sum_{i=1}^{n-1}\alpha_i a_i$,  \fbox{\dots\ à voir \dots\ }}

NB: un calcul  plus \gnl mais plus compliqué, sans l'hypothèse que $\gZ$
est intègre,  est donné dans l'implication \emph{1} $\Rightarrow$ \emph{3} du \tho XVI-2.2 du livre modules)

%%%%%%%%%%%%%%%%%%%%%%%%%%%%%%%%%%%%%%%%%
%:subsec{Formes réduites de matrices sur un domaine de Bezout}
\subsec{Formes réduites de matrices sur un domaine de Bezout}  

On considère un domaine de Bezout $\gB$ (cas particulier de \ddpz). 

\begin{definition} \label{defiManipBez}
Une \emdf{matrice de Bezout} de taille $n$ est une matrice
carrée, égale à la matrice identité, à l'exception de quatre \coes 
en positions $(i,i)$, $(j,j)$, $(i,j)$ et $(j,i)$ avec $1\leq i< j\leq n$ et
$\dmatrix{a_{ii}&a_{ij}\cr a_{ji}&a_{jj}}=1$. Nous pouvons noter cette matrice
$\mathrm{Bz}(n,i,j;a_{ii},a_{ij},a_{ji},a_{jj})$.
Par exemple (avec des cases vides en guise de zéros)%
\index{matrice!de Bezout}
$$\preskip.4em \postskip.4em
B=\mathrm{Bz}(6,2,4;u,v,s,t)=\cmatrix
{1&  &  &  &  &  \cr
  &u&  &v&  &  \cr
  &  &1&  &  &  \cr
  &s&  &t&  &  \cr
  &  &  &  &1&  \cr
  &  &  &  &  &1 
} \quad \hbox{avec}\quad ut-sv=1.
$$
Effectuer une  
\emdi{manipulation de Bezout} sur une matrice $M$ c'est la multiplier à gauche (manipulation de lignes) ou à droite (manipulation de colonnes)
par une matrice de Bezout. 
\end{definition}
%--------- fin definition ---------------------------------------------- 

Si l'on a dans un anneau une \egt matricielle \smashtop{$\lst{a\;  b}\Cmatrix{3.5pt}{u&v\cr s&t}=\lst{g\;  0}$}, alors pour une matrice $A\in\mat{\gA}{6,n}$ ayant des \coes $a$ et $b$ dans la ligne $k$ en positions
$(k,2)$ et $(k,4)$, effectuer le produit $A\,B$ revient à faire les manipulations de colonnes \und{simultanées} suivantes sur la matrice~$A$:
$$\preskip.4em \postskip.3em 
\formule {C_2\aff uC_2+sC_4 \\ [.5mm] %\quad \hbox{et}\quad  
C_4\aff vC_2+tC_4} 
$$
Le résultat sera (entre autres) que les  \coes en position $(k,2)$ et $(k,4)$ seront remplacés par $g$~et~$0$.

De la même manière, effectuer un produit $B'A $
par une matrice de Bezout du style  $B'=\tra B$ permettra de remplacer un couple de \coes $(a, b)$  situés sur une même colonne par un couple $(g,0)$.

%%%%%%%%%%%%%%%%%%%%%%%%%%%%%%%%%%%%%%%%%%%%%%%%%%%%%%%%%%%%%%%%%%%%
%: subsubsec{Réduction de Hermite (manipulations de colonnes, éch 
\subsubsec{Réduction de Hermite (manipulations de colonnes, échanges de lignes)}
Soit une matrice $F\in\MM_{m,n}(\gB)$. Il existe un entier $k\in\lrb{0..\inf(m,n)}$, une matrice  $C\in\GL_n(\gB)$ et une matrice de permutation $P\in\GL_m(\gB)$ telles que l'on ait
$$\preskip-.5em \postskip.3em
P\, F \, C =\blocs{.7}{.9}{.7}{.5}{$T$}{$0$}{$G$}{$0$} 
$$
avec  $T\in\MM_k(\gB)$ triangulaire inférieure de \deter $\neq 0$.

%: subsubsec{Forme réduite triangulaire (manipulations de colonnes et de lignes)}
\subsubsec{Forme réduite triangulaire (manipulations de colonnes et de lignes)}
Il existe un entier $k\in\lrb{0..\inf(m,n)}$, une matrice  $C'\in\GL_n(\gB)$ et une matrice  $L\in\GL_m(\gB)$ telles que l'on ait
$$\preskip.3em \postskip.4em
L\, F \, C' =\blocs{.7}{.9}{.7}{.5}{$T'$}{$0$}{$0$}{$0$} 
$$
avec  $T'\in\MM_k(\gB)$ triangulaire inférieure de \deter $\neq 0$. 
%

%%%%%%%%%%%%%%%%%%%%%%%%%%%%%%%%%%%%%%%%%%%%%%%%%%%%%%%%%%%%%%%%%%%%
%: subsubsec{Réduction de Smith en dimension $\leq 1$, (manipulations de lignes et de colonnes)}
\subsubsec{Réduction de Smith en dimension $\leq 1$, (manipulations de lignes et de colonnes)}
%:     theorem  thBezSmith
\begin{theorem} \label{thBezSmith}
Un domaine de Bezout de dimension $\leq 1$ est un anneau de Smith.
\end{theorem}
\begin{proof}
D'après le paragraphe précédent, il suffit de savoir traiter une matrice carrée (triangulaire) $M$  de \deter $d\neq 0$. Comme l'anneau quotient $\aqo\gB d$ est \zed et \ariz, toute matrice \hbox{sur $\aqo\gB d$} se ramène par \mnrs de lignes et de colonnes à une matrice en forme de Smith (fait \ref{factArizedSmith}).

On peut donc ramener, modulo $d$, la matrice~$M$ à une forme diagonale
de Smith, disons $\Diag(\ov{a_1},\dots,\ov{a_n})$.
Notons $\fd_k=\fD_k(M)$.

Puisque~$\gB$ est un domaine de Bezout, on a $\fd_i=\gen{d_i}$
pour des \elts $d_i$ de~$\gB$, \hbox{et $1=d_0\mid d_1\mid d_2 \mid \dots \mid d_n=d$}. On a donc des \elts $c_1$, \dots, $c_n$
tels \hbox{que $d_{k}=c_1\cdots c_k$ pour $k\in \lrbn$}. 

Quand on fait l'\eds $\gB\to \gB'=\aqo\gB d$, 
on obtient $\gen{\ov{a_1\dots a_k}}=\fD_k(\ov M)=\ov{\fd_k}$,
\cad $\gen{a_1\cdots a_k,d}=\gen{c_1\cdots c_k}$.\\
En fait, on a aussi $\gen{a_k,d}=\gen{c_k}$. Pour le voir,
on considère le \Bmo conoyau $\Coker(M)=K$. On a $dK=0$, donc 
$$
\preskip.4em \postskip.4em 
K= K/dK\simeq \aqo{\gB'}{\ov{a_1}} \oplus \cdots \oplus \aqo{\gB'}{\ov{a_n}} \simeq \aqo\gB{a_1,d}\oplus \cdots \oplus\aqo\gB{a_n,d}.
$$
En notant $a'_i=\pgcd(a_i,d)$, comme $\gen{\ov{a_i}}\subseteq \gen{\ov{a_{i-1}}}$, on a $a'_{i-1}\mid a'_i$. Cela montre\footnote{Ici on utilise le fait que deux matrices de même format qui ont des conoyaux isomorphes ont les mêmes \iddsz: ce sont les \idfs du conoyau.} que  $\fd_k=\gen{a'_1\cdots a'_k}$. Comme $\fd_k=\gen{c_1\cdots c_k}$, puisque $\gB$ est intègre et les $c_i$ non nuls, on obtient de proche en proche $\gen{a'_k}=\gen{c_k}$.

Les \mnrs modulo $d$ se relèvent en des \mnrs sur $\gB$, d'où une  \egt $LMC=M'$, où $\ov{M'}=\Diag(\ov{a_1},\dots,\ov{a_n})$. On voit que la $k$-ème ligne de $M'$ est multiple de $c_k=\pgcd(a_k,d)$. On peut donc écrire
$M'=\Diag(c_1,\dots,c_n)C_1$ avec $C_1\in\Mn(\gB)$. Et comme $\det(M')=d$, on obtient $\det(C_1)=1$. D'où enfin 
$LMC_2=\Diag(c_1,\dots,c_n)$ avec $C_2=CC_1^{-1}$. 
\end{proof}
%

%%%%%%%%%%%%%%%%%%%%%%%%%%%%%%%%%%%%%%%%%%%%%%%%%%%%%%%%%%%%%%%%%%%%
%\newpage
\section{Pseudo-bases et \pmasz}

\CAdre{.8}{Dans cette section, $\gA$ est un anneau intègre et l'on s'intéresse aux \alis entre modules qui sont sommes directes de modules \pro de rang $1$. Nous notons $\gK$ le corps de fractions}

Dans cette section on développe un \gui{calcul matriciel généralisé} pour les 
\alis entre \mptfs qui sont isomorphes à des sommes directes d'\ids \ivsz.

%:subsec{Pseudo-bases}
\subsec{Pseudo-bases}

Soit $E$ un \Amo \ptf isomorphe à une somme directe finie d'\ids \ivsz,  $E=E_1\oplus\cdots\oplus E_r$ avec $E_i\simeq \fe_i\in\Iff(\gA)$ (\ivsz). On peut 
regarder $E$ comme un sous-\Amo du \Kev $E_S$ obtenu en inversant les \elts du \mo $S=\gA\etl$. Comme $(E_i)_S$ est un \Kmo \pro de rang $1$, il est isomorphe à $\gK$ 
et l'on obtient \fbox{$E=\fe_1e_1\oplus\cdots\oplus\fe_re_r$}, o\`u \hbox{les $e_i\in E_S$} forment une base du \Kev $E_S$. \\
Notons que si $b\fe_1=a\fh_1$ pour deux \itfs $\fe_1$ et $\fh_1$ on a
l'\egt $E_1=\fe_1e_1=\fh_1 (\fraC a b e_1)=\fh_1e'_1$ ($\fe_1$ et $\fh_1$ sont deux \Amos isomorphes). Ainsi la \dcn $E=E_1\oplus\cdots\oplus E_r$ ne donne pas une écriture unique sous forme $\fe_i e_i$ pour les $E_i$.

On appelle \emdf{\pbaz} de $E$, un $r$-uplet 
$$\preskip.4em \postskip.4em 
\fbox{$\big((e_1,\fe_1),\dots,(e_r,\fe_r)\big)$} \hbox{ tel que } E=\fe_1e_1\oplus\cdots\oplus\fe_re_r, 
$$
où les $e_i$ sont des \elts de $E_S$ et les $\fe_i$ des \ifrs \ivs
de $\gA$.
Un $x\in E$ s'écrit alors de manière unique sous forme $\sum_ix_ie_i$ o\`u \hbox{les $x_i\in\fe_i$} (et donc $x_i e_i\in E$).

%:     factpba
\begin{fact} \label{factpba}
Donner une \pba $\cE=((e_1,\fe_1),\dots,(e_r,\fe_r))$ pour le \hbox{module $E$},
revient exactement à donner un \iso

\snic{\psi_\cE:\fe_1\times \cdots\times \fe_r\lora E,\; (\xr)\lmt\sum_{i\in\lrbr}x_ie_i,}

 dans lequel les images inverses des $e_i$ (après l'extension de $\psi_\cE$ à $E_S$) forment la base canonique de $\gK^{r}$
 ($\fe_1\times \cdots\times \fe_r\subseteq \gK^{r}$).
\end{fact}

\smallskip
\rem {Dans la suite, on pourrait imaginer que tous les modules \pros qui interviennent sont
des sous-\Amos de \Kevs $\gK^{n}$, o\`u $n$ est  supérieur ou égal au rang $r$ du module.}
\\
Dans un tel contexte les $e_i$ sont  soumis à la seule contrainte d'être des vecteurs \lint indépendants
 du \Kevz~$\gK^{n}$ considéré, et $E_S$ est simplement le
sous-\Kev de $\gK^{n}$ dont une base est donnée  par les vecteurs $e_i$.
Les $e_i$ sont donnés par  une matrice  $M\in\MM_{n,r}(\gK)$. On pourrait dans ce cadre aussi bien dire que la \pba $\cE$
est donnée par une matrice $M\in\MM_{n,r}(\gK)$ et une liste de $r$ \ifrs
\ivs de $\gA$.
\\
Il nous semble cependant préférable de garder pour nos explications un cadre plus \gnlz.
\eoe

%%%%%%%%%%%%%%%%%%%%%%%%%%%%%%%%%%%%%%%%%%%%%%%%%%%%%%%%%%%%%%%%%%%%%%%%%%%
%:\newpage
%\newpage
%:subsec{Matrice d'une \ali sur des pseudo-bases}
\subsec{Pseudo-matrices}

%:subsubsect{Matrice  d'une \ali  sur des pseudo-bases
\subsubsec{Matrice  d'une \ali  sur des pseudo-bases}
Rappelons que toute  cette section s'applique à un anneau intègre arbitraire~$\gA$ pour des modules qui admettent des \pbasz, \cad qui sont isomorphes à des sommes directes finies d'\ids \ivsz.
 
Soit  $\varphi:E\to H$ une \ali entre modules \pros de rangs $m$ et
$n$. Si \fbox{$\cE=\big((e_1,\fe_1),\dots,(e_m,\fe_m)\big)$} et \fbox{$\cH=\big((h_1,\fh_1),\dots,(h_n,\fh_n)\big)$} sont des \pbas de $E$ et $H$, l'application $\varphi_S:E_S\to H_S$ admet une matrice $\und A\in\gK^{n\times m}$ sur les $\gK$-bases $(e_1,\dots,e_m)$ et $(h_1,\dots,h_n)$ de $E_S$ et~$H_S$. 

\fbox{Les \coes $a_{ij}$ de cette matrice
ne sont pas  arbitraires dans~$\gK$}. 
\\
En effet, vue l'\egt $\varphi_S(e_j)=\sum_{i\in\lrbm}a_{ij}h_i$, l'inclusion $\varphi_S(E)\subseteq H$
est satisfaite \ssi pour tout $j\in\lrbm$, tout $x_j\in\fe_j$ et tout $i\in\lrbn$, on \hbox{a $a_{ij}x_j\in\fh_i$}. 

Autrement dit,
pour tous $i,j$, on doit avoir 

\smallskip \centerline{\fbox{$a_{ij}\fe_j\subseteq \fb_i$}, \hbox{ \cade }
\fbox{$a_{ij}\in\fh_i\fe_j^{-1}$}.}

%:     definota notapma
\begin{definota} \label{notapma} (Avec les notations ci-dessus)
\begin{enumerate}
\item %1
On appelle \emdf{matrice de $\varphi$ sur les \pbas $\cE$ et $\cH$}
la donnée 
$$
A=(\fh_1,\dots,\fh_n;\fe_1,\dots,\fe_m;\und A)=(\und\fh;\und\fe;\und A),\; \hbox{o\`u }\und A =(a_{ij})_{ij}\in\MM_{n,m}(\gK).
$$
Cette donnée satisfait les 
inclusions $a_{ij}\fe_j\subseteq \fh_i$.
On  note \smashbot{\fbox{$A=\scM_{\cE,\cH}(\varphi)$}}.
\item %2
On appelle \emdf{\pmaz} toute \gui{matrice} $(\und\fh;\und\fe;\und A)$ de ce type, i.e. satisfaisant les 
inclusions $a_{ij}\fe_j\subseteq \fh_i$.
\item %3
Pour des listes $\und \fe$ et $\und \fh$ fixées,  l'ensemble formé par ces \pmas est noté  \fbox{$\MM_{\und\fh;\und\fe}(\gA)$}.
\item %4
Cet ensemble est muni d'une structure naturelle de \Amoz: pour additionner $A$ et $B$, on additionne $\und A$ et $\und B$.
\item %5
On définit \egmt le produit de deux \pmas de formats convenables:

\snic{\hbox{Si } B=(\und\fc;\und\fh;\und B) \hbox{ et } A=(\und\fh;\und\fe;\und A), \;\;\hbox{alors }\; B\cdot A := (\und\fc;\und\fe;\und B\cdot \und A)}

\end{enumerate}
\end{definota}
%  end{definota}  

Dans le cas envisagé au point \emph{1} ci-dessus, les modules $\MM_{\und\fh;\und\fe}(\gA)$
et $\Lin_\gA(E,H)$ 
sont isomorphes.

 Les \pmas ont donc pour \gui{indices de colonnes} des couples $(j,\fe_j)$, pour \gui{indices de lignes} des couples $(i,\fh_i)$, et les \coes $a_{ij}$ satisfont les inclusions $a_{ij}\fe_j\subseteq \fh_i$.
On pourra visualiser les choses comme suit, par exemple pour une \pma $A$ de format $3\times 4$:
 
\CAdre{.7}{
$$\preskip-.5em \postskip.4em 
{A=\;\;\bordercmatrix [\lbrack\rbrack]{
       & \fe_1 & \fe_2  & \fe_3  & \fe_4    \cr
\fh_1 &  a_{11} &  a_{12} &  a_{13} &  a_{14}      \cr
\fh_2 &  a_{21} &  a_{22} &  a_{23} &  a_{24}      \cr
\fh_3 &  a_{31} &  a_{32} &  a_{33} &  a_{34}      \cr
}\quad a_{ij}\fe_j\subseteq \fh_i.} 
$$}

\medskip  Notez que le module $\MM_{n,m}(\gA)$ des matrices usuelles correspond au cas particulier o\`u les \ids $\fe_i$ et $\fh_j$ sont tous égaux à $\gen{1}$.

En l'absence de modules de référence $E$ et $H$, si $(e_1,\dots,e_m)$ et $(h_1,\dots,h_n)$ désignent les bases canoniques
de $\gK^{m}$ et $\gK^{n}$, la \pma $A$ est vue comme la matrice d'une \Ali $\varphi:E\to F$, o\`u 

\snic{ 
\begin{array}{cccccccccc} 
E&=&\fe_1\times \cdots\times \fe_m&=&e_1\fe_1+\cdots+ e_m\fe_m&\subseteq& \gK^{m}  &   \hbox{et}   \\[.3em] 
F&=&\fh_1\times \cdots\times \fh_n&=&h_1\fh_1+\cdots+ h_n\fh_n&\subseteq& \gK^{n} \,.  &    \end{array}
}

%:subsubsect{Calcul matriciel \gnez, quelques résultats
\subsubsec{Calcul matriciel \gnez, premiers résultats.}  

%:     Theorem{thDeter}
\begin{thdef} \label{thDeter} \emph{(Déterminant d'une \ali sur des pseudo-bases)} Avec les notations précédentes, on suppose que $E$ 
et~$H$ ont même rang~$n$ et l'on considère la matrice~\hbox{$A=\scM_{\cE,\cH}(\varphi)$} d'une \Ali $\varphi:E\to H$. On a donc $A\in\MM_{\und\fh;\und\fe}(\gA)$.
\\
Notons $\fh=\prod_{i\in\lrbn}\fh_i$ et $\fe=\prod_{i\in\lrbn}\fe_i$.
Le \deter de la matrice~$\und A$, qui est un \elt de $\gK$, est noté 
\fbox{$\det_{\cE,\cH}(\varphi)$} 
et on l'appelle \emdf{le \deter de $\varphi$ sur les \pbas $\cE$ et $\cH$}. Ci-après $\Delta:=\det_{\cE,\cH}(\varphi)=\det(\und A)$. 

\smallskip 
\begin{enumerate}
\item On a $\Delta\fe\subseteq \fh$, i.e. $\Delta\in\fh\fe^{-1}$. En particulier si $\fe=\fh$ on a $\Delta\in\gA$.
\item \Propeq
\begin{enumerate}
\item L'\ali $\varphi$ est un \isoz.
\item On a l'\egt $\Delta\fe=\fh$.
\item L'\ali  $\varphi$ est surjective.
\end{enumerate}
Dans ce cas, la matrice $A$ est une \emdf{\pma \ivz}.\\
Son inverse $A^{-1}=(\und \fe;\und\fh;{\und A}^{-1})$ est un \elt de $\MM_{\und\fe;\und\fh}(\gA)$: c'est la matrice  de l'\iso $\varphi^{-1}$ sur les \pbas $\cH$ et $\cE$.
%
%\begin{enumerate} \setcounter{enumi}{2}
%
\item L'\ali $\varphi$ est injective \ssi $\Delta\neq 0$. 
\item Lorsque $E=H$ et $\,\cE=\cH$,  $\Delta$ est un \elt de $\gA$, il ne dépend que de~$\varphi$, et
c'est le \deter de $\varphi$ au sens \gui{usuel}, noté $\det(\varphi)$.
\end{enumerate} 
\end{thdef}
%---------fin theorem----------------------------------------------
%\rem Étant donné que $\Vi^n E\simeq \fe$ et $\Vi^n H\simeq \fh$, on peut voir que la \pma à une ligne et une colonne $(\fh,\fe,\Delta)$
% représente l'\ali $\Vi^n\varphi:\Vi^n E\to \Vi^n H$ sur les \pbas déduites de $\cE$ et $\cH$.
%\eoe 
%%
\begin{proof}
 \emph{1.}  Puisque $a_{ij}\fe_j\subseteq \fh_i$ et $\Delta = \sum_{\sigma\in \cP_n} \mathrm{sign}(\sigma)  ( a_{1j_1}\cdots a_{n j_n})$, on obtient $\Delta\fe\subseteq \fh$.
 
%\smallskip 
  \emph{2a} $\Leftrightarrow$ \emph{2c.} \'Evidemment si $\varphi$ est un isomorphisme, alors $\varphi$ est un surjective.
  Si~$\varphi$ est surjective, alors $\varphi$ est scindable ([Modules, \hbox{\tho~XIII-2.1}]), donc $H\simeq E\oplus \Ker\varphi$, et puisque les modules~$E$ et~$H$ ont le même rang, $\varphi$  est un isomorphisme.
  
%\smallskip 
\emph{2a} $\Rightarrow$ \emph{2b.}  Si $\psi$ est l'\iso réciproque,
il est représenté par une matrice $B\in \MM_{\und\fe;\und\fh}(\gA)$
sur les \pbas $\cH$ et $\cE$, et l'on a $AB=(\und\fh;\und\fh;\In)$. Si on applique le point  \emph{1}  avec $\det(B) = \fraC{1}{\Delta}$, on obtient $\fraC{1}{\Delta} \fh \subseteq \fe$ et donc $  \fh \subseteq \Delta \fe$. 

%\smallskip 
\emph{2b} $\Rightarrow$ \emph{2a.}  Si $\Delta\fe=\fh$, nous allons montrer que la matrice $\und A^{-1}= \fraC1\Delta\,{\wi {\und A}}$
définit une \pma  de $\MM_{\und\fe;\und\fh}(\gA)$. Ainsi cette \pma définit un \aliz~\hbox{$\psi:H\to E$} sur les \pbas $\cH$ et $\cE$, et comme $\varphi\circ \psi$ et $\psi\circ \varphi$ sont représentées par la matrice identité (sur les \pbas convenables),~$\psi$ et~$\varphi$ sont deux \isos réciproques.

Il faut montrer \fbox{$ \fraC{1}{ \Delta }\ \wi a _{ij}\ \fh_j \subseteq \fe_i$}. On note que de manière \gnle 
%  equation label {eqAdjpma}
\begin{equation} \label {eqAdjpma}
  \wi a _{ij}\prod\nolimits_{r \neq i} \fe_r \,\subseteq\, \prod\nolimits_{s\neq j}\fh_s,
  \;\hbox{ ou encore }\;
  \fbox{$\wi a _{ij} \,\fh_j\,\fe\subseteq  \fe_i\, \fh$}. 
\end{equation}
%  end-equation
D'o\`u le résultat puisque $\Delta\fe=\fh$.

%\smallskip 
\emph{3.} Les modules $E$ et $H$ sont sans torsion et peuvent être vus comme des sous-\Amos de $\gK^{n}$. Donc  $\varphi$ est injective  \ssi $\varphi_S:\gK^{n}\to\gK^{n}$ est injective. Et $\varphi_S$ est injective \ssi  $\Delta \neq 0$.

%\smallskip 
\emph{4.} Le \deter d'un \endo se comporte bien par \edsz. En particulier lorsque l'on passe du \Amo $E$ au \Kev $E_S$ le \deter ne change pas.
Or la matrice de $\varphi$ sur la \pba $((e_1,\fe_1),\dots,(e_n,\fe_n))$ n'est autre que la matrice de $\varphi_S$ sur la base $(e_1,\dots,e_n)$.  
\end{proof}
%%%%%%%
  
\comm On notera que les deux choses suivantes sont totalement indépendantes
des \gui{vecteurs} $e_j$ et $h_i$ qui interviennent dans les \pbas $\cE$ et $\cH$.

\begin{itemize}
\item Le fait que la matrice est bien celle d'une \ali de $E$ dans~$H$: ceci est contr\^olé par les inclusions $a_{ij}\fe_j\subseteq \fh_i$.
\item Le fait que l'\ali est un \isoz: ceci est contr\^olé par
l'\egt $\Delta\fe=\fh$. 
\end{itemize}
\vspace{-2.4em}\eoe

%%%%%%%%%%%%%%%%%%%%%%%%%%%%%%%%%%%%%%%%%%%%%%%%%%%%%%%%%%%%%%%%%%%
%:subsubsec{Changement de pseudo-bases, matrice d'un tel changement}
\subsubsec{Matrice de changement de pseudo-bases}

Considérons deux \pbas différentes, \smashbot{$\cE=\big((e_1,\fe_1),\dots,(e_r,\fe_r)\big)$} 
\linebreak 
et $\cE'=\big((e'_1,\fe'_1),\dots,(e'_r,\fe'_r)\big)$,
pour un même module $E$, avec les $e_i$ et les~$e'_i$ dans le \Kev $E_S$.  

La \pma $A=\scM_{\cE',\cE}(\Id_E)$ s'appelle la \emdf{matrice de passage de~$\cE$ à~$\cE'$}: la matrice $\und A \in\MM_r(\gK)$  a pour colonnes  les $e'_i$ exprimés sur les  $e_i$.
C'est  la matrice de passage
de la base $(e_1,\dots,e_r)$ à la base $(e'_1,\dots,e'_r)$ dans~$E_S$.

L'important, pour que cela fonctionne au niveau du \Amo $E$, est que les \coes $p_{ij}$ de la matrice de passage $P$ vérifient $p_{ij}\fe'_j\in\fe_i$ et que le
déterminant $\Delta$ vérifie $\Delta\fe'=\fe$.

\medskip 
En fait, chaque fois que l'on a une \pma \iv  $P$  pour deux familles
$(\fe_1,\dots,\fe_n)$ et $(\fe'_1,\dots,\fe'_n)$,
si l'on considère un \Amo $E$ qui admet une \pba du type $\cE=\big((e_1,\fe_1),\dots,(e_n,\fe_n)\big)$, la matrice~$P$ peut être vue comme une matrice de changement de \pba pour~$E$. Il suffit pour cela de considérer les vecteurs $e'_j$
qui sont donnés par les colonnes de la matrice $\und P$ (i.e., $e'_j=\sum_{i}p_{ij}e_i$). \\
En effet, considérons le module $E'=\bigoplus_{j}e'_j\fe'_j$, 
et $\cE'=\big((e'_1,\fe'_1),\dots,(e'_n,\fe'_n)\big)$, la matrice $P$ est celle d'une \ali  $\varphi:E'\to E$ sur les \pbas $\cE'$ et $\cE$.
Mais, vu que les $e'_i$ sont donnés par les colonnes de $\und P$, \hbox{on a $\varphi_S=\Id_{E_S}$}. Ceci montre que $\varphi$ est une \ali d'inclusion de $E'$ dans~$E$. Enfin, puisque cette inclusion est un \isoz, elle est surjective et l'on obtient  $E=E'$ et $\varphi=\Id_E$.  

%%%%%%%%%%%%%%%%%%%%%%%%%%%%%%%%%%%%%%%%%%%%%%%%%%%%%%%%%%%%%%%%%%%
%:subsec{\Idds des \pmasz}
\subsec{\Idds}

On a défini le \deter d'une \ali sur des \pbasz.
On généralise les choses comme suit, en application du point~\emph{1}
du \thref{thDeter}.

%:     definition  defiiddspma
\begin{definition} \label{defiiddspma} \emph{(\Idds d'une \pmaz)}
\begin{enumerate}
\item Pour une \pma carrée $A=(\und\fh;\und\fe;\und A)$ on définit
son \emdf{\id \deterz}, noté 
$\fde(A)$: c'est  l'\id 
$${\fbox{$\fde(A)=\det(\und A)\, \fe\, \fh^{-1}$},\hbox{ o\`u }\fe=\prod\nolimits_{j}\fe_j \hbox{ et } \fh=\prod\nolimits_{i}\fh_i.}
$$
\item Soient $\beta= [\beta_1,\ldots ,\beta_r] \subseteq\lrbn$ et 
$\alpha= [\alpha_1,\ldots ,\alpha_r] \subseteq\lrbm$ des listes extraites en ordre croissant. Nous notons
$A_{\beta,\alpha}$ la \pma extraite sur les lignes (d'indices dans) $\beta$ et les colonnes
$\alpha$. \Prmt

\snic{A_{\beta,\alpha}=(\fh_{\beta_1},\dots,\fh_{\beta_r};\fe_{\alpha_1},\dots,\fe_{\alpha_r},\und A_{\beta,\alpha}).}

\snii
On dit que l'\id 

\snic{\fm_{\beta,\alpha}(A):=\fde(A_{\beta,\alpha})=\det(\und A_{\beta,\alpha}) (\prod_{i =0}^{ r} \fe_{\alpha_i})(\prod_{j=0}^r \fh_{\beta_j})^{-1}}

\snii
est \emdf{l'\id mineur
d'ordre $r$ de la matrice $A$ extrait sur les lignes
$\beta$ et les colonnes
$\alpha$}.
\item Pour une \pma arbitraire, et $r\leq \inf(m,n)$ on définit \emdf{l'\idd d'ordre $r$ de $A$}, noté $\fD_r(A)$, comme la somme des \ids mineurs d'ordre $r$ de $A$.\\
 Pour $r>\inf(m,n)$ on définit $\fD_r(A)=0$.
 \\
 Pour $r\leq 0$ on définit $\fD_r(A)=\gen{1}$.
\end{enumerate}
\end{definition}
%%%%%%%%%%%%%%%%%%%%%%%%%%%%%%%%%%%%%%%%%

%:     fact factProdDetPma
\begin{fact} \label{factProdDetPma}
 Si $A$ et $B$ sont deux \pmas carrées telles que le produit~$AB$ est défini, on a $\fde(AB)=\fde(A)\fde(B)$.\\
  Une \pma carrée est inversible \ssi son \id\deter est égal à $\gen{1}$.
\end{fact} 
\begin{proof}
On laisse le calcul du premier point \alecz. \\
Le deuxième point est donné dans le \thref{thDeter}.
\end{proof}
%
%:     thdef factIddPma
\begin{thdef} \label{factIddPma} \emph{(\Idds d'une \aliz)}
On continue avec la notation $A=(\fh_1,\dots,\fh_n;\fe_1,\dots,\fe_m;\und A)$. 
\begin{enumerate}
\item Lorsque $m=n$, l'\id \deter de la matrice $A=\scM_{\cE,\cH}(\varphi)$ ne dépend pas des \pbas choisies pour $E$ et $H$. On l'appelle \emdf{l'\id \deter de l'\ali $\varphi$}, noté $\fde(\varphi)$. 
On a les \eqvcs 
\begin{itemize}
\item $\fde(\varphi)=\gen{1}$ \ssi $\varphi$ est un \isoz,
\item $\fde(\varphi)\neq \gen{0}$ \ssi $\varphi$ est injective. 
\end{itemize}
\item Les \idds de la matrice $A$ ne dépendent pas des \pbas choisies pour $E$ et $H$. On les appelle les \emdf{\idds de l'\ali $\varphi$}, on les note $\fD_r(\varphi)$.
\item L'\ali  est surjective \ssi $\fD_n(\varphi)=\gen{1}$
\ssi elle admet un inverse à droite. Elle est injective \ssi $\fD_m(\varphi)\neq 0$. Elle admet un inverse à gauche \ssi  $\fD_m(\varphi)=\gen{1}$. 
\item Soit $s\in\gA\etl$  tel que les modules $E[1/s]$ et $H[1/s]$ sont  libres sur $\gA[1/s]$. Notons $\varphi_s:E[1/s]\to H[1/s]$ l'extension de $\varphi$ par   $\gA\to\gA[1/s]$. Alors pour chaque $r$
on a:
$\fD_r(\varphi)\gA[1/s]=\fD_r(\varphi_s)$  (\idds usuels définis pour une \ali entre modules libres). 
\item 
On a les inclusions
%---  equation eqfact.idd inc ---
\begin{equation}\label{eqfact.idd inc}\relax
\{0\}=\fD_{1+\min(m,n)}(\varphi) \subseteq\cdots\subseteq\fD_1(\varphi)
\subseteq\fD_0(\varphi)=\gen{1}=\gA
%---------------------end equation--------------
\end{equation}
Plus \prmt pour tous $k,r\in\NN$ on a une inclusion
%----- equation--eqfact.idd inc2-----------
\begin{equation}\label{eqfact.idd inc2}\relax
\fD_{k+r}(\varphi)\subseteq\fD_{k}(\varphi)\,\fD_{r}(\varphi)
\end{equation}
%---------------------end equation--------------
  
%
\end{enumerate}
\vspace{.1em}
\end{thdef}
\begin{proof} Le point \emph{1} résulte du fait \ref{factProdDetPma}
et du \thref{thDeter}.
Les autres points (comme d'ailleurs le point \emph{1}) sont conséquences du point \emph{4} (qui est facile) et des résultats analogues concernant les \alis entre modules libres. En effet, un \id est caractérisé par ses \lons en des \ecoz\footnote{Puisque l'anneau est intègre, l'\id est même l'intersection de ses localisés, vus comme sous-\Amos de $\gK$.}. Et les \prts d'injectivité et de surjectivité pour des \alis \egmtz.
\end{proof}

Un cas important est celui des \ids mineurs d'ordre $1$, 
$$
{\fbox{$\fm_{ij}(A)=a_{ij}\fe_j\fh_i^{-1}$},}
$$
qui  interviendront 
dans les manipulations de matrices conduisant à des formes réduites agréables.

%%%%%%%%%%%%%%%%%%%%%%%%%%%%%%%%%%%%%%%%%%%%%%%%%%%%%%%%%%%%%%%%%%%%
%: subsubsec{Pseudo-matrices par blocs}
\subsec{Pseudo-matrices par blocs} 
Les calculs par blocs se font pour les \pmas aussi bien que pour les 
matrices usuelles.

Considérons une \pma $M=(\fh_1,\dots,\fh_n;\fe_1,\dots,\fe_m;\und M)$ et des entiers
$p\in\lrb{2..n-1}$, $q\in\lrb{2..m-1}$. Décomposons $\und M$ en quatre blocs
$\und A\in\gK^{p\times q}$,
$\und B\in\gK^{p\times (m-q)}$, $\und C\in\gK^{(n-p)\times q}$ et $\und D\in\gK^{(n-p)\times (m-q)}$,  
$$\preskip.4em \postskip.4em
\und M=\blocs{.7}{1.0}{.8}{.6}{$\und A$}{$\und B$}{$\und C$}{$\und D$}\,
$$
ce qui donne quatre \pmas $A=(\fh_1,\dots,\fh_p;\fe_1,\dots,\fe_q;\und A)$, etc.  Alors on note
$$\preskip-.2em \postskip.4em
  M=\blocs{.7}{1.0}{.8}{.6}{$  A$}{$  B$}{$  C$}{$  D$}\,.
$$

%%%%%%%%%%%%%%%%%%%%%%%%%%%%%%%%%%%%%%%%%%%%%%%%%%%%%%%%%%%%%%%%%%%
%: subsec{Pivot de Gauss pour les \pmasz}
\subsec{Pivot de Gauss pour les \pmasz}

%:     definota notaMsc
\begin{definota} \label{notaMsc}~%
\begin{enumerate}
\item [--] Une \pma est dite  \emdf{strictement carrée}
lorsqu'elle est de la forme $A=(\und \fe;\und\fe;\und A)$. On  note l'\alg des \pmas correspondante par 

\centerline{\fbox{$\MM_{\fe_1,\dots,\fe_n}(\gA)=\MM_{\und\fe}(\gA)$}.}
\item [--] On appelle \emdf{\pma \elrz} une \pma strictement carrée $A$ telle que $\und A$ est une matrice \elr usuelle (sur~$\gK$). 
\\
Si $a_{ij}\in\gK$ est le \coe non nul en dehors de la diagonale, on doit donc avoir 
$a_{ij}\,\fe_j\subseteq \fe_i$.
\item [--] On appelle \emdf{\mnr de colonnes (resp. de lignes)} sur une \pma  la postmultiplication (resp. la prémultiplication) par une \pma \elrz.
\end{enumerate}
\end{definota}

\exls

\vspace{-1.4em}
\emph{Manipulation de colonnes.}
Considérons une \pma
$L=\;\;\bordercmatrix [\lbrack\rbrack]{
       & \fa & \fb    \cr
\fc &   a    &  b 
}\,.$

\snii
Si $\fm_{11}(L)\supseteq \fm_{12}(L)$, \cad
si $a\fa\supseteq b\fb$, on a une \pma \elr $P=\;\;\bordercmatrix [\lbrack\rbrack]{
       & \fa & \fb    \cr
\fa &   1    &  -\fraC b a \cr
\fb &   0    & 1  \cr
}$, 
de sorte que l'on obtient, comme pour le pivot de Gauss classique
$$\preskip-.4em \postskip.4em 
L P = \;\;\bordercmatrix [\lbrack\rbrack]{
       & \fa & \fb    \cr
\fc &   a    &  0 
} 
$$
Ainsi, dans le processus \gui{pivot de Gauss} ce qui tient lieu,
pour une \pmaz,  de \gui{$a_{11}$
divise $a_{12}$} dans le cas des matrices usuelles, c'est maintenant
\gui{$\fm_{11}\supseteq \fm_{12}$}, ou encore~\gui{$\fm_{12}\fm_{11}^{-1}\subseteq \gA$}, ce qui se dit aussi 
\gui{$\fm_{11}$ divise $\fm_{12}$}.

\vspace{-1em}
\emph{Manipulation de lignes.}
Considérons une \pma 
{$L=\;\;\bordercmatrix [\lbrack\rbrack]{
       & \fa      \cr
\fc_1 &   a     \cr
\fc_2 &   b     
}\,.$}
\\
Si $\fm_{11}(L)\supseteq \fm_{21}(L)$, \cad
si $a\fc_2\supseteq b\fc_1$, on a une \pma \elr $P=\;\;\bordercmatrix [\lbrack\rbrack]{
       & \fc_1 & \fc_2    \cr
\fc_1 &   1    &  0 \cr
\fc_2 &   -\fraC b a    & 1  \cr
}$, 
de sorte que l'on obtient%, comme pour le pivot de Gauss classique
\,\, $ P L = \;\;\bordercmatrix [\lbrack\rbrack]{
       & \fa      \cr
\fc_1 &   a   \cr
\fc_2 &   0   \cr
}$.

Ainsi, dans le processus \gui{pivot de Gauss} ce qui tient lieu,
pour une \pmaz,  de \gui{$a_{11}$
divise $a_{21}$} dans le cas des matrices usuelles, c'est maintenant
\gui{$\fm_{11}\supseteq \fm_{21}$}, ce qui se dit aussi 
\gui{$\fm_{11}$ divise $\fm_{21}$}.
\eoe

\medskip 
Voici une autre \gnn naturelle d'un résultat classique pour les matrices usuelles.

%:  Lem du mineur inv {lem.min.inv}
\CMnewtheorem{lemmininvpma}{Lemme du mineur inversible pour les \pmasz}{\itshape}
\begin{lemmininvpma}\label{lem.min.inv}\index{Lemme du mineur inversible pour les \pmasz}~\\
\emph{(Pivot généralisé)} On considère une \pma 
$$\preskip.4em \postskip.4em 
A=(\fh_1,\dots,\fh_n;\fe_1,\dots,\fe_m;\und A). 
$$
On suppose que l'\id mineur $\fm_{\lrbk,\lrbk}(A)$ est égal à $\gen{1}$.
\begin{enumerate}
\item \emph{(Pivot sur les colonnes)}  Il existe une  \pma \iv 
$$\preskip.4em \postskip.4em 
C=(\fe_1,\dots,\fe_m;\fh_1,\dots,\fh_k,\fe_{k+1},\dots,\fe_{m},\und C)
$$
décrite  par blocs  comme 
$C=\blocs{.7}{1.0}{.7}{.96}{$C_1$}{$C_2$}{$0$}{$\I_{m-k}$}$
%$C=\cmatrix{
%C_1 & C_2     \cr
%0 & \I_{m-k}}$ 
avec $C_1={A_{\lrbk,\lrbk}}^{-1}$
telle que la \pmaz~$AC$ soit de la forme 
$\blocs{.7}{1.0}{.7}{.6}{$\I_k$}{$0$}{$B$}{$D$}$.
%$\cmatrix{
%\I_k & 0     \cr
%B & D}.$
%
\item \emph{(Pivots sur les lignes et les colonnes)} \\
La \pma $A$ est \eqve à une \pma
$$\preskip.4em \postskip.4em 
A'=(\fh_1,\dots,\fh_n;\fh_1,\dots,\fh_k,\fe_{k+1},\dots,\fe_{m};\und{A'})=\blocs{.7}{1.0}{.7}{.6}{$\I_k$}{$0$}{$0$}{$A''$} 
$$
et l'on a $\fD_r(A'')=\fD_{k+r}(A)$ pour tout $r\in\ZZ$.
\\
En particulier si $\fD_{k+1}(A)=\gen{0}$, on a une forme réduite 
$$\preskip.4em \postskip.4em 
L\,A\,C=\blocs{.7}{1.0}{.7}{.6}{$\I_k$}{$0$}{$0$}{$0$}. 
$$
\end{enumerate}
NB: Si $A$ possède un mineur d'ordre $k>0$ égal à $\gen{1}$, on obtient des résultats analogues après permutation de lignes et/ou de colonnes.

\end{lemmininvpma}
%--------- fin lemme du mineur inversible pseudo matrices --------------- 

%----begin{proof---------------
\begin{proof} On recopie la \dem usuelle pour les matrices usuelles.

\emph{1.} 
Notons \smash{$A=\blocs{.7}{1.0}{.7}{.6}{$A_1$}{$A_2$}{$A_3$}{$A_4$}$} o\`u $\und{A_1}\in\MM_k(\gK)$. \\
La matrice $C$ est la matrice par blocs
$C=\blocs{.7}{1.0}{.7}{.9}{$C_1$}{$C_2$}{$0$}{$\I_{m-k}$}$  où $C_1={A_1}^{-1}$ \hbox{et $C_2=-C_1A_2$}. 
On obtient la matrice $A\,C$ donnée dans l'énoncé.

\emph{2.}  On continue en multipliant à gauche par la matrice
strictement carrée définie par blocs 
$L=\blocs{.7}{.6}{.7}{.6}{$\I_k$}{$0$}{$-B$}{$\I_{\ell}$}\in\MM_{\fh_1,\dots,\fh_n}(\gZ)$ avec $\ell=n-k$.
On obtient pour $L\,A\,C$ la \pma donnée dans l'énoncé. \\
Enfin $\fD_r(A'')=\fD_{k+r}(A')=\fD_{k+r}(A)$ pour \hbox{tout $r\in\ZZ$}.
\end{proof}
%----end{proof------------------

\rem
On a obtenu $A'=L\,A\,C$ avec $C$ \iv et $L$ strictement carrée \ivz.
Si en outre on a $\fe_i=\fh_i$ pour $i\in\lrbk$, la matrice $C$ est \egmt
strictement carrée.
\eoe

%%%%%%%%%%%%%%%%%%%%%%%%%%%%%%%%%%%%%%%%%%%%%%%%%%%%%%%%%%%%%%%%%%%%

%\penalty-2500
\section[Réduction de Hermite]{Domaines de Prüfer, \slisz, réduction de Hermite}

\CAdre{.8}{Dans cette section, $\gZ$ est un \ddp à \dve explicite, avec $\gK$ pour corps de fractions.}

%: subsec{Systèmes linéaires, généralités}
\subsec{Systèmes linéaires, généralités}

On considère un \sli $ A  X = B$ de $n$ \eqns à $m$ inconnues sur $\gZ$.
Dans la situation usuelle, tous les \ids indexant les lignes et les colonnes
des matrices $A$, $B$ (données) et $X$ (inconnue) sont égaux à $\gen{1}$, les inconnues (\coos de $X$) sont dans $\gZ$ \hbox{et $\und A\in\MM_{n,m}(\gZ)$}.
Plus \gnltz, on considère des \pmas $A=(\fh_1,\dots,\fh_n;\fe_1,\dots,\fe_m;\und A)$ et~$B=(\und \fh; \fb;\und B)$, et la matrice inconnue $X$ est dans
$\MM_{\und \fe;\fb}(\gZ)$.

Comme on le voit sur la \demz, le premier résultat vaut sur un anneau intègre arbitraire.
%:     proposition  propplg
\begin{proposition} \label{propplg} \emph{(Principe local-global)}
Soient $s_1$, \dots, $s_\ell$ \com dans $\gZ$. Le \sli $AX=B$ admet une solution $X\in\MM_{\und \fe;\fb}(\gZ)$ \ssi il admet une solution dans chaque $\MM_{\und \fe;\fb}(\gZ_{s_i})$. 
\end{proposition}
\begin{proof}
La condition est évidemment \ncrz. Soit $X_i=\fraC1{s_i^k}Y_i$ une solution dans $\MM_{\und \fe;\fb}(\gZ_{s_i})$ avec $Y_i\in\MM_{\und \fe;\fb}(\gZ)$, de sorte que $AY_i=s_i^k B$. On écrit $\sum_{i=1}^\ell u_is_i^k=1$ avec les $u_i\in \gZ$. Alors $X=\sum_{i=1}^\ell u_iY_i$ est une solution dans $\MM_{\und \fe;\fb}(\gZ)$.
\end{proof}

Le point \emph{1}  du \tho suivant, bien que cas très particulier, est intéressant car le calcul de la solution se fait dans $\gK$.
Le point \emph{2} reprend dans le langage des \pmas le point~\emph{3} du \tho XII-3.2 dans [ACMC]. 

%:     theorem  thsliDDP
\begin{theorem} \label{thsliDDP} \emph{(Systèmes linéaires sur un \ddpz, I)}
\begin{enumerate}
\item Cas particulier: on donne une \pma carrée avec \hbox{$\fde(A)\neq \gen{0}$} et l'on suppose que les \ids mineurs d'ordre $1$ du vecteur colonne $B$
sont tous contenus dans $\fde(A)$. Alors la solution $\und X$ calculée dans $\gK$ fournit une pseudomatrice $X\in\MM_{\und \fe,\fb}(\gZ)$, autrement dit $x_j\fb\in \fe_j$ pour $j\in\lrbn$.
\item Cas \gnlz: le \sli admet une solution \ssi les \idds de $[\,A\mid B\,]$ sont égaux à ceux de $A$.
\end{enumerate}
\end{theorem}
\begin{proof} Le point \emph{1} peut se démontrer par calcul direct en utilisant la règle de Cramer. Il est clair aussi qu’il résulte du point \emph{2}.  \\
\emph{2}.
Considérons des \eco $s_i$ de $\gZ$ tels que, après \lon en chacun des $s_i$, chacun des \ids $\fh_j$, $\fe_i$ et $\fb$ est principal, \cad libre de rang $1$. Ceci nous ramène au cas des matrices usuelles, pour lesquelles l'\eqvc est établie dans [ACMC, XII-3.2].
\end{proof}
%

%: subsec{Manipulations du type Bezout}
\subsec{Manipulations du type Bezout} 

Les manipulations de colonnes (resp. de lignes) de type Bezout que nous examinons correspondent à des changements de \pba pour le module de départ (resp. d'arrivée).

On désire faire appara\^{\i}tre un \coe nul de la même fa{\c c}on que ce que l'on fait dans le cas des domaines de Bezout avec des matrices usuelles.

\smallskip 
\emph{On commence par un cas particulier.}
\\
On traite la \pma particulière suivante de format $1\times 2$ 
$$\preskip.1em \postskip.4em 
L=\;\;\;\bordercmatrix [\lbrack\rbrack]{
       & \fa & \fb    \cr
\fc &   a    &  a 
}. 
$$
Plut\^ot que de donner tout de suite la solution, expliquons d'o\`u elle vient. Cela nous prend une quinzaine de lignes, mais cela éclaire la chose.
On a une suite exacte scindable

\snic{0\to\fa\cap\fb\vvvvers{\tra[\,-1\;1\,]~} \fa\oplus\fb\vvvers{[\,1\;1\,]~} \fa+\fb\to 0}

qui donne un \iso $(\fa+\fb)\oplus(\fa\cap\fb)\vers{\psi} \fa\oplus\fb$.
\\
Plus \prmtz, regardons $\fa\oplus\fb$ comme le sous-\Zmo $E$ de~$\gK^{2}$
muni de la \pba $\cE=\big((e_1,\fa),(e_2,\fb)\big)$ ($(e_1,e_2)$, base canonique de~$\gK^{2}$), i.e. 
$$\preskip.4em \postskip.4em 
E=e_1\fa\oplus e_2\fb=E_1\oplus E_2. 
$$
Il nous faut maintenant trouver une \pba $\cH=\big((f_1,\fa+\fb),(f_2,\fa\cap\fb)\big)$ de $E$, \cad une matrice $B\in\MM_2(\gK)$ convenable,
dont les colonnes exprimeront $f_1$ et $f_2$ sur la base $(e_1,e_2)$. 
%En prenant $f_2=-e_1+e_2$ \hbox{l'\id $\fa\cap\fb$} peut être vu comme un sous-module~$F_2=f_2(\fa\cap\fb)$
%de~$M$.
%\hbox{Et $\fa+\fb$} est isomorphe à un \sul $F_1=f_1(\fa+\fb)$ de $F_2$ 
%dans~$M$.
% Il s'agit maintenant de rendre explicite le \gui{changement de base} 
% qui fait passer de $E_1\oplus E_2$ à $F_1\oplus F_2$ pour le module~$M$.
% Ceci revient à calculer la matrice de l'\iso $\psi$.

 Pour cela, on utilise l'\egt $(\fa:\fb)_\gZ+(\fb:\fa)_\gZ=\gen{1}$. On a donc 
$$\preskip.4em \postskip.4em 
s\in (\fb:\fa)_\gZ\;\hbox{ et }\;t\in (\fa:\fb)_\gZ\;\hbox{ avec } \;s+t=1, 
$$
d'o\`u $s(\fa+\fb)\subseteq \fb$,  
$t(\fa+\fb)\subseteq \fa$ et $s\fa+t \fb= \fa\cap\fb$. 
D'o\`u la matrice de passage de $\cE$ à $\cH$ avec $f_1=te_1+se_2$ et $f_2=-e_1+e_2$.
$$\preskip.4em \postskip.3em
B=\scM_{\cH,\cE}(\Id_M)=\;\;\bordercmatrix [\lbrack\rbrack]{
       & \fa+\fb & \fa\cap\fb    \cr
\fa &   t    &  -1 \cr
\fb &   s    & 1  \cr
},
$$
avec pour inverse
$$\preskip-.3em \postskip.4em
{\scM_{\cE,\cH}(\Id_M)=\;\;\bordercmatrix [\lbrack\rbrack]{
              & \fa & \fb    \cr
\fa+\fb &   1    &  1 \cr
\fa\cap\fb &  -s    & t \cr
}.}
$$
Notons que l'\egt $\fa\,\fb=(\fa+\fb)(\fa\cap\fb)$ nous assurait que 
la première matrice était celle d'un \isoz, donc il n'y a aucune surprise 
à voir que la matrice inverse est bien celle d'une \Zliz.
\\
En résumé $B$ est une \pma \iv et l'on obtient le résultat souhaité:
$$\preskip-.1em \postskip.4em 
LB= \;\;\bordercmatrix [\lbrack\rbrack]{
              & \fa+\fb& \fa\cap\fb   \cr
\fc &   a    &  0 
} 
$$

%%%%%%%%%%%%%%%%%%%%%%%%%%%%%%%%%%%%%%%%%%%%%%%%%%%%%%%%%%%%%%%%%%%%%%%%%%%
\smallskip  
\emph{Le cas \gnlz.} 
\\
On veut maintenant traiter la \pma \gnle de format $1\times 2$ 
$$\preskip.2em \postskip.4em 
L=\;\;\bordercmatrix [\lbrack\rbrack]{
       & \fa & \fb    \cr
\fc &   a    &  b 
}. 
$$
On se demande d'abord si le pivot de Gauss fonctionne (\cad est-ce que $a\fa\supseteq b\fb$ ou $a\fa\subseteq b\fb$?), auquel cas on fait appel au pivot de Gauss. 
\\
Si ce n'est pas le cas, il suffit de se ramener au cas particulier précédent, en considérant les \ids $\fa'=b^{-1} \,\fa$ et $\fb'= a^{-1} \,\fb$.% 
\footnote{Dans le cas où l'on mène tous les calculs avec des \ids entiers, on considère
 un \elt $x\in\gK$ tels que les \ifrs $\fa'=x b^{-1} \,\fa$ et~$\fb'=x a^{-1} \,\fb$ soient entiers, i.e. $x\in b\fa^{-1}\cap a\fb^{-1}$. Simplement,  on peut prendre pour $x$  un multiple commun des numérateurs de $a$ et~$b$, ce qui donne $x b^{-1}$ et $x a^{-1}\in\gA$, et a fortiori  $\fa'$ et~$\fb'$
 entiers.}
On a alors une \pma \iv
$$\preskip-.2em \postskip.4em 
P=\;\;\bordercmatrix [\lbrack\rbrack]{
    & \fa' & \fb'    \cr
\fa &  b    &  0 \cr
\fb &  0    &  a  \cr
} \;\;\hbox{ avec }\; LP=\;\;\bordercmatrix [\lbrack\rbrack]{
       & \fa' & \fb'    \cr
\fc &   ab    &  ab 
}\,. 
$$
En fin de compte ce sont donc les \ids $\fa'$ et $\fb'$ qui seront traités
selon la méthode donnée dans le cas particulier étudié en premier.
En particulier, on prendra $s$ et $t\in\gA$ tels que $s+t=1$, $s\fa'\subseteq \fb'$  et 
$t\fb'\subseteq \fa'$ (i.e. $sa\fa\subseteq b\fb$ \hbox{et $tb\fb\subseteq a\fa$}).

%%%%%%%%%%%%%%%%%%%%%%%%%%%%%%%%%%%%%%%%%%%%%%%%%%%%%%%%%%%%%%%%%%%
%:subsec{Réduction de Hermite sur les colonnes}
\subsec{Réduction de Hermite: images et noyaux d'\alisz}

Dans cette sous-section, on montre comment la réduction de Hermite des \pmas permet de démontrer sous forme très concrète les \thos de structure \gnls concernant aussi bien les \slis que  les \mpfs sur les \ddpsz.  

%:     theorem  thRedHerDDP
\begin{theorem} \label{thRedHerDDP} \emph{(Réduction de Hermite)}\\
Soit $\gZ$ un \ddp de corps de fractions $\gK$
pour lequel on a un test de \dve et un \algo d'inversion
des \itfsz.
\\
Soit $A=((\fh_i)_{i\in\lrbn}),(\fe_j)_{j\in\lrbm}),\und A)$ avec $\und A=(a_{ij})\in\MM_{n,m}(\gK)$ une \pma
sur $\gZ$ (i.e.,  $a_{ij}\fe_j\subseteq \fh_i$ pour tous $i,j$).
\\
On peut calculer une \pma \iv 
$$\preskip.4em \postskip.4em 
C=((\fe_j)_{j\in\lrbm},(\fe'_j)_{j\in\lrbm},\und C)\hbox{  avec  }\und C=(c_{j,k})\in\Mm(\gK) 
$$
(on regarde $C$ comme opérant un \cpb au départ), et une matrice de permutation $\und L\in\Mn(\gZ)$ (la \pma $L$ correspondante opère un \cpb à l'arrivée), telles que la matrice $L\,A\,C$ soit en forme réduite de Hermite:
\smashbot{$H=\blocs{.7}{.9}{.7}{.5}{$T$}{$0$}{$G$}{$0$}\,$,} $T$ carrée triangulaire inférieure de \deter
 non nul. 
\end{theorem}

\begin{proof}
Cela résulte des manipulations de colonnes de type pivot de Gauss et de type Bezout décrites précédemment. Tout se passe comme pour une matrice sur un anneau de Bezout intègre, à ceci près que l'on doit remplacer les matrices usuelles par des \pmasz.
\end{proof}

Notons que le test de \dve est utilisé uniquement pour vérifier
que dans certains cas, des \mnrs du type pivot de Gauss sont admissibles.

On pourrait donc se passer de cette hypothèse en faisant fonctionner l'\algo de réduction de Hermite uniquement avec des manipulations de type Bezout.
 
Cependant, comme déjà indiqué, nous faisons toujours l'hypothèse que
l'anneau $\gZ$ est un \ddp à \dve explicite (au sens \cofz).

Les \thos qui suivent sont tous des corollaires faciles du \thref{thRedHerDDP}
concernant les \pmasz.
%:     theorem cor0thRedHerDDP
\begin{theorem} \label{cor0thRedHerDDP} \emph{(Image d'une \pmaz)}
 L'image de toute \pma est  isomorphe à une somme directe d'\ids \ivsz.
\end{theorem}
\begin{proof}
On effectue une réduction de Hermite de la matrice. La multiplication à droite par une \pma \iv ne change pas le module image.  Il suffit de traiter l'image d'une \pma $(\und \fh;\und{\fe'};H)$ o\`u \smashbot{$H=\blocs{.7}{.9}{.7}{.5}{$T$}{$0$}{$G$}{$0$}$} est en forme réduite de Hermite. C'est aussi l'image de la \pma  
\smashbot{$(\und\fh;\fe'_1,\dots,\fe'_r;H_1)$ o\`u
$H_1=\blocs{.7}{0}{.7}{.5}{$T$}{}{$G$}{}$}\,. Comme $\det(T)\neq 0$, l'\ali
$$\preskip.2em \postskip.4em 
\;\;\fe'_1\times \cdots\times \fe'_r\;\lora\;\fh_1\times \cdots\times \fh_n 
$$
définie par $H_1$ est injective. Donc $\Im A=\Im H=\Im H_1\simeq \fe'_1\times \cdots\times \fe'_r$.
\\
\Prmtz, supposons que l'on traite la \pma $A$ d'une \Zli $\varphi:E\to H$ sur des \pbas $\cE=\big((e_1,\fe_1),\dots,(e_m,\fe_m)\big)$ et~$\cH=\big((h_1,\fh_1),\dots,(h_n,\fh_n)\big)$ de $E$ et $H$. Alors l'image de $\varphi$ est le module somme directe $\fe'_1h'_1\oplus \cdots \oplus\fe'_rh'_r$, où les $h'_i$ sont les vecteurs exprimés par les colonnes de la matrice $H_1$ sur les $h_i$ (qui sont des \elts de $\gK\otimes_\gZ H$). 
\end{proof}
%
%:     theorem  cor1thRedHerDDP
\begin{theorem} \label{cor1thRedHerDDP} \emph{(Images)}
\begin{enumerate}
\item Tout \Zmo\ptf est isomorphe à une somme directe
d'\ids \ivsz, autrement dit, il possède une \pbaz.
\item L'image d'une   \ali entre \Zmos\ptfs  est un \mptf dont on peut calculer une \pbaz. 
\item 
\begin{itemize}
\item Tout sous-\Zmo \tf de $\gK^{n}$ est \ptfz.
\item Si $ g_1$, \dots, $g_m\in\gK^{n}$ et si~$\fe_1$, \dots, $\fe_m$ sont des \itfs non nuls, on peut calculer une \pba du \Zmo $E=\sum_j g_j\fe_j$. 
\end{itemize}
\end{enumerate}
\end{theorem}
\begin{proof}
\emph{1.} Puisqu'un \mptf est isomorphe à l'image d'une matrice de \prn usuelle, cela résulte du \thref{cor0thRedHerDDP}.

\emph{2.} Il s'agit pour l'essentiel d'une reformulation du  \thref{cor0thRedHerDDP}.
En effet, vu le point \emph{1}, tout \mptf possède une \pbaz. Alors une \ali entre \mptfs est représentée par une \pma 
sur les \pbas des modules.

\emph{3.} Les deux points sont clairement \eqvsz. Le second résulte
du \thref{cor0thRedHerDDP}  car $E$ est l'image d'une \pma
 $(\gen{d},\dots,\gen{d};\fe_1,\dots,\fe_m;\und G)$, où~$\und G$ a pour colonnes les $g_j$ exprimés sur la base canonique de $\gK^n$, et où
 l'\elt  $d\in\gK$
 est choisi de façon que  $g_{ij}\fe_j\subseteq d \gZ$ pour tous les $(i,j)$. 
\end{proof}
%

%:     theorem  cor2thRedHerDDP
\begin{theorem} \label{cor2thRedHerDDP} \emph{(Noyaux)}\\
Soit $\varphi:E\to H$ une  \ali entre \mptfs sur un \ddpz. Le noyau $\Ker\varphi$ est facteur direct dans~$E$ (a fortiori \ptfz).
Plus \prmtz, on peut calculer une \pba de $E$ dont les derniers termes forment une \pba de~$\Ker\varphi$.
\end{theorem}
\begin{proof}
Comme dans le \tho précédent, il suffit de faire la \dem pour le noyau d'une \pmaz, et l'on prend les notations du \thref{thRedHerDDP}. Le changement de \pba donné par la matrice~$C$
 fournit une nouvelle \pba  $\big((u_1,\fe'_1),\dots,(u_m,\fe'_m)\big)$
(les $u_j$ sont exprimés par les colonnes de la matrice $C$), pour laquelle l'\ali est exprimée par la nouvelle \pma $H$. Pour cette matrice il est clair que le noyau est fourni par la somme directe des
sous-modules $u_j\fe'_j$ pour les numéros $j$ des colonnes nulles de $H$.
\end{proof}
%

%: subsec{Pseudo-matrices surjectives}
\subsec{Pseudo-matrices surjectives}

La proposition suivante reprend le point \emph{3} du \thref{factIddPma}, qui \gns un \tho bien connu pour les matrices usuelles. 

Dans le point \emph{3} du \thref{factIddPma}, l'\ali $\varphi$, ici représentée par~$A$, admet un inverse à droite $\psi$ du fait qu'elle est surjective et que l'image est un module \proz. 
Le module source de $\varphi$ est  égal à $\Im\psi\oplus \Ker\varphi$. Lorsque $\Ker\varphi$ est une somme directe de modules \pros de  rang $1$, la matrice $A$ peut être complétée en la \pma inversible $B$ qui représente l'\ali $\varphi\frt{\Im\psi}\oplus\Id_{\Ker\varphi}$, laquelle a pour source $\Im\psi\oplus \Ker\varphi$ et pour image $\Im\varphi\oplus \Ker\varphi$.

Une fois que l'existence de la \pma $B$ est assurée par un argument \gnlz, il reste le \pb \algq d'un calcul aussi simple que possible de cette \pmaz.

Nous sommes intéressés ici par une \dem matricielle  \algq du résultat
%, éventuellement simplifiée 
dans le contexte des \ddpsz, pour lesquels les calculs sur les \pmas sont simplifiés. 
%Ici on demande un calcul \gui{clair} de l'inverse à droite d'une \pma surjective sous forme d'une \pmaz.

%:     proposition  thpmasurj
\begin{proposition} \label{thpmasurj}
 \emph{(Pseudo-matrices surjectives)} \\
Soit $\gZ$ un \ddp  et soit $n\in\lrbm$.\\
Soit $A=(\fh_1,\dots,\fh_n;\fe_1,\dots,\fe_m;\und A)\in\MM_{\und\fh,\und\fe}(\gZ)$. 
%On choisit $(\fh_{n+1},\dots,\fh_m)$
%de telle sorte que $\bigoplus_{k=1}^m\fe_k\simeq \bigoplus_{k=1}^m%\fh_k$.\\  
Cette \pma représente une \ali surjective \ssi  $\fD_n(A)=\gen{1}$. Dans ce cas, $A$ peut être complétée en une \pma \iv   
$$\preskip.3em \postskip.4em 
B=(\fh_1,\dots,\fh_m;\fe_1,\dots,\fe_m;\und B). 
$$
La matrice $AB^{-1}$ est alors en forme de Smith:
$$\preskip.3em \postskip.4em 
AB^{-1}=\;\;\bordercmatrix [\lbrack\rbrack]{
       & \fh_1 & \fh_2  & \cdots & \fh_n& \fh_{n+1}& \cdots& \fh_m   \cr
\fh_1 &  1      &  0 &  \cdots &  0  &  0   & \cdots &  0   \cr
\fh_2 &  0      &  1 &  \ddots       &  \vdots &  \vdots  &   &  \vdots\cr
~\vdots  &  \vdots &  \ddots   & \ddots &  0 &  \vdots &  &  \vdots \cr
\fh_n &   0 & \cdots  &  0&    1 &  0& \cdots &  0   
}. 
$$
Et la \pma extraite sur les $n$ premières colonnes de $B^{-1}$ est un inverse à droite de $A$. 
%Soit $A=(\fh_1,\dots,\fh_n;\fe_1,\dots,\fe_m;\und A)\in\MM_{\und\fh,\und\fe}(\gZ)$. Cette \pma représente une \ali surjective \ssi  $\fD_n(A)=\gen{1}$. Dans ce cas la \pma $A$ admet un inverse à droite $C=(\fe_1,\dots,\fe_m;\fh_1,\dots,\fh_n;\und c)\in\MM_{\und\fe,\und\fh}(\gZ)$.
\end{proposition}
\begin{proof} On démontre l'implication difficile.\\
D'après le \thref{thRedHerDDP}, on peut calculer une \pma \iv 
$$\preskip.4em \postskip.4em 
C=((\fe_j)_{j\in\lrbm};\fe'_1,\dots,\fe'_n,\fh_{n+1},\dots,\fh_m;\und C)\hbox{  avec  }\und C=(c_{j,k})\in\Mm(\gK) 
$$
%(on regarde $C$ comme opérant un \cpb au départ) 
telle que la matrice $A\,C$ soit en forme réduite de Hermite:
{$H=\blocs{.7}{.9}{.7}{0}{$T$}{$0$}{ }{ }\,$,} où $T$ est carrée triangulaire inférieure. On a $\fD_n(T)=\fD_n(A)=\gen{1}$, donc~$T$
est \ivz. En fait les \ids mineurs d'ordre $1$ diagonaux $\fm_{ii}(T)$ ont leur produit égal à $\gen{1}$, donc sont tous égaux à $\gen{1}$, ce qui autorise la méthode du pivot par \mnrs de colonnes. On termine en inversant la matrice diagonale obtenue. On a ainsi 
ramené~$T$ à la forme $(\fh_1,\dots,\fh_n;\fh_1,\dots,\fh_n;\In)$.  Ceci donne une \pma \iv $C'$ telle que $A\,C\,C'=\blocs{.7}{.9}{.7}{0}{$\In$}{$0$}{ }{ }\,$ et l'on pose $B:=(C\,C')^{-1}$.
\\
En regardant $B$ comme une \pma par blocs 
\smash{$B=\blocs{1.6}{0}{.7}{.9}{$B_1$}{}{$B_2$}{ }\,$}, on obtient
$ 
\blocs{.7}{.9}{.7}{0}{$\In$}{$0$}{ }{ } \;
\blocs{1.6}{0}{.7}{.9}{$B_1$}{}{$B_2$}{ }\, = A$ et $B_1=A$.
%\hum{On demande une \dem purement calculatoire, avec éventuellement une formule magique. Si c'est le cas, la formule magique est sans doute plus \gnle et pourrait être donnée dès la section 2. 
%On mettrait donc cette proposition comme précision juste après le \thref{factIddPma}.
%}
\end{proof}
%

%%%%%%%%%%%%%%%%%%%%%%%%%%%%%%%%%%%%%%%%%%%%%%%%%%%%%%%%%%%%%%%%%%%%
%: subsec{Double réduction de Hermite: sur les lignes et les colonnes}
\subsec{Double réduction de Hermite, conoyaux}

%:     theorem  thRedHerDDP2
\begin{theorem} \label{thRedHerDDP2} \emph{(Double réduction de Hermite)}\\
Mêmes hypothèses qu'au \thref{thRedHerDDP}. %\\
On peut calculer une \pma \iv $C$ 
 opérant un \cpb au départ, et une \pma \iv $L$ opérant un \cpb à l'arrivée, telles que  
\smashbot{$L\,A\,C=\blocs{.7}{.9}{.7}{.5}{$T$}{$0$}{$0$}{$0$}$\,,} où $T$ est  triangulaire inférieure
de \deter non nul. 
\end{theorem}
\begin{proof}
On fait une réduction de Hermite sur les colonnes, suivie d'une réduction
de Hermite sur les lignes.
\end{proof}

Le \tho suivant est une forme précise, \gui{pseudo-matricielle}, du \tho XIV-2.5 de [Modules].
%:     theorem  cor3thRedHerDDP
\begin{theorem} \label{cor3thRedHerDDP} \emph{(Conoyaux)}\\
Soit $\varphi:E\to H$ une  \ali entre \mptfs sur un \ddpz. 
\begin{enumerate}
\item Le saturé de $\Im\varphi$ dans $H$ est en facteur direct dans $H$. Plus \prmtz, on peut calculer une \pba de $H$ dont les premiers termes forment une \pba du saturé de $\Im\varphi$.
\item Le module conoyau $\Coker\varphi$ est un \mpfz, somme directe de son sous-module de torsion et d'un \mptfz.
\item Le sous-module de torsion de $\Coker\varphi$ est le quotient d'un module \pro par un sous-module \pro de même rang. Il est isomorphe au conoyau d'une \pma
carrée injective.
\end{enumerate}
En particulier tout \mpf sans torsion est \proz.
\end{theorem}
\begin{proof}
Le point \emph{1} résulte du \thref{thRedHerDDP2} et du fait que lorsque l'on sature l'image d'une \pma carrée injective, on obtient le module d'arrivée en entier. Le reste suit facilement.
\end{proof}
\rem Dans le \thref{cor3thRedHerDDP}, la structure du sous-module de torsion et celle du sous-\mptf ne sont pas décryptées de manière complètement satisfaisante. En particulier, ce \thoz, bien que \cofz, ne donne pas de solution \algq générale aux deux \pbs suivants:
\begin{itemize}
\item déterminer si deux \mpfs de torsion sont ou ne sont pas isomorphes;
\item déterminer si deux \mptfs sont ou ne sont pas isomorphes. 
\end{itemize}

\smallskip 
On verra dans la section \ref{CasDim=1} que lorsque le \ddp $\gZ$ est de dimension $1$, on sait répondre à la première question, mais pas \ncrt à la seconde,
qui se ramènera à la question de savoir si un \itf $\fe$ est principal. Ce \pb n'admet pas de solution \algq \gnlez, même pour les \doksz.
\eoe

%%%%%%%%%%%%%%%%%%%%%%%%%%%%%%%%%%%%%%%%%%%%%%%%%%%%%%%%%%%%%%%%%%%%%%%%%%%
%:subsec{Systèmes linéaires}
\subsec{Systèmes linéaires, via la réduction de Hermite}

Le \thref{thRedHerDDP} donne le corollaire suivant pour le traitement des \slisz.

On notera que les tests proposés ici sont explicites parce que le \dve dans~$\gZ$ est explicite.

%:     theorem  cor4thRedHerDDP
\begin{theorem} \label{cor4thRedHerDDP} \emph{(Systèmes linéaires sur un \ddpz, II)}
\\
On considère un \sli $ A  X = B$ de $n$ \eqns à $m$ inconnues sur $\gZ$.
Dans la situation usuelle, tous les \ids indexant les lignes et les colonnes
des matrices $A$, $B$ (données) et $X$ (inconnue) sont égaux à $\gen{1}$, les inconnues (\coos de $X$) sont dans $\gZ$ \hbox{et $\und A\in\MM_{n,m}(\gZ)$}.
\\
Les résultats s'appliquent aussi bien avec des \pmas $A=(\und\fh;\und\fe;\und A)$ et $B=(\und \fh; \fb;\und B)$, en adaptant le format de la colonne  $X$ aux données: les inconnues sont des \elts $x_j\in\gK$ soumis
aux conditions $x_j\fb\in \fe_j$. \\
On considère une double réduction de Hermite $L\,A\,C=\blocs{.7}{.9}{.7}{.5}{$T$}{$0$}{$0$}{$0$}$, 
avec la matrice $\und T\in\MM_r(\gK)$ triangulaire inférieure injective. Les \pmas $C=(\und\fe,\und{\fe'},\und C)$ et $L=(\und{\fh'},\und{\fh},\und L)$ sont \ivsz, et l'on introduit les nouvelles inconnues $Y=C^{-1}X$.
\begin{enumerate}
\item La solution \gnle du \sys \hmg sans second membre $AX=0$ est donnée par $y_1=\cdots=y_r=0$. Ceci donne pour $Y$ le sous-\Zmo 
$$\preskip.4em \postskip.4em 
K=\so 0 \times \cdots\times \so 0\times \fe'_{r+1}\times \cdots\times \fe'_m 
$$
et pour $X$, le \Zmo $C(K)$ paramétré par 
$$\preskip.3em \postskip.3em 
(y_{r+1},\dots,y_m)\in \fe'_{r+1}\times \cdots\times \fe'_m, 
$$
c'est aussi l'image de la \pma extraite de $C$ sur les $m-r$ dernières colonnes.
\item Le \sli avec second membre $AX=B$ admet une solution \ssi 
le sous-\Zmo de $\gK^{n}$  codé par $LB$
est contenu dans l'image de la \pma 
$$\preskip.0em \postskip.4em 
\Big(\,\fh'_1,\dots,\fh'_n;\fe'_1,\dots,\fe'_r;
\Cmatrix{3pt}{T\cr0}\Big)
. 
$$
\begin{itemize}
\item Les contraintes données par l'annulation des $n-r$  dernières \coos de $\und L\,\und B$  signifient exactement que le \sys $\und A\, \und X=\und B$ a une solution dans $\gK^{m}$. 
\item  Les contraintes sur les $r$ premières \coos de $\und L\,\und B$ peuvent être testées de proche en proche, vue la forme triangulaire de $\und T$. 
\item On peut donc tester l'existence d'une solution et calculer une solution particulière lorsqu'il en existe une. 
\end{itemize}
%
%\item 
%
\end{enumerate}
\end{theorem}
\begin{proof}
Il faut juste préciser la question du test
pour une solution particulière: toute solution donne par soustraction une solution avec les \gui{inconnues auxiliaires} nulles, donc le test peut porter uniquement sur l'équation $TZ=B'$ avec $Z$ et $B'$ donnés par troncatures de~$Y$ et 
$LB$ aux $r$ premières \coosz. Comme $T$ est injective, la solution, si elle existe, est unique.
\end{proof}
\comm L'énoncé du corollaire sur les \slis est un peu plus simple dans le cas d'un \sli usuel, mais il n'est pas \gui{nettement plus simple}, car les 
\ids $\fe'_j$ et $\fh'_i$ introduits par la double réduction de Hermite,
ainsi que les \coes de $\und T$ dans $\gK$ et non dans $\gZ$, semblent des
ingrédients inévitables.
\eoe

%%%%%%%%%%%%%%%%%%%%%%%%%%%%%%%%%%%%%%%%%%%%%%%%%%%%%%%%%%%%%%%%%%%%%%%%%
%\newpage
\section[Réduction de Smith en dimension 1]{Calculs modulaires et réduction de Smith  en dimension $1$} \label{CasDim=1}

\smallskip 
\CAdre{.8}{Dans cette section, $\gZ$ est un \ddp à \dve explicite,
de dimension $\leq 1$.}

%:     definition  defiSmithpma
\begin{definition} \label{defiSmithpma}
Soit  $A=(\fh_1,\dots,\fh_n;\fe_1,\dots,\fe_m;\und A)\in\MM_{\und\fh,\und\fe}(\gZ)$ et $p=\inf(m,n)$.
La \pma $A$ est dite \emdf{en forme de Smith} si la matrice $\und A=(a_{ij})$ a tous ses \coes nuls hormis éventuellement des \coes $a_{ii}$
et si les \ids mineurs d'ordre $1$ $\fc_i:=\fm_{ii}(A)=a_{ii}\fe_i\fh_i^{-1}$ satisfont les relations attendues:
$$\preskip.4em \postskip.4em 
\fc_1\supseteq \fc_2\supseteq \cdots\supseteq \fc_p.    
$$
\end{definition}

%  theorem thSmithDedekind
\begin{theorem} \label{thSmithDedekind}
Soit $\gZ$ un \dok explicite, i.e., un \ddp \noe à \dve  explicite.
Toute \pma sur $\gZ$ est \eqve à une \pma en forme de Smith. 
\end{theorem}
\begin{proof}
L'\algo est le même que pour la réduction de Smith d'une matrice sur un anneau principal, en remplaçant les matrices usuelles par des \pmasz. On peut se reporter à [Modules, Algorithme IV-2.2].  
\end{proof}

%: subsec{Sur l'unicité d'une réduite de Smith}
\subsec{Sur l'unicité d'une réduite de Smith usuelle}

Lorsqu'une matrice $M\in\MM_{n,m}(\gA)$ admet une réduction de Smith
dont la partie carrée est égale à $\Diag(a_1,\dots,a_n)$, le module  $\Coker M$ est isomorphe à la somme directe des $\aqo\gA{a_i}$:
$$\preskip.4em \postskip.4em
\Coker M\;\simeq\;\aqo\gA {a_1}\oplus\cdots\oplus \aqo\gA {a_n}  \; \hbox{ avec des \ids  } \; 
 \gen{a_1}\supseteq \cdots\supseteq \gen{a_n}.
$$
Les \ids $\fe_i=\gen{a_i}$ sont uniquement déterminés (fait \ref{factMV-9.1}). 

Plus \prmtz, si aucun des $a_i$ n'est \ivz, la liste
$(\fe_1,\dots,\fe_n)$ ne dépend que du module $\Coker M$. 

Cela ne signifie pas \ncrt que pour deux matrices \eqves en forme de Smith $A=(a_{ij})$ et $A'=(a'_{ij})$ les  $a'_{ii}$ soient associés aux $a_{ii}$.

Néanmoins, c'est vrai lorsque $\gA$ est intègre, et le lemme suivant montre que c'est aussi le cas pour un anneau \zedz.     

%:     lemma  lemzedassoc
\begin{lemma} \label{lemzedassoc}
Soient $a$ et $b\in\gA$  \zedz. 
Si $\gen{a}=\gen{b}$, alors $a$ et~$b$ sont \emdf{associés}: $a=wb$
pour un \elt $w\in\Ati$.
\end{lemma}
\begin{proof}
On a $a=ub$ et $b=va$, donc $a(uv-1)=0$.
Comme $\gA$ est \zed on a un \idm $e$ tel que $uv-1$ est nilpotent 
modulo~$e$ et inversible modulo $1-e$. Dans la composante $\aqo\gA e$, on a $uv=1+x$ avec $x$  nilpotent, donc $u$ est \ivz. Dans la composante $\aqo\gA {1-e}$, $a=0$ donc $b=0$ et l'on a $a=1.b$. En bref $a=wb$ avec $w=u.(1-e)+1.e$ \ivz.
\end{proof}

On montre maintenant un lemme qui complète en quelque sorte le lemme précédent.

%:     lemma  lemzedsmithass
\begin{lemma} \label{lemzedsmithass}
Soit $\Diag(\an)$ une matrice en forme de Smith sur un anneau $\gA$, et des $u_i\in\Ati$ tels que $\prod_iu_i=1$. Définissons $a'_i=a_iu_i$, alors les deux matrices $\Diag(\an)$ et $\Diag(a'_1,\dots,a'_n)$ sont \elrt \eqvesz.
\end{lemma}
\begin{proof}
Il suffit de traiter le cas suivant, avec $b=ca$, $a'=au$, $uv=1$, $b'=bv$,
$A=\cmatrix{a&0\\0&b}$ et $B=\cmatrix{a'&0\\0&b'}$.
$$\preskip.4em \postskip.4em \mathrigid 2mu 
\Cmatrix{2.5pt}{a&0\\0&b} \vers 1 \Cmatrix{2.5pt}{a&a'-a\\0&b}
\vers 2 \Cmatrix{2.5pt}{a'&a'-a\\b&b}\vers 3 \Cmatrix{2.5pt}{a'&a'-a\\0&b'}\vers 4 \Cmatrix{2.5pt}{a'&0\\0&b'}
$$
$1$:  $C_2\aff C_2+ (u-1) C_1$.\\
$2$: $C_1\aff C_1+C_2$.\\
$3$: $L_2\aff L_2 - cvL_1$, avec $b-cv(a'-a)=b-ca+cva=cva=vb=b'$.\\
$4$:  $C_2\aff C_2+ (v-1) C_1$.\\
En bref $\Cmatrix{2.5pt}{1&0\\-cv&1}A\Cmatrix{2.5pt}{1&u-1\\0&1}\Cmatrix{2.5pt}{1&0\\1&1}\Cmatrix{2.5pt}{1&v-1\\0&1}=B$.
\end{proof}
%

%\newpage
%: subsec{La structure des \mpfs de torsion, et quelques conséquences matricielles}
\subsec{La structure des \mpfs de torsion, et quelques calculs matriciels}

Rappelons la \dem de la proposition suivante [Modules, XVI-4.6].

%:     proposition{propMpftorsionPrdim1}
\begin{proposition} \label{propMpftorsionPrdim1}
Un \mpf de torsion $P$ sur un \ddp  $\gZ$ de dimension $1$
est isomorphe \`a une somme directe
$$\preskip.4em \postskip.4em
\gZ/ \fa_1\oplus\cdots\oplus \gZ/ \fa_n,  \; \hbox{ avec des \ids \ivs } \; 
 \fa_1\supseteq \cdots\supseteq \fa_n.
$$
\end{proposition}
%--------- fin corollary ----------------------------------------------  
Si les $\fa_i$ sont $\neq \gen{1}$ cette \'ecriture est unique\footnote{Dans la mesure où nous supposons que $\gZ$ est à \dve explicite, on pourra  se débarrasser des \ids égaux à $\gen{1}$.}. Plus \gnltz, deux écritures de ce type de même longueur sont identiques (fait \ref{factMV-9.1}). 
\begin{proof}
Le module $P$ est annul\'e par un \elt $\delta \neq 0$, parce qu'il est de torsion et  \tf (il suffit que $\delta $ annule chacun des \gtrsz). 

L'anneau quotient $\gZ'=\aqo\gZ \delta $ est \zed et \ariz, donc est un anneau de Smith ([Modules, XVI-4.2]).

Ainsi, on a un \iso de  \Zmos
(ce sont aussi \hbox{des $\gZ'$-modules}) du type suivant 
$$\preskip.3em \postskip.3em
P=P/\delta P\simeq \aqo{\gZ'}{\ov{a_1}}\oplus \cdots\oplus \aqo{\gZ'}{\ov{a_n}}
\simeq \aqo{\gZ}{\delta ,a_1}\oplus \cdots\oplus \aqo{\gZ}{\delta ,a_n},
$$
avec les inclusions 
${\gen{\delta ,a_1}\supseteq \cdots\supseteq \gen{\delta ,a_n}.}$

\emph{En pratique, cela se passe \prmt comme suit}. Le module $P$ est présenté par une matrice $M\in\MM_{n,m}(\gZ)$ avec $\fd_n:=\fD_n(M)=\fF_0(P)\neq \gen{0}$. 
\\
Pour un $\delta\neq 0$ dans $\fd_n$, on calcule sur $\gZ'$ une réduite de Smith par \mnrs 
$\ov L\,\ov M\, \ov C=\blocs{.7}{.9}{.7}{.0}{$D$}{$0$}{ }{ }$ avec
$D=\Diag(\ov{a_1},\dots,\ov{a_n})$
 et~$\gen{\ov{a_1}}\supseteq \cdots\supseteq \gen{\ov{a_n}}$. 
 On en déduit la structure de $P$ comme ci-dessus, de sorte
 que $\fD_k(M)=\fF_{n-k}(P)=\gen{a_1\cdots a_k,\delta}=\fa_1\cdots\fa_k$
 et $\fa_k=\gen{a_k,\delta}$ pour chaque $k\in\lrb{0..n}$. 
\end{proof}

On en déduit le corollaire suivant.

%:     corollary  corpropMpftorsionPrdim1
\begin{corollary} \label{corpropMpftorsionPrdim1}
Sur un  \ddp  de dimension $\leq 1$, la structure d'un \mpf de torsion $P$ est
caractérisée par ses \idfs $\fF_k(P)$. 
\end{corollary}
\begin{proof}
En effet, dans la proposition \ref{propMpftorsionPrdim1}, notons $\fd_i=\fF_{n-i}(P)$. On a alors  $\fd_{1}=\fa_1$,
et plus \gnlt $\fd_{k}=\fa_1\cdots \fa_{k}$  pour $k\in\lrb{2..n}$ (fait \ref{factidfSmith}); donc $\fa_{k}=(\fd_{k}:\fd_{k-1})$ car 
$\fa_{k}\neq \gen{0}$.
\end{proof}

Le corollaire suivant est plus surprenant.

%:     corollary  cor2propMpftorsionPrdim1
\begin{corollary} \label{cor2propMpftorsionPrdim1}
Soit $M\in\MM_{n,m}(\gZ)$ avec $\fD_n(M)\neq 0$. Soit $\delta\neq 0$ dans $\fD_n(M)$.  On considère l'\anar \zed $\gZ'=\aqo\gZ \delta$. On note $\ov x$ pour $x\mod \delta$.  On calcule par \mnrs %(vues sur $\gZ$) 
une réduite de Smith de $\ov M$, disons 
$$\preskip.4em \postskip.4em
\ov{M'}\equiv\blocs{.7}{.9}{.7}{.0}{$D$}{$0$}{ }{ } \mod \delta  \hbox{ avec }D=\Diag(\ov{a_1},\dots,\ov{a_n})\hbox{ et }\gen{\ov{a_1}}\supseteq \cdots\supseteq \gen{\ov{a_n}}.$$
Alors:
\begin{enumerate}
\item $\gen{a_1\cdots a_i,\delta}=\gen{a_1,\delta}\cdots\gen{a_i,\delta}$ pour $i\in\lrb{2..n}$.
\item $\fD_n(M)=\gen{a_1,\delta}\cdots\gen{a_n,\delta}$.
\item  
{\mathrigid 2mu 
$\big(\geN{\prod_{i\in I\cup J} a_i,\delta}:\geN{\prod_{i\in I} a_i,\delta}\big)=\geN{\prod_{i\in J} a_i,\delta}$} dans les cas suivants:\\[.3em]
$I=\lrbk$ et $J=\lrb{k+1..\ell}\subseteq \lrbn$, ou vice versa.
\item  
$\big(\prod_{i\in I\cup J} \ov{a_i}:\prod_{i\in I} \ov{a_i}\big)=\prod_{i\in J} \ov{a_i}$ dans les cas suivants:\\[.3em]
$I=\lrbk$ et $J=\lrb{k+1..\ell}\subseteq \lrbn$, ou vice versa.
\end{enumerate}
 
\end{corollary}
\begin{proof}  \emph{1} et \emph{2}. Cas particulier de ce qui a été expliqué dans la \dem de la proposition \ref{propMpftorsionPrdim1}.
\\
Le point \emph{3}  résulte des \egts $\fd_{i}=\fa_1\cdots \fa_{i}$, car $\gZ$ est un \ddpz.
\\
Le point  \emph{4} en résulte parce 
que lorsqu'on a trois \ids
$I\subseteq J\subseteq K$ dans un anneau $\gA$, on a l'\egt
$$\preskip.4em \postskip.4em 
(J:K)_{\gA} \mod I = (J/I:K/I)_{\gA/I}. 
$$ 
\vspace{-2em}

\end{proof}

\hum{apparemment, les points \emph{1}, \emph{3} et \emph{4}
ne sont pas utilisés dans la suite}

Et maintenant, oh surprise!, on obtient la \gui{réduction de Smith avec \pbasz}
d'une matrice carrée de \deter non nul  sans effort. 

%:     proposition  proppmaSmith
\begin{proposition} \label{proppmaSmith}
Soit $M\in\Mn(\gZ)$ une matrice avec  $d=\det(M)\neq 0$  ($\gZ$ \ddp de dimension $1$). Cette matrice est \eqve à une \pma 
de la forme
$$\preskip.0em \postskip.4em 
M_1=\;\;\bordercmatrix [\lbrack\rbrack]{
       & \fa_1 & \fa_2  & \cdots & \fa_n   \cr
\gen{1} &  1      &  0 &  \cdots &  0      \cr
\gen{1} &  0      &  1 &  \ddots       &  \vdots \cr
~\vdots  &  \vdots &  \ddots   & \ddots &  0 \cr
\gen{1} &   0 & \cdots  &  0&    1    
} 
$$
avec $\fa_1\supseteq \fa_2\supseteq \cdots\supseteq \fa_n $,
$\prod_i\fa_i=\gen{d}$
et chaque $\fa_i$ de la forme  $\gen{a_i,d}$.
\end{proposition}
\begin{proof} 
On reprend ce qui a été expliqué dans la \dem de la proposition \ref{propMpftorsionPrdim1}  \hbox{(éventuellement avec $\delta=d$)}.
En examinant ce qui se passe dans le calcul de la réduction de Smith modulaire par \mnrsz, on obtient sur $\gZ$ une \egt et une congruence
$$\preskip.4em \postskip.4em 
LMC=M'\equiv\Diag(a_1,\dots,a_n) \mod \delta\;. 
$$
Considérons la \pmaz\footnote{On vérifie \imdt que c'est bien une \pmaz.}  $E:=\bordercmatrix [\lbrack\rbrack]{
       & \gen{1} &    \cdots & \gen{1}   \cr
\fa_1 &     &      &         \cr
~\vdots  &    &  M'   &  \cr
\fa_n &    &   &        
}$.
%(\gen{1},\dots,\gen{1};\fa_1,\dots,\fa_n;{M'})$. 
\\
Puisque $M'=\Diag(a_1,\dots,a_n) \mod \delta$, elle est de la forme
$$\preskip.3em \postskip.4em
E\;\;= 
\;\;\bordercmatrix [\lbrack\rbrack]{
       & \gen{1} & \gen{1}  & \cdots & \gen{1}   \cr
\fa_1 &  a_1+\delta x_{11}  &  \delta x_{12} &  \cdots &  \delta x_{1n}      \cr
\fa_2 &  \delta x_{21}      &  a_2+\delta x_{22} &          &  \vdots \cr
~\vdots  &  \vdots &     &  \ddots &  \vdots  \cr
\fa_n &  \delta x_{n1} & \cdots  & \phantom{\ddots}  &    a_n+\delta x_{nn}    
}. 
$$
Oh merveille! c'est une \pma de changement de \pbaz.
En effet, $\prod_{i=1}^n\fa_i=\gen{d}$ avec $d=\det(M)=\det(LMC)=\det(M')=\det(\und E)$, donc $\fde(E)=\gen{1}$.
Son inverse est la \pma
$$\preskip.1em \postskip.4em 
E^{-1}=\;\;\bordercmatrix [\lbrack\rbrack]{
       & \fa_1 &  \cdots & \fa_n   \cr
\gen{1} &         &    &         \cr
~\vdots  &    &  ~{M'}^{-1}   &   \cr
\gen{1} & \phantom{\ddots}    & \phantom{{2}^{2^{2^{2}}}}   &       
}  
$$
Posons $C_1=CE^{-1}$. On obtient $LMC_1=M_1$ de la forme annoncée, où $L$ et $C_1$ sont des \pmas de changements de base ($L$ et $M$ sont des matrices usuelles vues comme des \pmas ayant les \ids $\gen{1}$ en indices de lignes et en colonnes).
\end{proof}

Voici maintenant le résultat analogue pour une matrice
\gui{usuelle} plus large que haute. Il nous faut faire ici appel à la proposition \ref{thpmasurj} concernant les \pmas surjectives en espérant qu'elle soit efficace sur un \ddp de dimension $1$.

%:     theorem  thpmaSmith
\begin{theorem} \label{thpmaSmith} 
Soit $M\in\MM_{n,m}(\gZ)$ une matrice avec  $\fD_n(M)\neq 0$.
Soit un \elt $\delta$ non nul de $\fD_n(M)$. 
Cette matrice est \eqve à une \pma  de la forme
$$%\preskip.0em \postskip.4em 
M_1=\;\;\bordercmatrix [\lbrack\rbrack]{
       & \fa_1 & \fa_2  & \cdots & \fa_n&\cdots & \fa_m   \cr
\gen{1} &  1      &  0 &  \cdots &  0    &  \cdots &  0  \cr
\gen{1} &  0      &  1 &  \ddots   &\vdots&    &  \vdots \cr
~\vdots  &  \vdots &  \ddots   & \ddots &  0 &  \cdots &  \vdots\cr
\gen{1} &   0 & \cdots  &  0 &    1&\cdots &  0    
} 
$$
avec $\fa_1\supseteq \fa_2\supseteq \cdots\supseteq \fa_n $,
$\prod_{i=1}^n\fa_i=\fD_n(M)$,
 chaque $\fa_i$ de la forme  $\gen{a_i,\delta}$ et $\bigoplus_{k=1}^m\fa_k\simeq \gZ^m$.
\end{theorem}

\begin{proof} 
On reprend ce qui a été expliqué dans la \dem de la proposition \ref{propMpftorsionPrdim1}.
En examinant ce qui se passe dans le calcul de la réduction de Smith modulaire par \mnrsz, on obtient sur $\gZ$ une \egt et une congruence
$$\preskip.4em \postskip.4em 
LMC=M'\equiv\blocs{.7}{.9}{.7}{.0}{$D$}{$0$}{ }{ } \mod \delta$$
 avec
$D=\Diag(\ov{a_1},\dots,\ov{a_n})$ et $\gen{\ov{a_1}}\supseteq \cdots\supseteq \gen{\ov{a_n}}$. 

Considérons la \pmaz\footnote{On vérifie \imdt que c'est bien une \pmaz.}  $E:=\bordercmatrix [\lbrack\rbrack]{
       & \gen{1} &    \cdots&    \cdots&    \cdots & \gen{1}   \cr
\fa_1 &     &      &         \cr
~\vdots  & & \phantom{A}  &  M'   &  \phantom{A} &\cr
\fa_n &    &   &        
}$.
%(\gen{1},\dots,\gen{1};\fa_1,\dots,\fa_n;{M'})$. 
\\
Elle est de la forme
$$\preskip.3em \postskip.4em
E\;\;= 
\;\;\bordercmatrix [\lbrack\rbrack]{
       & \gen{1} & \gen{1}  & \cdots & \gen{1}   & \cdots & \gen{1}  \cr
\fa_1 &  a_1+\delta x_{11}  &  \delta x_{12} &  \cdots &  \delta x_{1n}  &\cdots&  \delta x_{1m}  \cr
\fa_2 &  \delta x_{21}      &  a_2+\delta x_{22} &    &&&  \vdots   \cr
~\vdots  &  \vdots &     &  \ddots &&&  \vdots  \cr
\fa_n &  \delta x_{n1} & \cdots  & \phantom{\ddots}  &    a_n+\delta x_{nn} &\cdots&  \delta x_{nm}   
}. 
$$
Oh merveille! c'est une \pma surjective.
En effet, 
$$\preskip.4em \postskip.4em \ndsp
\prod_{i=1}^n\fa_i=\fD_n(M)\hbox{  avec  }\fD_n(M)=\fD_n(LMC)=\fD_n(M')=\fD_n(\und E), 
$$
donc $\fD_n(E)=\gen{1}$.
D'après la proposition  \ref{thpmasurj}, on peut la compléter en une \pma \iv $B=(\fa_1,\dots,\fa_m;\gen{1},\dots,\gen{1};\und B)$ telle que la \pma $M'B^{-1}$ soit la matrice $M_1$ de l'énoncé.
Avec $C_1=CB^{-1}$, on obtient $LMC_1=M_1$, où $L$ et $C_1$ sont des \pmas de changements de base ($L$ et $M$ sont des matrices usuelles vues comme des \pmas ayant les \ids $\gen{1}$ en indices de lignes et en colonnes).
\end{proof}
%

%\subsection{Réduction de Smith des \pmas sur les \doksz}
%
%\hum{Je vais écrire un petit quelque chose}
%

%\newpage

\section*{Références}

\begin{description}
\item[\hbox{[ACMC]}] {\sc Lombardi H. \& Quitté C.}
 \emph{Alg\`ebre Commutative. Méthodes constructives.} 
 Calvage\&Mounet (2011). 
\item[\hbox{[CACM]}] English version of [ACMC]. Springer (2015).
\item[\hbox{[Cohen]}] {\sc Cohen H.}  \emph{Advanced topics in computational number theory}. Graduate texts in mathematics 193.  Springer-Verlag (1999).
\item[\hbox{[Modules]}] {\sc Díaz-Toca G.-M., Lombardi H. \& Quitté C.} \emph{Modules sur les anneaux commutatifs}. Calvage\&Mounet (2014).  
\item[\hbox{[MRR]}] {\sc Mines R., Richman F. \& Ruitenburg W.} {\it  A Course in
Constructive Algebra.} Universitext. Springer-Verlag, (1988).

\end{description}

\setcounter{tocdepth}{3}
 
\tableofcontents

\end {document}